\documentclass[10pt]{article}
\usepackage{amssymb,amsmath,amsfonts} 

\usepackage[colorlinks=true]{hyperref}
\usepackage{comment}

\renewcommand{\marginpar}[2][]{} %attivarlo per eliminare i commenti a margine

\newtheorem{theorem}{Theorem}[section]

\newtheorem{lemma}{Lemma}[section]
\newtheorem{proposition}{Proposition}[section]

\newtheorem{remark}{Remark}[section]

\newcommand{\diam}{{\rm diam}}
\newcommand{\io}{{\infty}}

\newcommand{\real}{ {\mathbb R}   }
\newcommand{\torus}{ {\mathbb T}   }
\newcommand{\integer}{ {\mathbb Z}   }
\newcommand{\complex}{ {\mathbb C}   }
\newcommand{\bB}{ {\mathbb B}   }
\newcommand{\cA}{ {\mathcal A}   }
\newcommand{\cB}{ {\mathcal B}   }
\newcommand{\cE}{ {\mathtt E}   }
\newcommand{\cP}{ {\mathcal P}   }

\newcommand{\cH}{ {\mathcal H}   }

\newcommand{\cD}{ {\mathcal D}   }

\renewcommand{\Im}{\, {\rm Im}\,}

\newcommand{\eproof}{\qed}%{\hfill $\Box$}
\newcommand\beq[1]{ \begin{equation}\label{#1} }
\newcommand{\eeq}{ \end{equation} }
\newcommand{\beqno}{ \[ }
\newcommand{\eeqno}{ \] }
\newcommand\beqa[1]{ \begin{eqnarray} \label{#1}}
\newcommand{\eeqa}{ \end{eqnarray} }
\newcommand{\beqano}{ \begin{eqnarray*} }
\newcommand{\eeqano}{ \end{eqnarray*} }

\newtheorem{definition}{Definition}[section]
\newcommand\dfn[1]{ \begin{definition}\label{#1} \rm}
\newcommand\edfn{ \end{definition} }
\newcommand{\proof}{\par\medskip\noindent{\bf Proof\ }}
\newcommand\equ[1]{{\rm (\ref{#1})}}
\newcommand{\nl}{{\smallskip\noindent}}
\newcommand{\giu}{{\medskip\noindent}}
\newcommand{\Giu}{{\bigskip\noindent}}
\newcommand{\noi}{{\noindent}}
\newcommand{\qed}{\hskip.5truecm
\vrule width 1.7truemm height 3.5truemm depth 0.truemm
\par\Giu}
\newcommand{\qedeq}{\hskip.5truecm
\vrule width 1.7truemm height 3.5truemm depth 0.truemm}
 \newcommand\casialt[3]{ \left\{  \begin{array}{ll}
 {#1} & \mbox{ {\rm if} ${#2}$} \\
 {#3} & \mbox{ {\rm otherwise}}
 \end{array} \right.}

\newcommand{\x}{\xi}

\newcommand{\e}{\varepsilon}

\renewcommand{\a }{\alpha }
\renewcommand{\b }{\beta }
\newcommand{\s }{\sigma }
\newcommand{\ii }{{\rm i} }
\renewcommand{\d }{\delta }

\newcommand{\g }{\gamma}
\newcommand{\f }{\varphi}

\renewcommand{\l }{\lambda }
\renewcommand{\L }{\Lambda }
\newcommand{\m }{\mu }

\renewcommand{\t }{\tau }
\renewcommand{\o }{\omega }
\renewcommand{\O }{\Omega }
\newcommand{\C}{\mathbb{C}}
\newcommand{\Z}{\mathbb{Z}}
\newcommand{\z }{\zeta }

\newcommand{\tetta}{\vartheta }
\newcommand{\checco}{\eta }
\newcommand{\K}{{\mathtt K}}

\newcommand{\Hpend}{H_{\rm pend}}

\newcommand{\ttM}{\mathtt M}
\newcommand{\ttL}{\mathtt L}
\newcommand{\ttH}{\mathtt H}
\newcommand{\ttD}{\mathtt D}
\newcommand{\ttm}{\mathtt m}
\newcommand{\cgot}{\mathfrak c}
\newcommand{\pgot}{\mathfrak p}

\newcommand{\Cgot}{\mathfrak C}
\newcommand{\Pgot}{\mathfrak P}
\newcommand{\Fgot}{\mathfrak F}
\newcommand{\Ggot}{\mathfrak G}
\def\R{\mathbb R}
\def\T{\mathbb T}

\def\muci{\kappa}  % esponente del logaritmo

\def\const{{\, c\, }}
\def\dst{\displaystyle}
\def\bks{\, \backslash\, }
\def\meas{{\rm\, meas\, }}
\def\Tp{{T_K^\perp}}
\def\TNp{T_N^\perp}
\def\vain{\ \rightsquigarrow\ }
\newcommand\eqby[1]{\stackrel{\equ{#1}}{=}}
\newcommand\leby[1]{\stackrel{\equ{#1}}{\le}}
\newcommand\proiezione{\, {\mathtt p}}
\newcommand\riscala{\varsigma_k}

\title{
KAM Theory for secondary tori  \\
\ \\
}

\begin{document}

\author{ 
%Nota di Luca Biasco e Luigi Chierchia \\{\bf Scientific chapter:} {\sl Mathematical analysis.} \\
\footnotesize L. Biasco  \& L. Chierchia
\\ \footnotesize Dipartimento di Matematica e Fisica
\\ \footnotesize Universit\`a degli Studi Roma Tre
\\ \footnotesize Largo San L. Murialdo 1 - 00146 Roma, Italy
\\ {\footnotesize biasco@mat.uniroma3.it, luigi@mat.uniroma3.it}
\\ 
}

\maketitle

{\small 
\tableofcontents
}

%\date{}

%\newpage

\begin{abstract}

\noindent
\begin{itemize}\item[(i)] 
In \cite{BClin} (Rend. Lincei Mat. Appl. {\bf 26} (2015), 1--10; see also arXiv:1503.08145 [math.DS]) the following result has been announced:

\nl
{\bf Theorem} {\sl
Consider a real--analytic nearly--integrable mechanical system with potential $f$, namely,  
a Hamiltonian system with  real-analytic Hamiltonian  
$$H(y,x)=\frac12 \sum_{i=1}^n y_i^2  +\e f(x)\ ,$$  
$(y,x)\in\real^n\times\torus^n$ being standard action--angle variables.
For ``general non--degenerate'' potentials $f$'s there exists $\e_0,a>0$ such that, if $0<\e<\e_0$, then the Liouville measure of the complementary of $H$--invariant  tori is smaller than $\e|\log \e|^a$.}

\nl
In this paper we provide a proof of such result.

\item[(ii)] The class  of ``general non--degenerate'' potentials $\cP_s$ (defined in  \S 1) is, for any given $s>0$,  an open and dense subset of real--analytic functions on $\torus^n$ having holomorphic extensions on $\{x\in\complex^n|\, |\Im x_i|<s\}$,  the topology being that induced by the weighted Fourier norm $|f|_s:=\sup_{k\in \integer^n} |f_k| e^{|k|s}$. 
The class $\cP_s$ is  also of ``full measure'' in  natural ways (compare Proposition~\ref{grizman} proven in Appendix~\ref{pogba}). 

\item[(iii)] The above Theorem is based on an extension of KAM Theory to a suitable $\e$--dependent neighborhood
of simple resonances $\{y\in B: y\cdot k=0\}$, with $|k|\le K\sim |\log \e|^b$, for a suitable $b>0$ and any given ball $B\subset \real^n$. The main issue is giving a quantitative analytic description of the integrability structure of the averaged Hamiltonian at simple resonances  suitable for application of KAM methods. Such analytic  properties are summarized in the ``Structure Theorem'' of \S 4, whose proof occupies the main part of this paper, namely, Sect's 5, 6 and Appendix~C, where action--angle variables for generic parameter--depending systems are discussed. 

\item[(iv)] In view of the Structure Theorem one can then apply simultaneously (for  $|k|\le K$) explicit classical KAM measure estimates (as given, e.g., in \cite{BCKAM}) and conclude the proof of the  main Theorem.

\end{itemize}

\end{abstract}

\section{Functional setting and Main Theorem}

%, probability measures and Fourier--projections} 

\nl 
In this paper we consider real--analytic functions, which are ``non--degenerate'' in a suitable sense\footnote{It would be easier to consider larger function spaces of smooth functions. However, the natural (both from the theoretical and applicative point of view) and most challenging setting is, we believe, that of real--analytic potentials.}. 

\nl
Let $s>0$ and 
%let $\torus^n_s$  the complex neighborhood of $\torus^n$ given by $\{\f\in \complex^n/(2\pi\real)^n$ such that $|\Im \f_j|\le s\}$ and 
consider the real--analytic functions on $\torus^n$ having zero average and finite ``sup--Fourier norm''
\begin{equation}\label{enorme}
|f|_s:=\sup_{k\in \integer^n} |f_k| e^{|k|s}<\infty\ ,
\end{equation}
where $f_k$ denotes Fourier coefficients and, as usual, $|k|$, for integer vectors, denotes the 1-norm $\sum |k_j|$. Denote by $ \cA_s^n$ the Banach space of such functions. 

\nl
Let $\integer^n_\sharp$ denote the set of integer vectors $k\neq 0$ in $\integer^n$ such that the  first non--null  component is positive:
\beq{coscia}
 \integer^n_\sharp:=
 \big\{ k\in\integer^n:\ k\neq 0\ {\rm and} \ k_j>0\ {\rm where}\ j=\min\{i: k_i\neq 0\}\big\}\ ,
 \eeq
and denote by ${\integer}^n_*$ {\sl the generators of one--dimensional maximal lattices}, namely, the set of  vectors $k\in  \integer^n_\sharp$ such that the greater common divisor of their components is 1, namely
\beq{cippalippa}
\integer^n_*:=\{k\in\integer^n_\sharp:\ {\rm gcd} (k_1,\ldots,k_n)=1\}\ ;
\eeq
then, the list of one--dimensional maximal lattices is given by the sets $\integer k$ with $k\in \integer^n_*$.

\nl
We can, now, decompose the Fourier expansion of
any function $f\in  \cA_s^n$ as sum of real--analytic functions of one--variable, which are the projection of $f$ onto the one--dimensional maximal lattice $\integer k$ (with $k\in \integer^n_*$), as follows:
%\footnote{Beware of notation: $\xi\to F^k(t)$ is a periodic function of one variable, whose Fourier coefficients, for $j\in \integer$, are given by $F^k_j=f_{jk}$.}:
\beq{dec}
f(x)= \sum_{k\in \integer^n_*} F^k(k\cdot x)\ ,\quad {\rm where}\quad
F^k(t):=\sum_{j\in \integer\backslash\{0\}} f_{jk}e^{ij t}
\eeq
$f_{jk}$ being the Fourier coefficient of $f$ with Fourier index $jk\in\integer^n$. 
Notice that, since $f\in  \cA_s^n$, the functions $F^k$ belong to $ \cA_{|k|s}^1$. 

\begin{definition} {\bf (The class $\cP_s$ of non--degenerate potentials) } \label{albertone}
Let $s>0,0<\d\leq 1$ and let   
\beq{miguel}
K_s(\d):=\const \max\big\{ 1\ ,\ \frac{1}{s}\ ,\ \frac{1}{s}\ \log\frac{1}{s \, \d} \big\}	\,,
\eeq
where $\const>1$ is
 a suitably large constant to be chosen 
 below (see \eqref{larocca}), depending only on $n$.

\nl 
 Let $\cP_s(\d)$  be the set of functions in $ \cA_s^n$ 
%there exists $ 0<\b<1$  such that 
%for every $k\in\integer^n_*$
such that,  for $k\in\integer^n_*$, the following holds

\giu
{\rm (P1)}
$\displaystyle{  |f_k|\geq \d |k|^{-\frac{n+3}2}\ e^{-|k|s} \ \ {\sl if}\ \  
|k|> K_s(\d)\,;}$

\giu
{\rm (P2)}
$\displaystyle{\  
%\min_{  |k|\leq K_s(\d) }\ \ 
\min_{\xi\in \real} \ \big( |\partial_\xi F^k(t)|+|\partial^2_\xi F^k(t)|\big) >0 \ \ {\sl if}\ \  
 |k|\le  K_s(\d)\,;}$

\giu
{\rm (P3)}   $F^k(t_1)\neq F^k(t_2)$
for every $0\leq t_1<t_2<2\pi$ such that
$\partial_t F^k(t_1)=\partial_t F^k(t_2)=0$ and
$|k|\leq  K_s(\d)$.
 
\giu
Finally, $\displaystyle\cP_s:=\bigcup_{\d>0} \cP_s(\d)$.
\end{definition}

\giu
An example of function $f\in {\cal P}_s(\d)$, as it is immediate to verify, is given by
\beq{esempietto}
f(x)=2\d \sum_{k\in \integer^n_*} e^{-|k|s}\,  \cos (k\cdot x)\ ,
\qquad{\rm i.e., }\qquad
f_k=\casialt{\dst\d e^{-|k|s}}{\pm k\in\integer^n_*}{0}\ .
\eeq

\nl
The class $\cP_s$ is ``general'' in several ways:  from  a probabilistic, topological and measure theoretical points of view.

\nl
To describe the probabilistic point of view, let us 
denote by
$\ell_\io^n$  the Banach space of complex sequences $z=\{z_k\}_{k\in \integer^n_\sharp  }$ with finite sup--norm 
$|z|_\io:=\sup_{k\in\integer^n_\sharp } |z_k|$.  
The map 
\begin{equation}\label{miserere}
j:f\in  \cA_s^n \to \big\{f_k e^{|k|s}\big\}_{k\in \integer^n_\sharp}\in \ell_\io^n
\end{equation}
is an isomorphism of Banach spaces\footnote{Recall that since the functions in $\cA^n_s$ are {\sl real}--analytic one has the reality condition $f_k=\bar f_{-k}$.},
which allows to identify functions in $ \cA_s^n$ with points in $\ell_\io^n$ and the Borellians  of $ \cA_s^n$ with those of $\ell_\io^n$. 
Now, consider the probability measure given by the standard normalized Lebesgue--product measure 
on the unit closed ball of $\ell_\io^n$, namely,   the unique probability measure $\mu$ on the Borellians of $\{z\in \ell_\io^n: |z|_\io\le 1\}$ such that,
given   Lebesgue measurable  sets $A_k$ in the unit complex disk $A_k\subseteq D:=\{w\in\complex:\ |w|\le 1\}$
with $A_k\neq D$ only for finitely many $k$, one has
 $$
\mu \Big(\prod_{k\in\Z^n_\sharp} A_k\Big)=
\prod_{\{k\in \integer^n_\sharp:\, A_k\neq D\}} 
{\rm meas}(A_k)\,
$$
where ``meas'' denotes  the  normalized Lebesgue measure on the unit complex disk $D$. Denote by $\bB$ the closed ball of radius one in $\cB_s^n$ and by $\cB$ the Borellians in $\bB$.
Then,  {\sl the isometry $j$ in \eqref{miserere}
naturally induces a  measure
$\mu_s$ on the Borellians $\cB$}.

\nl 
The properties of $\cP_s$ are collected in the following proposition, whose simple proof is given in Appendix~\ref{pogba}.

\begin{proposition}\label{grizman}
Let $s>0$. The set $\cP_s\subseteq \cA_s^n$ contains an open dense set, is prevalent, $\cP_s\, \cap \, \bB\in \cB$ and $\m_s(\cP_s\cap \bB)=1$.
\end{proposition}

\giu
Let $\|\cdot\|$ be the standard Euclidean norm on $\R^n.$ Then, one has the following 

\begin{theorem}\label{bounty}
Let $s>0$ and  let $\O$ be a bounded region in $\real^n$ with $n\ge 2$.
Let $f\in \cP_s$ and consider the Hamiltonian
\beq{storia}
H:= \frac12 \|y\|^2 +\e f(x)\ ,
\eeq
There exist $\e_0>0$ and $\muci>0$ such that, for any $0<\e<\e_0$,  the measure of the set of $H$--trajectories in $\O\times \torus^n$, which do not lie on an invariant Lagrangian (Diophantine) torus, is bounded by  $\e |\log \e|^\muci$.
\end{theorem}
We remark that the constants $\e_0$ and $\muci$ depend only on  
$n$, $s$, and $F^k$ with $k\in \integer^n_*$, $|k|\le K_s(\d)$.

\begin{remark}\label{abbada}
In proving Theorem \ref{bounty} we will assume that
\begin{equation}\label{bada}
\O=B_1(0):=\{ \|y\| <1\}\qquad\text{and}\qquad   |f|_s=1\,.
\end{equation}
This is not restrictive since we can always consider
a large enough ball $B_R(0)\supseteq \O$
and rescale the action and the time in order to obtain \eqref{bada}
(suitably renaming $\e$ and $f$). 
\end{remark}

\section{Geometry of resonances}

In this section we construct a covering of  $\real^n$ (thought of as frequency space) by three regions: a non--resonant region
$\Omega_0$, a neighborhood of simple resonances $\Omega_1$ and a region $\Omega_2$ of ``small'' measure   containing all other resonances.  The sets $\O^1$ and $\O^2$ will be described in terms of linear maps $L_k$, $k\in\integer^n$,   that depend on a given resonance $\{y\cdot k=0\}$: such maps $L_k$ will later be associated to generating functions $S(J,x)= x\cdot L_k J$, whose corresponding symplectic maps  have the role of ``straighten out'' the geometry.

\nl
Fix $k\in \integer^n\backslash\{0\}$ with gcd$(k_1,\ldots,k_n)=1$. 
 Then, there exists a matrix 
$A_k\in\ {\rm Mat}_{n\times n}(\Z)$ such that\footnote{Here, $k$ is  a row vector. Normally we do not distinguish between row and column vectors since it will be clear from context. The notation  $|M|_\infty$, with $M$ matrix or vector, denotes the maximum norm $\max_{ij}|M_{ij}|$ or, respectively, $\max_i |M_i|$. }
\begin{equation}\label{scimmia}
A_k=\binom{\hat A_k}{k}\,,\ \ \ \hat A_k=
\hat A_k\in{\rm Mat}_{(n-1)\times n}(\Z)\,,
\ \ \
\det A_k=1\,,\ \ \ 
|\hat A_k|_\infty\leq |k|_\infty\ ,
\end{equation}
and\footnote{\noindent
Recalling that 
for  any $n\times n$ matrix $M$, one always  has 
$|\det M|\leq n^{n/2} |M|_\infty^n$, \equ{atlantide} follows by the D'Alembert expansion of determinants (with $c=(n-1)^{(n-1)/2}$).}
\begin{equation}\label{atlantide}
|A_k^{-1}|_\infty\leq 
%(n-1)^{(n-1)/2}
\const |k|_\infty^{n-1}\,;
\end{equation}
the existence of such a matrix is guaranteed by an elementary result of linear algebra based on Bezout's Lemma (see Lemma~\ref{alfonso} in appendix~\ref{appendicedialfonso}).

\nl
We then define a linear map\footnote{Without further notice, we shall always identify linear maps with the associated matrices.}  $L_k:\R^n\to\R^n$  by setting
%, for all  $J=(\hat J,J_n)\in \real^{n-1}\times \real$, 
\begin{equation}\label{LPD}
L_k:J=(\hat J,J_n)\in\real^{n-1}\times \real  \mapsto
L_k J:=
J_n k 
+ \proiezione^\perp_k \hat A^T_k\hat J   \,,
\end{equation}
where 
$\proiezione^\perp_k$
is the orthogonal
projection  on the subspace perpendicular  to\footnote{Explicitely,
for  $y\in\R^n$,
\begin{equation}\label{bardo}
\proiezione^\perp_k y:= y- \frac{1}{\kappa}(y\cdot k) k\,.
\end{equation}
}
$k$. Observe that $L_k$ can be also written as composition of two linear maps:
\beq{LPDbis}
L_k= A_k^T   U_k 
\eeq
where  $U_k $ acts as the identity on the first $(n-1)$ components and:
\begin{equation}\label{centocelle}
U_k : (\hat J,J_n)\mapsto \big(\hat J,J_n-
\frac1\kappa (\hat A_k k)  \cdot \hat J\big)
\,, \qquad {\rm i.e.} \quad
U_k =\left(\begin{matrix}
{\rm I}_{n-1} & 0\cr -\kappa^{-1} \hat A_k k & 1 \cr
\end{matrix}\right)
\ ,
\end{equation}
where ${\rm I}_{m}$ denotes the $(m\times m)$--identity matrix and where 
\beq{cappetto}
\kappa:=\|k\|^2=k_1^2+\cdots+k_n^2\in {\mathbb N}\ .
\eeq 
Note that, by \eqref{scimmia}, we 
get\footnote{Here and in the following, we shall denote by ``$\const$'' suitable constants (which, in general, differ from formula to formula) depending only on $n$.} 
\beq{verruca}
\| U_k\|\ , \|U_k^{-1}\|\leq \const\ .
\eeq
Elementary properties of $L_k$ are the following\footnote{Eq. \equ{diego} follows from \equ{scimmia}, \equ{LPDbis} and \equ{centocelle}.
Identities \equ{cicogna},
\equ{farfalla}, \equ{raffinatoeffendi}  and 
the first bound in \equ{LPDter}
follow directly from  definition \equ{LPD}; the bound on $\|L_k^{-1}\|$ follows by observing that $\|L_k^{-1}\|=\|U_k^{-1}A_k^{-T}\|\le \|U_k^{-1}\|\, \|A_k^{-T}\|$ and that $\|U_k^{-1}\|\le \const $ and that
$\|A_k^{-T}\|\le \const \|k\|^{n-1}$ (as it follows by bounding the norm of a matrix by a constant times the maximum of its entries
and using  the co--factor representation  for the inverse of $A_k$, and taking into account \equ{scimmia}). 
}:
\beqa{diego}
&& \det L_k=1\ ,\\
&&  L_kJ\cdot k =\kappa J_n \ , \label{cicogna}
%\qquad {\rm and}\qquad (L^T x)_n= \frac{k\cdot x}{\|k\|^2}\ ,\ \ \ serve???? 
\\
&&\nonumber\\
&& \|L_k J\|^2= \kappa J_n^2 + 
\| \proiezione_k^\perp \hat A_k \hat J\|^2  
\label{farfalla}\ ,\\
&& \hat \O\subseteq \real^{n-1}\ ,\ a>0\  \implies\ \ 
L_k\big(\hat \O\times
 (-\frac{a}{\kappa},\frac{a}{\kappa})\big)= \{|y\cdot k|< a\}
\cap  \big\{y: \proiezione_k^\perp y\in \proiezione_k^\perp  \hat A_k^T \hat \O\big\} \ ,
\label{raffinatoeffendi}
\\
&&\nonumber\\
&& \|L_k \|\le \const \|k\| \ ,\qquad {\rm and}\qquad \|L_k^{-1} \|\le \const \|k\|^{n-1} \ .
\label{LPDter}
\eeqa
Further more interesting properties of $L_k$ are given in the following simple
\begin{lemma}
\label{zioeffendi}
{\rm (i)} The map $\proiezione^\perp_k \hat A_k:\real^{n-1}\to k^\perp$ is a linear isomorphism.

\noindent
{\rm (ii)} For any $a>0$, 
$L_k\big(\real^{n-1}\times (-\frac{a}{\kappa},\frac{a}{\kappa})\big)=\{|y\cdot k|<a\}$. 
\end{lemma}

\proof
(i): Let $\hat A_k=\left( \begin{matrix} a^1 \\ a^2\\ \vdots \\a^{n-1} \end{matrix}\right)$,  with $a^i\in \integer^n$. Then $\proiezione_k^\perp \hat A_k \ \real^{n-1}= {\rm span} \{\proiezione_k^\perp a^1,...,\proiezione_k^\perp a^{n-1}\}$. But the vectors $\proiezione_k^\perp a^i$ are linearly independent  (if $0= \sum \l_i \proiezione_k^\perp a^i = 
\proiezione_k^\perp (\sum \l_i \ a^i)$, then there exists $c$ such that $\sum \l_i \ a^i=c k$, which implies that $\l_1=\cdots=\l_{n-1}=c=0$ since $A_k={\hat A_k \choose k}$ has determinant one and hence the vectors $a^1,...,a^{n-1},k$ are linearly independent). The claim then follows from the rank--nullity theorem of linear algebra.

\nl
(ii): Let $W:=L_k(\real^{n-1}\times (-a/\kappa,a/\kappa))$. From \equ{cicogna} it follows that $W\subseteq \{y|\ |y\cdot k|<a\}$.
Now, let $y\in \real^n$ be such that $|y\cdot k|<a$ and define $J_n:= (y\cdot k)/\kappa$. Then, $|J_n|<a/\kappa$ and, furthermmore, $y- J_n k\in k^\perp$. Thus, by part (i) of this Lemma, there exists $\hat J\in \real^{n-1}$ such that $y-J_n k= \proiezione_k^\perp \hat A_k \hat J$, hence, $y=L_k (\hat J,J_n)$,  proving that 
$ \{y|\ |y\cdot k|<a\} \subseteq W$ and, thus, $ \{y|\ |y\cdot k|<a\} = W$.
\qed

\nl
Set\footnote{Recall  \eqref{cippalippa}.} 
\begin{equation}\label{cippalippaK} 
\Z^n_{*,K}:=\{ k\in\Z^n_*\,, \ \ |k|\leq K\}\,.
\end{equation}
\nl
To quantify neighborhoods of resonances, as standard, we introduce two $\e$--dependent  Fourier cut-offs $K,\K$ and an $\e$--dependent ``width'' $\a$  of simple resonance, by setting: 
\begin{equation}\label{lapparenza}
K:= \log^2 \frac{1}{\e}\,,
\qquad\qquad \K:=K^2
\,,\qquad \qquad 
\a:=\sqrt\e K^{\nu+1}\ ,
%\log^2 \frac{1}{\e}\ , NO\ FARE  \ K:=\log^{\nu_0+1} \frac{1}{\e} \,,\ \nu_0\geq 2, Forse\ non\ serve\ farlo!  LUCA
\end{equation}
where 
\beq{gnu}
\nu>n+1
\eeq
is a suitable constant to be fixed later.
Let us assume that $\e$ is so small that
(recall \eqref{LPDter})
\begin{equation}
\|L_k^{-1} \|\le  K^n \ .
\label{LPDquater}
\end{equation}
Set
\begin{equation}\label{arrosticini}
\hat D:=B_{K^n} (0)
\end{equation}

\nl
For $k\in \Z^n_{*,K}$,  define 
%RICONTROLLARE \eqref{shevket} LA  E IL LEMMA \ref{cb500f}  LUCA
\begin{eqnarray}
&&\hat Z_k :=
\left\{ \hat J\in\hat D \  : \ 
\min_{l\in \Z^n_{*,\K}\,,\   l \notin \Z k}
\big| \big(\proiezione_k^\perp \hat A_k^T \hat J\big)\cdot l\big|
\geq 3\a\K \frac{\|l\|}{\|k\|}   \right\}
\,,
\ \ 
Z_k^\sharp :=
\hat Z_k\times(-\frac{\a}{2\kappa},\frac{\a}{2\kappa})\subseteq \R^n\,,\quad\quad
\label{shevket}
\\
\ \nonumber\\
&&Z_k :=
\hat Z_k\times(-\frac{\a}{\kappa},\frac{\a}{\kappa})\subseteq \R^n\,,
\qquad {\rm and} \qquad  
Z_k' :=
\big( \hat D\setminus \hat Z_k\big)\times(-\frac{\a}{\kappa},\frac{\a}{\kappa})\subseteq \R^n\,,
\label{caspio}
\\
\ \nonumber\\
&&\O^0:=\{ \|y\|<1\, :\, \min_{k\in \Z^n_{*,K}}|y\cdot k|\geq \a/2  \}
\,,\ \ \ 
\O^1:=\!\!\bigcup_{k\in\Z^n_{*,K}}\!\! L_k Z_k^\sharp \,,\ \ \
\O^2:=\!\!\bigcup_{k\in \Z^n_{*,K}}\!\! L_k Z_k'\,.
\label{neva}
\end{eqnarray}
\begin{remark}\label{orhan}
The set $\O^0$ is a non--resonant set. The set $\O^1$, by \equ{raffinatoeffendi} and Lemma~\ref{zioeffendi}, is seen to be a suitable neighborhood of simple resonances  $\{y\cdot k= 0\}$ with $k\in \Z^n_{*,K}$; finally, $\O^2$ is a neighborhood of order--two or higher resonances. Next Lemma clarifies and quantifies these observations. 
\end{remark} 

\begin{proposition}\label{cb500f}
{\rm (i)} \ \ \ 
 $\O^0\cup\O^1\cup \O^2\supseteq B_1(0)$.

\nl
{\rm (ii)}
The set $\O^0$ is $(\a/2,K)$ completely non--resonant\footnote{\label{amarcord}We follow the terminology introduced in \cite{poschel}: given $\a,K>0$ and a sublattice $\L\subseteq \integer^n$, one says that  $D\subseteq  \real^n$ (or $D\subseteq \complex^n$) is ``$(\a,K)$ non--resonant modulo $\L$'' if $|y\cdot k|\ge \a$ for all $y\in D$ and $k\in \integer^n\bks \L$ with $|k|\le K$; if
$\L=\{0\}$ is the trivial lattice, then $D$ is said to be $(\a, K)$ completely non--resonant.}, i.e., 
\begin{equation}\label{cipollotto}
y\in \O^0 \qquad\Longrightarrow\qquad
|y\cdot k|\geq \a/2\,,\ \ \ \forall\, 0<|k|\leq K\,.
\end{equation}
{\rm (iii)} For each $k\in\Z^n_{*,K}$, the set $L_k Z_k$  is $(2\a\K/\|k\|,\K)$ non--resonant modulo $\integer k$, i.e., 
\begin{equation}\label{cipollotto2}
y\in L_k Z_k \qquad\Longrightarrow\qquad
|y\cdot l |\ge 2\a\K/\|k\|\,,
\ \ \ \forall\, \ l\in \integer^n\ ,\ l\notin \Z k\ , \    |l|\leq \K\ .
\end{equation}
{\rm (iv)}  There exists a constant $\const>0$
depending only on $n$ such that:
\begin{equation}\label{teheran4} 
\meas (\O^2) \le \const \a^2  K^{n^2-n-1}\K^{n+2}
\ .
\end{equation}
\end{proposition}
\nl
\begin{remark}
\label{cbf1000}
Since $\a$ has been chosen as $\sqrt\e K^\const$ (see \equ{lapparenza}), from \equ{teheran4} it follows that
the region of second or higher order resonances $\O^2$ has Lebesgue measure smaller than $(\const \e |\log \e|^\const)$ and therefore no further analysis on $\O^2$ is needed.
\end{remark}

\begin{remark}\label{laforcella}
In the definition  of $\hat Z_k$ in \eqref{shevket} one could also use a 
smaller set $\hat D_k\subset \hat D,$
such that Proposition \ref{cb500f} still holds true. We have chosen a unique $\hat 
D$ for every $k$, just for simplicity. 
\end{remark}

\begin{remark}
\label{spuntature}
The geometry of resonances here
is different from the geometry of resonances (in the convex case) as discussed, e.g.,  in  \cite{poschel}. In fact, in \cite{poschel}
more resonances are disregarded in the non--resonant set,
namely, the resonances with  $|k|\leq (1/\e)^a$, $a>0$.
Furthermore,  
the neighborhood of  simple resonances  in \cite{poschel} has width $\e^b,$ $0<b<1/2$, 
%which is larger than $R_k$.
and, as a consequence, the set of double resonances has measure greater than
$\e^{2b}$, which is a set not negligible for our purposes.
On the other hand, in  Nekhoroshev's
theorem one can average out the perturbation up to an exponentially
high order $e^{- {\rm const\, }(1/\e)^a},$ while we will get only $\e^{|\log\e|}$.
\end{remark}

\proof {\bf of Proposition \ref{cb500f}}.

\nl
(i):  If $y\notin \O^0$ with $\|y\|<1$,  there exists a $k\in \integer^n_{*,K}$ such that $|y\cdot k|<\a/2$; but then, in view of Lemma~\ref{zioeffendi}, $y$ belongs to
\begin{eqnarray*}
&&\{y'\in \real^n: |y'\cdot k|<\a/2\}\cap B_1(0)= L_k\big(\real^{n-1}\times (-\a/2\kappa,\a/2\kappa)\big)  \cap B_1(0)
\\
&&=L_k\bigg(\Big(\real^{n-1}\times (-\a/2\kappa,\a/2\kappa)\Big) \cap L_k^{-1}
 B_1(0)\bigg) 
\\
&&\stackrel{\eqref{LPDquater}}\subseteq
 L_k\bigg(\Big(\real^{n-1}\times (-\a/2\kappa,\a/2\kappa)\Big) \cap 
 B_{K^n}(0)\bigg)
 \\
&&\stackrel{\eqref{arrosticini}}\subseteq
 L_k\Big( \hat D\times (-\a/2\kappa,\a/2\kappa)\Big) 
\\
&&\subset
L_k(Z_k^\sharp\cup Z'_k)
\subset \O^1\cup \O^2\,.
\end{eqnarray*}

\nl
(ii): Let $y\in\O^0$ and  let $0<|k|\leq K$. Then,  there exist $j\in\Z\setminus\{ 0\}$
and $k'\in\Z^n_{*,K}$
with\footnote{Indeed, $j=\pm\, {\rm gcd}\{k_1,...,k_n\}$.} $k=jk'$, so that 
$$
|y\cdot k|=|j||y\cdot k'|\geq |y\cdot k'|
\stackrel{\eqref{neva}}
\geq \a/2\,,
$$
proving \eqref{cipollotto}.

\nl
(iii): Let   $y=L_kJ=L_k(\hat J,J_n)$ for some $k\in \integer^n_{*,K}$,  $\hat J\in \hat Z_k$ and $|J_n|<\a/\kappa$. Let, also,  $l\in\integer^n$, $l\notin \integer k$ with $|l|\le \K$. As above, there exists 
$j\in\Z\setminus\{ 0\}$ and 
$l'\in\Z^n_{*,\K}$ such that $l=jl'$. Then, 
%By the same argument used in (i), we see that \equ{shevket} implies that $\big| \big(\proiezione_k^\perp \hat A_k^T \hat J\big)\cdot l\big| \geq 3\a {K}{\|k\|}$ for all $l\notin \integer k$ with $|l| \le K$.
\beqano
|y\cdot l| &=&| L_kJ\cdot l| \stackrel{\equ{LPD}}= 
\Big| J_n \, k\cdot l + \big(\proiezione_k^\perp \hat A_k^T \hat J\big)\cdot l\Big|
\ge  \Big| \big(\proiezione_k^\perp \hat A_k^T \hat J\big)\cdot l\Big|- |k\cdot l|  |J_n|\\
&=&  |j|\, \Big| \big(\proiezione_k^\perp \hat A_k^T \hat J\big)\cdot l'\Big|- |k\cdot l|  |J_n|
\\
&\stackrel{\eqref{shevket}}\geq & 3\alpha\K \frac{\|l\|}{\|k\|} - \a\, \frac{\|l\|}{\|k\|} \ge 
2\alpha\K \frac{\|l\|}{\|k\|} 
\ge 2 \a\frac{\K}{\|k\|}\,,
\eeqano
proving \equ{cipollotto2}.

\nl
(iv): 
Then (denoting Lebesgue measure by ``$\meas$''), from the definition of $\O^2$ in \equ{neva} it follows:
\beqa{palle}
\meas (\O^2) &=& 
\meas\big(\bigcup_{k\in \integer^n_{*,K}} L_k Z'_k\big)
\leq  \sum_{k\in \integer^n_{*,K}} \meas
(L_k Z'_k)
\nonumber\\
&=&  \sum_{k\in \integer^n_{*,K}} |\det L_k| \meas(Z'_k)\stackrel{\equ{diego}}=
 \sum_{k\in \integer^n_{*,K}} \meas(Z'_k)\ .
\eeqa
Moreover\footnote{
Denoting Lebsgue measure on $\real^{n-1}$ again by ``$\meas$''.}
\beqa{sonolesettedivenerdi}
\meas(Z'_k) \leq \const \, \a\kappa^{-1}\,  \sum_{ l\in \integer^n_{*,\K} , l\notin \integer k} \meas
 \Big\{\|\hat J\|\leq K^n \ : \  \big| (\proiezione_k^\perp \hat A_k^T \hat J)\cdot l \big| < 3\a \K \frac{\|l\|}{\|k\|}\Big\}\,.
\eeqa
Now, denoting by 
\beq{v}
v_{k,l}:= \|k\|^2 l - (l\cdot k) k= \|k\|^2 \proiezione_k^\perp l\ ,
\eeq
we see that
\beq{setteecinque}
 | (\proiezione_k^\perp \hat A_k^T \hat J)\cdot l | < 3\a \K \frac{\|l\|}{\|k\|}
 \qquad \Longleftrightarrow \qquad
| \hat J\cdot \hat A_k v_{k,l}|<3\a \|l\|\,  \|k\| \K\ ,
\eeq
so that \equ{sonolesettedivenerdi} reads
\beq{setteequindici}
\meas(Z'_k) \leq \const \, \a\kappa^{-1}\,  \sum_{ l\in \integer^n_{*,\K} , l\notin \integer k} \meas
 \Big\{\|\hat J\|\leq K^n \ :\  | \hat J\cdot \hat A_k v_{k,l}|<3\a \|l\|\,  \|k\| \K\}\ .
\eeq
Now, observe that $v_{k,l}\in \integer^n\bks\{0\}$ (since $l\notin \integer k$) and that $v_{k,l}\in k^\perp$. But then 
$\hat A_k v_{k,l}\neq 0$ (indeed, $\hat A_k v_{k,l}=0$ implies that $A_k v_{k,l}= \binom{\hat A_k v_{k,l}}{k\cdot v_{k,l}}=0$, contradicting the invertibility of $A_k$), hence (since $\hat A_k v_{k,l}\in \integer^n\bks\{0\}$),
$\|\hat A_k v_{k,l}\|\ge 1$.
Thus, from \equ{setteequindici}  there follows\footnote{In general, fixed a positive integer $m$, there exists a constant $c>0$ such  that
for every  $w\in \real^m\bks\{0\}$ and $b>0$ one has: $\meas \{y\in \real^m: \|y\|\le r\ ,\ {\rm and}\  |y\cdot w|<b\}\le c r^{m-1} b/\|w\|$.}
\beq{venticello}
\meas(Z'_k\cap B_{R_k})\le  \const \, \a \,  \sum_{ l\in \integer^n_{*,\K} , l\notin \integer k} \a \K  \frac{\|l\|}{\|k\|} K^{n(n-2)}
\le \const \a^2 K^{n(n-2)} \K^{n+2} \|k\|^{-1}\ ,
\eeq
which, together with \equ{palle}, yields
\eqref{teheran4}.
\qed

\section{Normal Forms}\label{vongole}
In this section we describe a normal form lemma, which  allows to average out non--resonant Fourier modes of the perturbation on suitable non--resonant regions, and then apply it on $\O^0$ and $\O^1$. 

\nl
We remark that such normal form lemma is not standard  as, for technical reasons which will be clarified later, we shall need estimates in a complex domain {\sl very close} to the initial one.

\subsubsection*{Notation} 
\label{giallo}
Given a set $D\subseteq \R^m,$ $r>0$ 
we denote by $D_r\subseteq\C^m$  the complex  open neighborhood of $D$ formed by points  $z\in \complex^m$ such that $\|z-y\|<r$, for some $y\in D$.

\nl
Given $s>0,$ we denote by $\torus^n_s$ the open complex neighborhood of $\torus^n$ given by 
$$
\T^n_s:=\{ x\in\C^n\ \ :\ \  \max_{1\leq j\leq n}|{\rm Im}x_j|<s \}
/2\pi\Z^n
$$
Given a  real--analytic function
 $f:D_r\times \T^n_s\to\C$,
$f(y,x)=\sum_{k\in\Z^n}f_k(y) e^{\ii k\cdot x},$
we consider the weighted sup--norm 
\begin{equation}\label{laga}
|f|_{r,s}
%:= \sup_{y\in D_r} |f(y,\cdot)|_s
 := \sup_{k\in\Z^n}
 \big(
 \sup_{y\in D_r}|f_k(y)|e^{|k|s}
 \big)
 \, ;
\end{equation}
if the (real) domain need to be specified, we let:
\beq{romaostia}
|f|_{D,r,s}:=|f|_{r,s}\ .
\eeq

%\giu
%Given a sublattice $\L$ of $\Z^n$ and $\a>0, K\geq 1$ we say that 
%$D\subseteq \C^n$ is ($\a$,$K$) {\sl non--resonant modulo}
%$\L$ for the (integrable) Hamiltonian $h(y)$ if for every $\o=h'(y),$ $ y\in D$
%$$
%|\o\cdot k|\geq \a\qquad \forall\,  |k|\leq K\,,\ \ \  k\notin\L\,.
%$$ 

\nl
Given $f(y,x)=\sum_{k\in\Z^n}f_k(y)e^{\ii k\cdot x}$ and   a sublattice $\L$ of $\Z^n$, 
we denote by $\proiezione_\L$ the projection on the Fourier coefficients in $\L,$ namely
$$
\proiezione_\L f:=\sum_{k\in\L}f_k(y)e^{\ii k\cdot x}\,.
$$
and by $\proiezione_\L^\perp$ its ``orthogonal'' operator (projection on the Fourier modes in $\integer^n\bks\L$):
$$
\proiezione_\L^\perp f:=\sum_{k\notin\L}f_k(y)e^{\ii k\cdot x}\,.
$$
Finally, given $N>0$, we introduce   the following ``truncation'' and ``high--mode'' operators 
\begin{equation}\label{labestia}
T_N f:=\sum_{|k|\leq N}f_k(y)e^{\ii k\cdot x}\,,
\qquad
\TNp f:=\sum_{|k|>N}f_k(y)e^{\ii k\cdot x}
\,.
\end{equation}
For later use, we point the following elementary decay property of  analytic function with vanishing low modes: 
\begin{equation}\label{582}
T_N f=0\ ,\ 0<\s<s \qquad
\Longrightarrow\qquad
|f|_{r,s-\s}\leq e^{-N\s} |f|_{r,s}\,  .
\end{equation}
We are, now,  ready to state the Normal Form Lemma we need. In order not to introduce too many symbols we shall denote by $H=h+f$ (but {\sl without} $\e$) the Hamiltonian and by $\a$ and $K$ (which have been already fixed in \equ{lapparenza}) the non--resonance parameters, however the lemma applies to arbitrary $H$, $\a$ and $K$.

\begin{lemma}\label{pesce}{\bf (Normal Form Lemma)}
Let $r,s,\a$ be positive numbers, $K\ge 2$, $D\subseteq \real^n$,  and let $\L$ be a sublattice of $\integer^n$.
Let $H(y,x)=h(y)+f(y,x)$ be real--analytic on  
$ D_r \times \T^n_s$ with $|f|_{r,s}<\infty.$
Assume that $D_r$ is
($\a$,$K$) non--resonant modulo\footnote{In case $\L=\{0\}$, one also says that $D_r$ is completely $(\a,K)$ non--resonant.} $\L,$ namely
\begin{equation}\label{mentovare}
|h'(y)\cdot k|\geq \a\,,\qquad \forall\,
y\in D_r\,,\ k\notin\L\,, \ |k|\leq K
\end{equation}
 and that
\begin{equation}\label{gricia}
\tetta_* := \frac{2^9 n K^3}{\a r s}\, |f|_{r,s}  <1\,.
\end{equation}
Then,  there exists a real--analytic 
 symplectic change of variables
\begin{equation}\label{dadostar}
\Psi: D_{r_*}\times \T^n_{s_*} \to D_r \times \T^n_s \, \quad {\rm with} 
\qquad
r_*:=r/2\,,\ \ 
s_*:=s(1-1/K)
\end{equation}
such that 
\beq{senzanome}
H\circ\Psi=h+ f^\flat +
f_*\,,   \qquad f^\flat := \proiezione_\L f+ \Tp\proiezione_\L^\perp f
\eeq
with 
\begin{equation}\label{pirati}
%|\proiezione_\L f_*|_{r_*,s_*}\,,\  |\Tp\proiezione_\L^\perp f_*|_{r_*,s_*}
|f_*|_{r_*,s_*}
\leq 2\tetta_* |f|_{r,s}\,,\qquad 
|T_{K}\proiezione_\L^\perp f_*|_{r_*,s_*} 
\leq
(\tetta_* /2)^K |f|_{r,s}\,.
\end{equation}
\end{lemma}
\proof\!: See Appendix~\ref{provapesce}.
\begin{remark}\label{cb500x} 
{\rm (i)} Having information on non--resonant Fourier modes up to order $K$, the best one can do is to average out
the non--resonant Fourier modes up to order $K$, namely, to ``kill" the term $T_K \proiezione_\L^\perp f$ of the Fourier expansion of the perturbation. This explains the ``flat'' term $f^\flat=\proiezione_\L f+ \Tp\proiezione_\L^\perp  f$ surviving in \equ{senzanome} and  which cannot be removed in general.
Now, think of the remainder term $f_*$ as 
$$f_*=\proiezione_\L f_* + \big( \Tp\proiezione_\L^\perp  f_* + T_{K}\proiezione_\L^\perp  f_*\big)\ ;$$
then,  $\proiezione_\L f_*$  is a $(\tetta_* |f|_{r,s})$--perturbation of the part in normal form (i.e., with Fourier modes in $\L$), while 
$ \Tp\proiezione_\L^\perp  f_*$ is, by \equ{582}, a term exponentially small with $K$ (see also below)
and 
$T_{K}\proiezione_\L^\perp  f_*$ is a very small remainder bounded by $(\tetta_* /2)^K |f|_{r,s}$.

\nl
{\rm (ii)}
The ``novelty'' of this lemma is that the bounds in \equ{pirati} hold on the large angle domain $\torus^n_{s_*}$ with $s_*=s(1-1/K)$. In particular, it will be important in our analysis 
(precisely in order to obtain \eqref{inbloom} below)  the first estimate in \equ{pirati}. The drawback of the gain in angle--analyticity strip is that the power of $K$ in the smallness condition 
\equ{gricia} is not optimal: for example in \cite{poschel} the power of $K$ is one (and $s_*=s/6$). 

\nl
{\rm (iii)} To compare with more standard formulations, such as the Normal Form Lemma in \S~2 of \cite{poschel}, write \equ{senzanome} as 
\beq{silvia}
H\circ \Psi=h+g+ f_{**}\qquad {\rm with}\quad \proiezione_\L g=g\ ,\quad \proiezione_\L f_{**}=0\ .
\eeq
Then,  $g=\proiezione_\L f + \proiezione_\L f_*$,
$f_{**}=\Tp\proiezione_\L^\perp f+\proiezione_\L^\perp f_*= T_K\proiezione_\L^\perp f_*+ \Tp\proiezione_\L^\perp (f_*+f)$ and the following bounds hold\beq{giovanni}
|g- \proiezione_\L f|_{r_*,s_*}\le 2\tetta_* |f|_{r,s}\ ,\qquad
|f_{**}|_{r_*,s/2}\le 2 e^{-Ks/2} |f|_{r,s}\ ,
\eeq
provided
\beq{francescone}
\tetta_* \le e^{-s}/2 \ ,\qquad \qquad
K\ge 2
\eeq
(which will be henceforth assumed).
To check \equ{giovanni}, notice that by  \equ{pirati} and 
\equ{582} (used with $N=K$, $s$ replaced by $s_*$ and $\s=\frac{s}2-\frac{s}K$ so that $s_*-\s=s/2$ and $e^{-K\s}=e^{-Ks/2} \cdot e^s$), one gets
\beqano 
|f_{**}|_{r_*,s/2} &\le & |T_K\proiezione_\L^\perp f_*|_{r_*,s_*} + |\Tp\proiezione_\L^\perp (f_*+f)|_{r_*,s_*-\s}\\
&\le& |T_K\proiezione_\L^\perp f_*|_{r_*,s_*} + |\Tp f_*|_{r_*,s_*-\s}
+|\Tp f|_{r_*,s-s/2}
\\
&\le& \tetta_*^K |f|_{r,s} + e^{-Ks/2}   
(e^s\tetta_* +1) |f|_{r,s} 
\\
&\leby{francescone} & 2 e^{-Ks/2} |f|_{r,s}\ .
\eeqano 
%
%STIMA LUCA
%\beqano 
%|f_{**}|_{r_*,s/2} &\le & |T_K\proiezione_\L^\perp f_*|_{r_*,s_*} + |\Tp f_*|_{r_*,s_*-\s}
%+|\Tp f|_{r_*,s/2}\\
%&\le& \tetta ^K |f|_{r,s} + e^{-Ks/2}   (
%e^s\tetta +1) |f|_{r,s} \le 3 e^{-Ks/2} |f|_{r,s}\ ,
%\eeqano 
%provided, e.g.,  $\tetta \le e^{-s}$ and $K>2$.

\nl
{\rm (iv)} In our applications $\a$ and $K$ (or $\K$) are as in \equ{lapparenza},
$r \stackrel{\sim}{>} \a/K$
and $f$ is replaced by $\e f$. Thus, 
\beq{zanzara}
\tetta_* \sim  |\log \e|^{-4(\nu-1)}\qquad \stackrel{\equ{gnu}}{\implies} \qquad   \tetta_*^K\ll \e^{|\log \e|}\ ,
\eeq
which is smaller than any power of $\e$ (but not exponentially small with $1/\e$).

\nl
{\rm (v)} 
If a set $D\subseteq \real^n$ is $(\a,K)$ non--resonant (mod $\L$) for $h=\|y\|^2/2$, then the complex domain $D_r$ is $(\a-rK,K)$ non--resonant (mod $\L$), provided\footnote{Indeed, if $y\in D_r$ there exists $y_0\in D$ such that $\|y-y_0\|<r$ and 
$|y_0\cdot k|\ge \a$ for all $k\in\integer^n\bks\L$, $|k|\le K$. Thus, $|y\cdot k|=|y_0\cdot k - (y_0-y)\cdot k|\ge |y_0\cdot k|-r K\ge \a-rK$.
} 
$rK<\a$.

\end{remark}

%%%%%%%%%%%%%%

%\subsection{Normal forms at resonances of order zero and one}

\noi
We now apply the Normal Form Lemma  to the Hamiltonian
$H$ in \equ{storia} in the non-resonant and simple resonant regions.

\subsection{Normal form in $\O^0$ (non--resonant regime)}
Recalling the definition of $\a$ given in \equ{lapparenza}, we set
\begin{equation}\label{timone}
r_{\{0\}}:=\frac{\a}{4K}=\frac14 \sqrt\e
K^\nu\,.
\end{equation}
We can apply Lemma~\ref{pesce} to $H$ in
\eqref{storia}
with\footnote{The set $\O^0$ is defined in \eqref{neva}. By Remark~\ref{cb500x}--(v)
and \eqref{cipollotto},
the domain $\O^0_{r_{\{0\}}}$ is $(\a/4,K)$ completely non--resonant.}: 
\beqano 
f\rightsquigarrow\e f\ ,\quad 
D\rightsquigarrow \O^0\ ,\quad
r\rightsquigarrow r_{\{0\}}\ ,\quad
\L\rightsquigarrow \{0\}\ , \quad
\a\rightsquigarrow \a/4\ ,
\eeqano
and 
\begin{equation}\label{ELP}
\tetta_* \rightsquigarrow \tetta_{\{0\}} :=
2^{11}n\frac{K^3 \e |f|_s}{\a r_{\{0\}} s}
\stackrel{\eqref{bada},\eqref{lapparenza},\eqref{timone}}=
\frac{2^{13} n}{s K^{2\nu-2}}\,.
\end{equation}
By \eqref{lapparenza}, $\tetta_{\{0\}}<1$, provided 
$\e$ is small enough depending on $s$ and $n$
(recall that $\nu>n+1$). 
Thus, there exists a symplectic change of variables
\begin{equation}\label{trota}
\Psi_{\{0\}}: \O^0_{r_{\{0\}}/2}\times \T^n_{s_*} \to 
\O^0_{r_{\{0\}}} \times \T^n_{s} 
\,,\qquad
s_*:=s(1-1/K)
\end{equation}
(recall \eqref{dadostar})
such that $H$  is transformed in 
\begin{equation}\label{prurito}
H_{\{0\}}:=H\circ\Psi_{\{0\}}
=\|I\|^2/2+\e g^{\{0\}}(I) +
\e f^{\{0\}}_{**}(I,\f) \,,\qquad
{\rm with}\quad
\langle f^{\{0\}}_{**}\rangle=0\,,
\end{equation}
where $\langle \cdot \rangle=\proiezione_{\{0\}} \cdot$ denotes the 
average with respect to the angles $\f$; by  \eqref{giovanni} and \equ{bada}, one has:
\begin{equation}\label{552}
\sup_{\O^0_{r_{\{0\}}}}| g^{\{0\}}-\langle f \rangle|
\leq
 2 \tetta _{\{ 0\}} 
\,,
\qquad
|f^{\{0\}}_{**} |_{r_{\{0\}}/2,s/2} 
\leq
2 e^{-Ks/2} 
\stackrel{\eqref{lapparenza}}=
2\, \e^{\frac{s}2 |\log \e|}
\end{equation}
provided $\e$ is  
small enough (depending on $s$ and $n$)
so that \eqref{francescone} is satisfied.

\subsection{Normal form in $\O^1$ (simple resonances)}
In order to construct normal forms near simple resonances, recall that $\O^1$ is the union of  sets\footnote{Recall the definitions given in  \eqref{LPD}, \eqref{caspio}, \eqref{neva}.}
$L_k Z_k$,
with $k\in \integer^n_{*,K}$, which are $(2\a \K/\|k\|,\K)$ non--resonant modulus the one--dimensional lattice $\integer k$;
compare Proposition~\ref{cb500f}, (iii). Therefore,  fixed $k\in\Z^n_{*,K}$, we let 
\begin{equation}\label{limone}
r_k:=\frac{\a}{\|k\|}= \frac{\sqrt\e
K^{\nu+1}}{\|k\|}\, ,
\end{equation}
and apply the normal form  Lemma~\ref{pesce}
with\footnote{
The symbol  ``$ a \rightsquigarrow b$'' reads ``with $a$ replaced by $b$''.
} 
\beq{boston}
f\rightsquigarrow\e f \ , \quad
D\rightsquigarrow D^k:=L_k Z_k\ , \quad 
r\rightsquigarrow r_k\ , \quad 
\a\rightsquigarrow \a \K/\|k\|\ ,
\quad
K\rightsquigarrow \K\ , \quad
\L\rightsquigarrow \Z k
\eeq 
and
\begin{equation}\label{ELPbis}
\tetta_* \rightsquigarrow \tetta_k 
\stackrel{\eqref{limone}}{:=}
2^{9}n\frac{\K^2\|k\|^2 \e |f|_s}{\a^2 s}
\stackrel{\eqref{bada},\eqref{lapparenza}}=
\frac{2^{9} n\K^2\|k\|^2}{s K^{2\nu+2}}\le 
\frac{2^{9} n}{s K^{2\nu-4}}
\, ;
\end{equation}
Notice that by Remark~\ref{cb500x}--(v)
and \eqref{cipollotto2},
the domain $D^k_{r_k}$ is $(2\a \K/\|k\|-r_k \K,\K)=(\a \K/\|k\|,\K)$  non--resonant modulus $\integer k$.
Again, by \eqref{lapparenza}, $\tetta_k<1$, provided 
$\e$ is small enough (depending on $s$ and $n$).
Thus,  there exists a symplectic change of variables\footnote{Recall \eqref{dadostar} and
Rematk~\ref{cb500x}, (iii).}
\begin{equation}\label{canarino}
\Psi_k: D^k_{r_k/2}\times \T^n_{s_*} \to 
D^k_{r_k} \times \T^n_{s} \,,
\qquad
s_*:=s(1-1/\K)
\end{equation}
such that $H$ in \eqref{storia} is transformed in 
\begin{equation}\label{Hk}
H\circ\Psi_k
=\|I\|^2/2+\e g^k(I,\f) +
\e f^k_{**} (I,\f)\,,
\end{equation}
where  
\begin{equation}\label{masiccio}
g^k=\proiezione_{{}_{k\Z}} g^k\,,\qquad
\proiezione_{{}_{k\Z}} f^k_{**}=0\,,
\end{equation}
with the following estimates holding
(recall \eqref{giovanni} and \eqref{bada}): 
\begin{equation}\label{552bis}
| g^k-\proiezione_{{}_{k\Z}} f|_{r_k/2,s_*}
\leq
2 \tetta 
\,,
\qquad
|f^{k}_{**} |_{r_k/2,s/2} 
\leq
2e^{-\K s/2}\stackrel{\eqref{lapparenza}}=
2\, \e^{\frac{s}2 |\log \e|^3}\ ,
\end{equation}
where
\beq{tetta}
\tetta:=\frac{2^{9} n}{s K^{2\nu-4}}\ ,\qquad (\nu>n+1)\ .
\eeq
Note that 
$g^k$  and $\proiezione_{{}_{k\Z}}f$ depend, effectively,  only on one angle  $t\in\T^1$: more precisely, setting
\begin{equation}\label{sugna}
F^k_j:=f_{jk}\,,\ \ \  	 G^k_j(I):=g^k_{jk}(I)\,,
\quad {\rm and}\quad 
F^k(t):=\sum_{j\in\Z} F^k_j e^{\ii j t}\,,\ \ \ 
 G^k(I,t):=\sum_{j\in\Z} G^k_j(I) e^{\ii j t}\,,
\end{equation}
we have (recall \eqref{dec})
\begin{equation}\label{palettabis}
\proiezione_{{}_{k\Z}} f(\f)=
F^k(k\cdot \f)\,,\qquad
g^k(I,\f)=G^k(I,k\cdot\f)\,.
\end{equation}
We also remark  that
since $f\in  \cA_s^n$, the functions 
$F^h$ belong to $ \cA_{|h|s}^1$ for every $h\in\Z^n_*,$ with 
\begin{equation}\label{trilli}
|F^h|_{|h|s}\leq |f|_s \eqby{bada}1\,,
\end{equation}
Analogously, by \eqref{552bis}
\begin{equation}\label{cristina}
|G^k-F^k|_{r_k/2,|k|s_*} \leq 2 \tetta \,.
\end{equation}
For later use (compare \equ{allblacks} and \equ{neftali} below), we point out that
\beqa{inbloom}
\frac{1}{|f_k|}|G^k_1-f_k|_{r_k/2}&\leby{cristina}&
\frac{2 \tetta e^{-|k|s_*}}{|f_k|}\stackrel{\eqref{canarino}}=\frac{2 e^{s|k|/\K} \tetta e^{-|k|s} }{|f_k|}
\leq
\frac{2  \tetta e^{(1-|k|)s} }{|f_k|}
\nonumber\\
&\stackrel{\rm (P1)}{\le}&
\frac{2  \tetta e^s |k|^{\frac{n+3}{2}}}{\d}
\leq
\frac{2  \tetta e^s K^{\frac{n+3}{2}}}{\d} 
\stackrel{\eqref{tetta}}=
\frac{2^{10} n e^s}{\d s K^{\frac{4\nu-n-11}{2}}}
\,,
\eeqa
which is small if 
$\e$ is small  (recall that $n\ge 2$ and $\nu>n+1$).

\subsubsection{The effective potential}

We now show that for $|k|$ large, the ``effective potential'' $G^k$ (defined in \eqref{sugna}) behaves, essentially, as  a cosine; compare, in particular,
Eq. \equ{allblacks} below.

\nl
Recalling the definition of $T_h$ given in\eqref{labestia},
we have that
$$T_1 F^k(\psi_n')=f_k e^{\ii \psi_n'}+f_{-k} e^{-\ii \psi_n'}= 2|f_k| \cos(\psi_n' +\psi_n^{(k)})$$
for suitable constants $\psi_n^{(k)}$.

\begin{remark}\label{coppa} 
{\sl We can assume, up to translation,  that $\psi_n^{(k)}=\pi$}.
%\beq{italiaspagna20}
%\psi_n^{(k)}=\pi
%\ ,\qquad  \forall \ K_s(\d)<|k|\leq K 
%\ .
%\eeq 
%compare \eqref{wasabi} below.
So, from now on, we assume that
\begin{equation}\label{sumo}
T_1 F^k (\psi_n')=
f_k e^{\ii \psi_n'}+f_{-k} e^{-\ii \psi_n'}= -2 |f_k| \cos(\psi_n')
%\,,\qquad \forall\, K_s(\d)<|k|\leq K
\,.
\end{equation}
\end{remark}
\begin{lemma}
\begin{eqnarray}
&&G^k(I,\psi'_n)=2|f_k|\Big( -\cos(\psi_n')+
 R^k(I,\psi'_n)
\Big)\,,
\nonumber
\\
&&|R^k(\sqrt{2\e/\kappa} A_k^TJ',\psi'_n)|_{D^k,r_k/2,|k|s/3}
\leq 
\frac{2|k|^{\frac{n+3}{2}} e^{-|k|s/4}}{\d}+
\frac{2^{10} n e^s}{\d s K^{\frac{4\nu-n-11}{2}}}
\,.
\label{allblacks}
\end{eqnarray}
\end{lemma}
\proof  Indeed, 
%$\forall\, K_s(\d)<|k|\leq K,$
\begin{eqnarray}
&&\frac{1}{|f_k|}
|F^k(\psi_n')+ 2|f_k|  \cos(\psi_n')|_{|k|s/3}
\eqby{sumo}
\frac{1}{|f_k|}
|T^\perp_1 F^k|_{|k|s/3}
=
\frac{1}{|f_k|}\sup_{|j|\geq 2} |f_{kj}|e^{|k|s|j|/3}
\nonumber
\\ 
&&
\stackrel{\eqref{bada}}\leq
\frac{1}{|f_k|}\sup_{|j|\geq 2}e^{-2|k|s|j|/3}
=
\frac{1}{|f_k|}e^{-4|k|s/3}
\stackrel{(P1)}\leq 
\frac{|k|^{\frac{n+3}{2}} e^{-|k|s/3}}{\d}\,.
\label{pizzetta}
\end{eqnarray}
Also, 
\begin{eqnarray}\label{pizzetta2}
&&\frac{1}{|f_k|}
\Big|T^\perp_1\big(G^k(I,\psi'_n)
-F^k(\psi_n')\big)\Big|_{D^k,r_k/2,|k|s/3}\stackrel{\eqref{cristina}}
\leq
\frac{2\tetta}{|f_k|}
\sup_{|j|\geq 2}
e^{-|j||k|(s_*-s/3)} 
\nonumber
\\
&&=
\frac{2\tetta}{|f_k|}
e^{-2|k|(s_*-s/3)} 
\stackrel{(P1)}\leq
\frac{2\tetta|k|^{\frac{n+3}{2}} e^{-|k|(2s_*-5s/3)}}{\d}
\stackrel{\eqref{canarino}}\leq
\frac{2\tetta|k|^{\frac{n+3}{2}} e^{-|k|s/4}}{\d}\,,
\end{eqnarray}
provided $\K\ge 24$.
Then \eqref{allblacks} follows by \eqref{inbloom},\eqref{pizzetta} and
\eqref{pizzetta2}. \qed

\nl
Moreover by \eqref{inbloom},  \eqref{pizzetta2}, we have also that
%again $\forall\, K_s(\d)<|k|\leq K$,
\begin{equation}\label{neftali}
|f_k|^{-1}
|G^k(I,\psi'_n)-
F^k(\psi'_n)|_{D^k,r_k/2,|k|s/3} 
\leq
\frac{4  \tetta e^s K^{\frac{n+3}{2}}}{\d} 
\stackrel{\eqref{tetta}}=
\frac{2^{11} n e^s}{\d s K^{\frac{4\nu-n-11}{2}}}
\,,
\end{equation}
which is small if 
$\e$ is small  (recall that $n\ge 2$ and $\nu>n+1$).

\subsubsection{Rescalings}
Recalling the definition of $K_s(\d)$ in \eqref{miguel}, we set
\begin{equation}\label{cenerentola}
%\d_k=\d_{k,s}:=
\d_k:=
\left\{ \begin{array}{ll} 1 &
\ \ {\rm if} \ \ |k|\leq  K_s(\d)\\
2|f_k| & \ \ {\rm if} \ \ |k|>  K_s(\d)
 \end{array}\right.\,.
\end{equation}
Note that by \eqref{bada}
\begin{equation}\label{sassa}
\d_k\leq 1\,.
\end{equation}

\nl
Define the conformally symplectic transformation 
\begin{equation}\label{prospettiva}
\Phi^{(0)}:(I',\f')\mapsto (I,\f)=
(\riscala I',\f')\,,
\end{equation}
where 
\beq{capocotta}
\riscala := \sqrt{\frac{2\d_k\e}{\kappa}}=
\sqrt{\frac{2\d_k\e}{\|k\|^2}}
 \ .
\eeq
Then, the flow of $H\circ\Psi_k$ (recall \eqref{Hk} and \equ{palettabis}) is equivalent to the flow of the Hamiltonian\footnote{See Lemma~\ref{hans}, (ii) in Appendix~\ref{noccetti}).}
$$
\frac1{\riscala} H\circ\Psi_k\circ\Phi^{(0)}(I',\f')
=\sqrt{\frac{\d_k\e}{2\kappa}} \|I'\|^2+
\sqrt{\frac{\kappa\e}{2\d_k}} \big( G^k(\riscala I',k\cdot\f')
+ f^k_{**}(\riscala I',\f')  \big) \,.
$$
Dividing such Hamiltonian by $\sqrt{\d_k\kappa\e/2}$ (which corresponds to  a time rescaling\footnote{See, again, Lemma~\ref{hans}.}),
we are lead to study the  Hamiltonian
\begin{equation}\label{disfare}
H_k(I',\f'):=\frac{1}{\kappa} \|I'\|^2+\frac{1}{\d_k}
\left(
 G^k(\riscala I',k\cdot\f')
+ f^k_{**}(\riscala I',\f')\right) \,,
\end{equation}
which is  defined on the domain
\begin{eqnarray}\label{fare}
&&D^{'k}_{r'_k}\times \T^n_{s_*}\,,\ \  \ 
{\rm with}\qquad
r'_k:=\frac{r_k}{2\riscala}
 \stackrel{\eqref{limone}}=\frac{K^{\nu+1}}{\sqrt {8\d_k} }\qquad
 {\rm and}
 \\
 &&D^{'k}:=\frac{1}{\riscala}D^k
\stackrel{\eqref{boston}}=
\frac{1}{\riscala}L_k Z_k
\stackrel{\eqref{caspio},\eqref{lapparenza},\eqref{capocotta}}=
L_k 
\left(
\frac{1}{\riscala}\hat Z_k\times\Big(-\frac{K^{\nu+1}}{\sqrt{2\d_k}\|k\|},
\frac{K^{\nu+1}}{\sqrt{2\d_k}\|k\|}\Big)
\right)\,.
\nonumber
\end{eqnarray}
Note that, by \eqref{552bis}
\begin{equation}\label{avezzano2}
|\d_k^{-1}f^k_{**}(\riscala I',\f')
 |_{D^{'k},r'_k,s/2} 
\leq
\frac{2}{\d}|k|^{\frac{n+3}{2}} e^{|k|s}e^{-\K s/2}
\stackrel{\eqref{lapparenza}}\leq 
\frac{1}{\d}e^{-\K s/4}=
\frac{1}{\d} \e^{\frac{s}4 |\log \e|^3}
\end{equation}
for $\e$ small enough.

\nl
Recalling \eqref{shevket} we set
\begin{equation}\label{foggydew}
D^{'k}_\sharp:=
\frac{1}{\riscala}L_k Z_k^\sharp\,.
\end{equation}
Note that 
\begin{equation}\label{foggydew2}
\Phi^{(0)}\big(D^{'k}_\sharp\times\T^n\big)=
L_k Z_k^\sharp\times\T^n\,.
\end{equation}
Then by \eqref{neva}
\begin{equation}\label{foggydew3}
\Phi^{(0)}\bigg(\bigcup_{k\in\Z^n_{*,K}}D^{'k}_\sharp\times\T^n\bigg)=
\bigcup_{k\in\Z^n_{*,K}}\Phi^{(0)}\big(D^{'k}_\sharp\times\T^n\big)=
\O^1\times\T^n\,.
\end{equation}

\section{The nearly--integrable structure at simple resonances}
\label{erezione}

%%%%%%%%%%%%%%%%%
%%%%%%%%%%%%%
%%%%%%%%%%%%%%%%
%%%%%%%%%%%%%%%%
%%%%%%%%%%%%%%%%
%%%%%%%%%%%%%%%%
%%%%%%%%%%%%%%%%
%%%%%%%%%%%%%%%%
%%%%%%%%%%%%%%%%
%%%%%%%%%%%%%%

%%%%%%%%%
%%%%%%%%%%
%%%%%%%%%
%%%%%%%%%%
%%%%%%%%%%%
%%%%%%%%%%%
%%%%%%%%%%%%
%%%%%%%%%%%%
%%%%%%%%%%%
%%%%%%%%%%%

\Giu
Given a bounded holomorphic function 
$f:D_{r}\times \T^n_s\to\C,$ with $D\subseteq \R^{n'}
$
we set
\begin{equation}\label{ionico}
\|f\|_{D,r,s}=\|f\|_{r,s}:=\sup_{D_r\times \T^n_s}|f|\,.
\end{equation}
The following relation between the two norms
$|\cdot|$ and $\|\cdot\|$ holds:
for $\s>0$, we have\footnote{
Since $\sum_{k\in\Z^n}e^{-|k|\s}=\Big(
\sum_{k\in\Z}e^{-|k|\s}\Big)^n
=
\Big(1+ 2\sum_{j\geq 1}e^{-j\s}\Big)^n
=
\Big(\frac{e^\s+1}{e^\s-1}\Big)^n
= \coth^n(\s/2).$}
\begin{equation}\label{battiato2}
|f|_{r,s}\leq \|f\|_{r,s}\leq \coth^n(\s/2)|f|_{r,s+\s}
\leq
(1+2/\s)^n
|f|_{r,s+\s}\,.
\end{equation}

\subsection{A class of Morse non-degenerate functions}

Let $s_0>0$ and let us consider a bounded holomorphic function
\begin{equation}\label{latte}
F^0:\T_{s_0}\to\C\,,\qquad \text{with}\qquad
\|F^0\|_{s_0}<\infty\,.
\end{equation}

\begin{definition}\label{morso}
 Let $\b,M>0.$ We say that $F^0$ as in \eqref{latte} is
 $(\b,M)$-Morse-non-degenerate if
 $\|F^0\|_{s_0}\leq M$ and 
 \begin{eqnarray}
\label{Lbeta1}
\min_{x} \ \big( |(F^0 )'(x)|+|(F^0 )''(x)|\big) 
&\geq& \b\,,
\\
\label{Lbeta2}
\min_{1\leq i<j\leq 2N } 
|F^0(x^0_i)-F^0(x^0_j)|
&\geq& \b\,,
\end{eqnarray}
where $x^0_i$, $1\leq i\leq 2N$
are the critical points of $F^0$
 in $(-\pi,\pi].$
\end{definition}

We note that, by \eqref{Lbeta1},
the function $F^0 $ has only non-degenerate critical points:
let us say
$N $ minima: $x^0_{2j-1},$
 and $N $ maxima: $x^0_{2j},$ 
 in $(-\pi,\pi],$ for some integer $N \geq 1$ and $1\leq j\leq N.$
 It is immediate to realize   that
$
N 
$
is uniformly bounded by a constant depending only on 
$s_0, M$ 
and the minimum appearing in \eqref{Lbeta1}.
\\
The corresponding critical energies are 
\begin{equation}\label{cinzia}
E^{0}_i:=F^0 (x^{0}_i)\,,\qquad
1\leq i\leq 2N \,.
\end{equation}
By \eqref{Lbeta2}, $E^{0}_i$ are all different.

\begin{definition}\label{sorellastre}
 We say that $F^0$ as in \eqref{latte} is
 $\g$-cosine-like\footnote{Actually we should say {\sl minus-cosine-like}} if
 \begin{equation*}
 \|F^0(x)+\cos x \|_{s_0}\leq \g\,,
 \qquad \text{for some}\qquad 
 0<\g\leq \frac14\min\{1,s_0^2\}\,.
\end{equation*}
\end{definition}

\begin{lemma}\label{sibillini}
 If $F^0$ is
 $\g$-cosine-like, then it is also
 $(\b,M)$-Morse-non-degenerate
 with
 $$
 \b=1/4\,,  \qquad
 M=\g+\cosh s_0\leq \frac14+\cosh s_0\,.
 $$
 Moreover $F^0$
 has only two non-degenerate critical points
(a maximum and a minimum).
\end{lemma}
\proof
We have, by Cauchy estimates,
$$
 |(F^0 )'(x)|+|(F^0 )''(x)|
 \geq |\sin x| +|\cos x| 
 -\frac{\g}{s_0}-2\frac{\g}{s_0^2}
 \geq 1-\frac{\g}{s_0}-2\frac{\g}{s_0^2}\geq \frac14\,.
$$
We can choose $M$ as above since
$\|\cos x\|_{s_0}=\cosh s_0.$
Regarding the last sentence
we note that for $x\in(-\pi,\pi]$
we have only two critical points,
a minimum in $(-\pi/6,\pi/6)$
and a maximum in $(-\pi,-5\pi/6)\cup
(5\pi/6,\pi].$ Indeed we have that, 
setting $g(x):=F^0(x)+\cos x$,
$(F^0)'(x)=\sin x+g'(x),$ so that
\begin{equation}\label{coratella}
(F^0)'(x)=\sin x+g'(x)
\geq \sin x-\g/s_0\geq \sin x -1/4\,.
\end{equation}
This implies that 
$(F^0)'(\pi/6)\geq 1/4,$ 
$(F^0)'(-\pi/6)\leq -1/4.$
Then, by continuity, there exists a critical
point of $F^0$ in $(-\pi/6,\pi/6).$
Moreover such point is a minimum and there are no other critical points in $(-\pi/6,\pi/6)$ since there $F^0$
is strictly convex:
$$
(F^0)''(x)=\cos x+g''(x)
\geq \sqrt 3/2-2\g/s_0^2\geq \sqrt 3/2 -1/2
>0\,.
$$
Similarly in $(-\pi,-5\pi/6)\cup
(5\pi/6,\pi]$ there is only one critical point,
which is a maximum.
Finally, by \eqref{coratella},
$(F^0)'(x)\geq 1/4$ for $x\in[\pi/6,5\pi/6]$
and, analogously, 
$(F^0)'(x)\leq -1/4$ for 
$x\in[-5\pi/6,-\pi/6];$
so that there are no other critical points.
\qed

\subsection{The Structure Theorem}\label{lupin1}

\nl
We start introducing a parameter 
\begin{equation}\label{vaccinara}
\theta\geq 0\,,
\end{equation}
that will be chosen in Section \ref{lupin3} as a function of 
$\e.$
We also say that a function $\theta \to E(\theta)\subseteq \R^m$
is {\sl decreasing w.r.t. $\theta$} if
$\theta\leq \theta'$ implies $E(\theta)\supseteq E(\theta').$

In light of \eqref{disfare} we are now 
going to study the 
behavior of  the ``effective Hamiltonian''
close to a simple resonance identified by 
a fixed $k\in\Z^n_*,$ namely we
are considering Hamiltonian of the form
\begin{equation}\label{golia}
\cH(I',\f'):=\frac1\kappa \|I'\|^2
+\Ggot(I',k\cdot \f')\,.
\end{equation}

\noindent
\textbf{Assumptions on the ``effective  potential'' $\Ggot $}\\
Consider the parameters
\begin{equation}\label{santostefano}
s_0\,,\,  r_0>\,0\,,\ \ 
r':=c n |k|_\infty r_0\,,
\end{equation}
$c>1$ being the constant defined in \eqref{verruca}, which depends only on $n$.
We will make the following assumptions:

\noindent
\textbf{(A1)} There exists 
\begin{equation}\label{sgravone}
F^0\in\cA^1_{s_0}
\end{equation}
such that
\begin{equation}\label{Ccristina5}
\|\Ggot  - F^0 \|_{\cD,r',s_0}
\leq 
\checco_*
\,,
\end{equation}
with\footnote{Recall the definition of 
$L_k$ given in \eqref{LPD}.}
\begin{equation}\label{calimero}
\cD=L_k\big(\hat D \times (-R_0 ,R_0 )
\big)\,,\ \ \ 
\hat D \subset \R^{n-1}\,,\qquad
R_0\geq 4+\cosh s_0
\,;
\end{equation}

\begin{flalign}\label{legna}
\textbf{(A2)} \ F^0\  \text{is} \ 
(\b,M)\text{-Morse-non-degenerate with}\ \ 
\max\{2\sqrt M,4\}\leq R_0\,. &&
\end{flalign} 
In alternative to \textbf{(A2)}     
we will assume, when it holds, the following
{\sl stronger} (recall Lemma \ref{sibillini} and \eqref{calimero}) hypothesis:
\begin{flalign}\label{drago}
\textbf{(A3)} \ F^0\  \text{is} \ 
\g\text{-cosine-like
with}\ \ 
\g\leq \cgot(s_0)
:=\cgot_*\min\{1,s_0^4\}\,, &&
\end{flalign}
where $0<\cgot_*\leq 1/4$ is a suitably small positive
constant to be chosen below
(see Lemma \ref{lapajata}).

\medskip

\noindent
Recalling \eqref{calimero} we set
\begin{equation}\label{calimero2}
\cD_\sharp:=L_k\big(\hat D \times (-R_0 /2,R_0 /2)
\big)\,,\ \ \ 
\hat D \subset \R^{n-1}
\,.
\end{equation}

\begin{theorem}[Integrable structure at simple resonances]\label{porretta}$\phantom{.}$\\

\smallskip\noindent
{{\rm \bf Part I.}} Assume that $\Ggot$  in \eqref{golia} satisfies \textbf{(A1)}
and \textbf{(A2)} or \textbf{(A3)}.   
Then there exist   a suitably large constant $\mathtt c>1,$ which, when \textbf{(A3)} holds, 
depends only on $n,s_0,r_0$,
otherwise it depends also on\footnote{In any case it is independent of $\hat D.$} $F^0,$ 
 such that if
\begin{equation}\label{tolomeo}
\checco_*
\leq 1/\mathtt c
\end{equation}
the following holds.
For every\footnote{$2N$ being the number of critical points of $F^0.$
Note that, when \textbf{(A3)} holds,
$N=1$
by Lemma \ref{sibillini}.
 } 
$0\leq i\leq 2N$ there exist
\\
i) disjoint open subsets\footnote{The set $\bigcup_i \Cgot^i(0)$ contains  $\cD_\sharp\times \T^n$ 
up to the connected components of the critical energy level containing critical points.}
$\Cgot^i(\theta)\subseteq \cD\times\T^n$ 
decreasing w.r.t. $\theta$, with
\begin{equation}\label{fujiyama}
\meas\Big( \big(\cD_\sharp\times \T^n\big) \setminus \bigcup_{0\leq i\leq 2N}
\Cgot^i(\theta) \Big) \leq \mathtt c \theta |\ln\theta|\,;
\end{equation}
ii) open subsets $\mathtt B^i(\theta)\subseteq\R^n$, decreasing w.r.t. $\theta$ with\footnote{
Where $c$ is the constant defined in \eqref{verruca}.}
\begin{equation}\label{pratola2}
{\rm diam} \big( \mathtt B^i(0)\big)
\leq 2c \big(R_0+{\rm diam}(\hat D)\big)\,,
\qquad \forall\, 0\leq i\leq 2N\,;
\end{equation}
iii) 
a symplectomorphism 
\begin{equation}\label{tokyo}
\Psi^i \ :\ (p,q)\in\mathtt B^i (0)\times\T^n \ \to \ \Cgot^i (0)
\ni(I',\f')\,,\qquad
\mbox{with}\quad
\Psi^i  \Big(\mathtt B^i (\theta)\times\T^n \Big)=
\Cgot^i (\theta)\,,\ \ \ \forall\,\theta\geq 0\,,
\end{equation}
with holomorphic extension 
\begin{equation}\label{blueeyes3}
\Psi^i \ :\ \big(\mathtt B^i (\theta)\big)_{\rho_*}\times\T^n_{\s_*}
\ \to \ \cD_{r'}\times\T^n_{s_0}
\,,\qquad
\mbox{with}\quad
\rho_*:=\frac\theta{\mathtt c |k|_\infty^{n-1}} \,,
\quad  \ \s_*:=\frac1{{\mathtt c}|k|_\infty^{n-1}
|\log \theta|}\,,
\end{equation}
such that 
 \begin{equation}\label{matteo}
\cH
\circ \Psi^i(p,q)=:
h ^{(i)}(p)\,.
\end{equation}
Moreover
\begin{equation}\label{vallinsu2}
\|\partial_{pp} h^{(i)} \|_{ \mathtt B^i (\theta), \rho_*}\leq \mathtt c/\theta\,,\qquad
\text{for }\ \ \ 0\leq i\leq 2N\,.
\end{equation}

\noindent
{{\rm \bf Part II.}}
Assuming 
\begin{equation}\label{gaetanga}
0<\checco_*\leq 1/\mathtt c\|k\|^{2n}\,,
\end{equation}
we have that
for every
\begin{equation}\label{gaetanga2}
0<\mu\leq 1/\mathtt c\|k\|^{2n}\,,
\end{equation}
 there exist
open subsets $\tilde{\mathtt B}^i (\mu)\subseteq
\mathtt B^i (0)$, decreasing w.r.t. $\mu,$
such that
\beq{gibboso2}
\meas(\mathtt B^i (0)\bks \tilde{\mathtt B}^i (\mu))\le 
\mathtt c \|k\|^{4n} \mu^{1/\mathtt c}
\eeq
and
\begin{equation}\label{gadamer}
\left| {\rm det}\left(\partial_{pp} h ^{(i)}(p)
    \right)  \right|
> \mu\,,\qquad
\forall\, \ 0\leq i\leq 2N \,,\ \ |k|\leq K\, , \qquad \forall \ p\in  \tilde{\mathtt B}^i (\mu)\ .
\end{equation}

\end{theorem}

\nl
The following two  sections are devoted to the proof
of Theorem \ref{porretta} part I  and part 
II,
respectively.

\section{Proof of Part I of the Structure Theorem}\label{riso}

In this section we will prove
Theorem \ref{porretta} part I.

\subsection{Critical points and critical energies of 
the ``unperturbed potential'' $F^0$}\label{farro}

%Up to a translation\footnote{BISOGNEREBBE MOSTRARE CHE SE $F^0$ E' $\g$-cosine-like RIMANE TIPO, $2\g$-cosine-like!!!!!!!}, 
% \begin{equation}\label{wasabik}
 %\text{we can  assume  that the absolute maximum of $F^0 $ is attained at $x=\pi.$}
%\end{equation}
We order the critical points of $F^0 $
(recall \eqref{sgravone}) 
 in the following way (where\footnote{Similarly we will set $E_0^{0}:=E_{2N}^{0}$ below.}
  $x^{0}_0:=x^{0}_{2N}-2\pi$) 
 \begin{equation}\label{recremisik}
x^{0}_0< x^{0}_1<x^{0}_2<\ldots <x^{0}_{2N-1}<x^{0}_{2N}\,,\qquad
 x^{0}_{2j-1}\  {\rm minimum}\,,\ \ \ 
  x^{0}_{2j}\  {\rm maximum}\,, \ \ 1\leq j\leq N \,.
\end{equation}
Fix $1\leq j\leq N$ and consider a minimum point $x^0_{2j-1},$
thanks to \eqref{Lbeta1} the function $F^0$ is strictly increasing, resp. strictly decreasing, in the interval $ [x^0_{2j-1},x^0_{2j}],$
resp. $ [x^0_{2j-2},x^0_{2j-1}],$
then we can invert $F^0$ on the above intervals 
obtaining two functions
\begin{equation}\label{andria}
X^0_{2j} : [E^0_{2j-1}, E^0_{2j}]\to  [x^0_{2j-1},x^0_{2j}]
\qquad {\rm and}\qquad
X^0_{2j-1} : [E^0_{2j-1}, E^0_{2j-2}]\to [x^0_{2j-2},x^0_{2j-1}]
\end{equation}
such that
$$
F^0 (X^0_i(E))=E\,,\qquad
X^0_i(F^0(\psi_n))=\psi_n
\,,\qquad \forall\, 1\leq i\leq 2N\,.
$$
Note that $X^0_i$ is increasing, resp. decreasing, if 
$i$ is even, resp. odd.

Set 
\begin{eqnarray}\label{timur}
&&E^{(i),0}_-:=E^0_i\,,\qquad
E^{(2j-1),0}_+:=
\min\{ E_{2j-2}^0, 
E_{2j}^0\}\ \ {\rm for}\ \ 1\leq j\leq N\,,
\nonumber
\\
&&E^{(2j),0}_+:=\min\{ E^0_{2j_-}, E^0_{2j_+}\}
\ \ {\rm for}\ \ 1\leq j< N\,,\qquad
E^{(2N),0}_+=E^{(0),0}_+=+\infty
\,,
\end{eqnarray}
where
\begin{equation}\label{21stcentury}
j_-:=\max\{ i<j\ \ {\rm s.t.}\ \ E^0_{2i}>
 E^0_{2j} \}\,,\qquad
 j_+:=\min\{ i>j\ \ {\rm s.t.}\ \ E^0_{2i}>
 E^0_{2j} \}\,.
\end{equation}

\subsection{The slow angle}
Let us, now,  perform the linear symplectic change of variables
$\Phi^{(1)}:(J',\psi')\mapsto (I',\f')$
generated by 
$S(J',\f'):=A_k\f'\cdot J',$:
\begin{equation}\label{talktothewind}
\Phi^{(1)}: (J',\psi')  \mapsto  (I'\f')=
(A_k^TJ',  A_k^{-1} \psi')=(k J'_n+\hat A_k^T \hat J',A_k^{-1} \psi')\,.
\end{equation}
Note that $\Phi^{(1)}$ does not mix actions with angles, its projection on the angles is a diffeomorphism of $\T^n$ onto $\T^n$,
and, most relevantly, 
\beq{canonico}
\psi'_n=k\cdot \f'
\eeq
is the canonical angle associated to the one-dimensional ``secular system'' near the simple resonance $\{y\cdot k=0\}$ 
(i.e., the one-dimensional system governed by the Hamiltonian obtained disregarding the small term $f_{**}^k$ in \equ{disfare}).

\nl
In the $(J',\psi')$--variables, we have:
\begin{equation}\label{colosseum3}
\tilde H (J',\psi'):=\cH\circ \Phi^{(1)}(J',\psi')
=
\frac{1}{\kappa} 
\| A_k^TJ' \|^2+
 \Ggot( A_k^TJ',\psi'_n)\,.
 \end{equation}
As for the $(J',\psi')$--domain,  we see that
the real $J'$--domain is given by
\begin{equation}\label{formentera}
\tilde D:= A_k^{-T}\cD=
 A_k^{-T} 
L_k\big(\hat D \times (-R_0 ,R_0 )
\big)
  \eqby{LPDbis} 
U_k\big(\hat D \times (-R_0 ,R_0 )
\big)
\,,
\end{equation}
with  $U_k$ in \eqref{centocelle}.
We also set
\begin{equation}\label{formentera2}
\tilde D_\sharp:=
U_k\big(\hat D \times (-R_0 /2,R_0 /2)
\big)
\,.
\end{equation}
Then, if we choose
\beq{pelotas}
\tilde r:=\frac{r'}{ n |k|_\infty}
\ ,\qquad
\tilde s:=
\frac{s_0}{c |k|_\infty^{n-1}}\,,
\eeq
(for a suitable $c$ depending only on $n$)
we see that $\Phi^{(1)}$ has holomorphic extension on the complex domain
(recall Lemma \ref{pergamena})
\begin{equation}\label{fare2}
\Phi^{(1)}: \tilde D_{\tilde r}\times \T^n_{\tilde s} \to \cD_{r'}\times \T^n_{s_0}
\end{equation}
indeed:  by \eqref{scimmia},  $\|A_k^T\|=\|A_k\|\leq n |k|_\infty$,
so that $\|A_k\| \tilde r\leq r'$, while, for every $1\leq i\leq n,$ 
$$
\sum_{1\leq j\leq n} |(A_k^{-1})_{ij} |
\leq n  |A_k^{-1}|_\infty
\stackrel{\eqref{atlantide}}\leq 
c |k|_\infty^{n-1}\,,
$$

\nl
Note that by \eqref{Ccristina5}
\begin{equation}\label{cristina2}
	\|\Ggot(A_k^TJ',\psi'_n)  - F^0(\psi'_n)
 \|_{\tilde D,\tilde r,s_0}
\leq 
\checco_* \,.
\end{equation}

\subsection{The auxiliary Hamiltonian}

A crucial role will be played by the
{\sl auxiliary Hamiltonian}
\begin{equation}\label{spigola}
H^*:= (J_n'')^2
+
F^*(J'',\psi''_n)\,,\qquad\quad
\text{where}\qquad
F^*(J'',\psi''_n):=
\Ggot( L_k J'',\psi''_n)
\,.
 \end{equation}
This Hamiltonian represents a  
one dimensional mechanical system depending on the parameter $\hat J''.$
The relation between $\tilde H$
(defined in \eqref{colosseum3})
and $H^*$ is the following: 
recalling \eqref{LPD},\eqref{LPDbis} and  the change $J'=U_k J''$, $U_k$ defined in \eqref{centocelle},
 it results
\begin{equation}\label{spigolo}
H^*(J'',\psi'')=
\tilde H(U_k J'',\psi_n'') 
- \frac{1}{\kappa} \|\proiezione_k^\perp \hat A_k^T \hat J''\|^2\,.
\end{equation}
\\
Recalling \eqref{formentera}, the potential $F^*$
in \eqref{spigola} is defined 
for  
\begin{equation} \label{sontuosa}
(J'',\psi_n'')\in D_{r_0}\times\T^1_{s_0}\,,\qquad \text{where}\qquad
D:=
\hat D \times (-R_0 ,R_0 )\,,
\end{equation}
\eqref{santostefano},\eqref{pelotas} and \eqref{verruca}. 
Note that, by \eqref{cristina2},
\begin{equation}\label{cristina7}
\|F^*  - F^0\|_{D,r_0,s_0}
\leq 
\checco_* \,.
\end{equation}

%%%%%%%%%%%%
%%%%%%%%
%%%%%%%%%%
%%%%%%%%%%%%%

%%%%%%%%%%%%%%%%

\subsection{A
special group 
of symplectic transformations}

%
%Recollecting we have defined all the parameters in \eqref{pgot2} 
%\textsl{independently} of $k$.
% Then we \textsl{simultaneously}, namely for every $|k|\leq K,$ apply 
% the results of Section \ref{gerusalemme}
%  to the Hamiltonian $H^*\rightsquigarrow H^*_k$ (defined in \eqref{spigola}
% and recall \eqref{mambo}).
%
%%%%%%%%%%%%%%
%%%%%%%%%%%%%

\nl
In the following,  symplectic transformations will  have a special form, namely, they will belong to a special group 
${\cal G}$, formed by symplectic transformations $\Phi$ satisfying
 \begin{equation}\label{grillo}
\hat I=\hat J\,,\quad I_n=I_n(J,\psi_n)\,,\quad
\hat\f=\hat\psi+\hat\f'(J,\psi_n)\,,\quad
\f_n=\f_n(J,\psi_n)\,, 
\end{equation} 
where, in general,  $\f,\psi$ may belong either to $\T^n$ or to $\R^n$. 
For  a transformation $\Phi$ as in \eqref{grillo}
we let
$\check \Phi$ 
denote the map
\begin{equation}\label{islands}
\check\Phi(J,\psi_n):=\big( 
\hat J, I_n(J,\psi_n), \f_n(J,\psi_n)
\big)\,.
\end{equation}
Some general properties of ${\cal G}$ are discussed in Appendix~\ref{noccetti}.

\subsection{An intermediate transformation}\label{orzo}

To simplify geometry, we now introduce a symplectic transformation that removes the dependence upon 
$J_n$ from the potential.  
Since \eqref{tolomeo} holds\footnote{
Note  that this implies   
\eqref{urea}.   
}, by Lemma~\ref{francesco} in Appendix~\ref{noccetti}, 
one can find a symplectomorphism
$\Phi^{(2{\rm bis})}\in \mathcal G$
satisfying
\begin{equation}\label{valentinesuite2}
\Phi^{(2{\rm bis})}: (J,\psi)\to(J'',\psi'')\,, \qquad
\hat J''=\hat J,\ \ J_n''=J_n+ a_* (\hat J,\psi_n),\ \  
\hat\psi''=\hat\psi+ b_*  (\hat J,\psi_n),
\ \psi_n''=\psi_n\,,
\end{equation}
 with (taking $\mathtt c$ large enough)
 \begin{equation}\label{flaviano2}
\Phi^{(2{\rm bis})} :\ D_{r_0/2}\times
 \T^{n}_{s_0/2}
 \ \to\ 
 D_{r_0}\times\T^n_{s_0}\,.
\end{equation}
 and
\begin{equation}\label{urina2}
\|a_* \|_{\hat D_k,r_0,s_0}\leq 4\checco_*/r_0\,,\ 
\|b_* \|_{\hat D_k,r_0/2,s_0}
\leq (16\pi+8)\checco_*/r_0^2\,,
\end{equation}
and such that
\begin{eqnarray}\label{guglielmo}
&&\Hpend(J,\psi_n):=
H^*\circ \Phi^{(2{\rm bis})} =\big(1+ b (J,\psi_n ) \big)\big(J_n-  J_n^*(\hat J )\big)^2
+ F (\hat J,\psi_n)\,, 
\\
&& F (\hat J,\psi_n)= F^0  (\psi_n )+
G (\hat J,\psi_n)\,.
\nonumber
\end{eqnarray}
Furthermore (see \equ{straussbis} below) 
\begin{eqnarray}\label{mazinga}
&&\|J_n^*\|_{\hat D ,r_0}
\leq 2\checco_*/r_0\leq \checco r_0
\,,\qquad 
  \|G \|_{\hat D,r_0,s_0}
 \leq \left( 1+4/r_0^2\right)\checco_*
 \leq \checco
  \,,
  \nonumber
  \\ 
  &&\|(1+|J_n-  J_n^*(\hat J )|)
  b (J,\psi_n)\|_{D,r_0/2,s_0}\leq 
\left(4+\frac{34}{r_0^2}\right)\checco_*
\leq \checco
\,,
\\
&&\| |J_n-  J_n^*(\hat J )|
  \partial_{J_n}b (J,\psi_n)\|_{D,r_0/2,s_0}\leq 
  \frac{48}{r_0^2}\checco_*
  \leq \checco\,,
\nonumber
\end{eqnarray}
where 
\begin{equation}\label{tolomeo2}
\checco:=
\left(4+\frac{48}{r_0^2}\right)\checco_*
\end{equation}

\Giu
{\bf Notations}
\nl
{\sl For brevity we introduce the following notations.}
\begin{equation}\label{pgot}
\pgot:=(n,
r_0,
s_0,\b,M)\,. %,\underline R)
\end{equation}
{\sl We say that}
\begin{equation}\label{pgot2}
a\lessdot b \ \ \ {\rm if}\ \ \ \exists \ C=C(\pgot)>0\ \ \ {\rm s.t.}\ \ \  
a\leq C b\,.
\end{equation}
{\sl We also say that, given $F^0$
satisfying} ({\bf A1}) {\sl and} ({\bf A2}),
\begin{equation}\label{pgot3}
a\lessdot_{F^0} b \ \ \ {\rm if}\ \ \ \exists \ C=C(F^0)>0\ \ \ {\rm s.t.}\ \ \  
a\leq C b\,.
\end{equation}
\begin{remark}\label{gianfilippo}
Note that 
$a\lessdot b$ implies
$a\lessdot_{F^0} b.$
 Note also that if {\rm {\bf (A3)}} holds,
 then, by Lemma \ref{sibillini},
 $\pgot$ reduces to $(n,r_0,s_0)$
\end{remark}

Let us assume that
\begin{equation}\label{parigi}
\checco\leq \checco_0=\checco_0(\pgot)\,,
\end{equation}
for a suitable small $\checco_0.$
By \eqref{Lbeta1}, for $\checco_0 $ small enough, 
we can continue  the critical points $x^0_j$ (defined in \eqref{recremisik}),
resp. critical energies $E^0_j,$ of $F^0$ obtaining
critical points
$x_j(\hat J),$
resp. critical energies $E_j(\hat J),$
 of $F(\hat J,\cdot),$
solving the implicit function equation\footnote{To find $x_j=:x^0_j+\chi_j$ we have to solve (for every $\hat J$) the equation $\partial_{\psi_n} F(\hat J,x^0_j+\chi_j)=0.$ Since
$$
\partial_{\psi_n} F(\hat J,x^0_j+\chi_j)=\partial^2_{\psi_n\psi_n}F^0(x^0_j)\chi_j
+O(\chi_j^2)+O(\checco_0)
$$
the equation reduces, by \eqref{Lbeta1}, to find $\chi_j$ solving the fixed point
$\chi_j=O(\chi_j^2)+O(\checco_0)$.
Moreover since $F$ is an analytic function of $\hat J\in \hat D_{r_0}$
the same holds for $\chi_j(\hat J).$
}
\begin{equation}\label{faciolata}
\partial_{\psi_n} F(\hat J, x_j(\hat J))=0
\end{equation}
 and then evaluating
 \begin{equation}\label{faciolata2}
F(\hat J, x_j(\hat J))=:E_j(\hat J)\,,
\end{equation}
respectively.
 Note that
 $x_j(\hat J),$
and $E_j(\hat J)$ are analytic functions of $\hat J\in \hat D_{r_0}.$
By \eqref{Lbeta1}
\begin{equation}\label{october}
\sup_{\hat J\in \hat D_{r_0}} |x_j(\hat J)-x_j^0 |\,,\ 
\sup_{\hat J\in \hat D_{r_0}} |E_j(\hat J)-E_j^0 |
\ \lessdot\, \checco\,.
\end{equation}
Therefore we note that 
$x_j(\hat J)$ and 
$E_j(\hat J)$ maintain the same order 
of $x^0_j$ and $E^0_j.$
In particular, recalling  \ref{timur}
\begin{eqnarray}\label{timur2}
&&E^{(i)}_-(\hat J):=E_i(\hat J)\,,\qquad
E^{(2j-1)}_+(\hat J):=\min\{ E_{2j-2}(\hat J), 
E_{2j}(\hat J)\}\ \ {\rm for}\ \ 1\leq j\leq N\,,
\nonumber
\\
&&E^{(2j)}_+(\hat J):=\min\{ E_{2j_-}(\hat J), 
E_{2j_+}(\hat J)\}\ \ {\rm for}\ \ 1\leq j< N\,,
\qquad
E^{(2N)}_+(\hat J)=E^{(0)}_+(\hat J)=+\infty\,,
\end{eqnarray}
where
$j_\pm$ were defined in \eqref{21stcentury}.

By \eqref{Lbeta1},\eqref{Lbeta2}  for $\checco_0$ small enough we get
  \begin{equation}\label{ladispoli3}
  \inf_{\hat J\in \hat D_{r_0}}
\min_{\psi_n\in\R}
\Big( |\partial_{\psi_n}F(\hat J,\psi_n)|+
|\partial_{\psi_n\psi_n}F(\hat J,\psi_n)|\Big)\geq \frac\b 2\,,\qquad
 \inf_{\hat J\in \hat D_{r_0}}
\min_{i\neq j}|E_j(\hat J)-E_i(\hat J)|\geq \frac\b 2\,.
\end{equation}

Reasoning as above, by \eqref{ladispoli3}, for  
 $\hat J\in \hat D$ (namely $\hat J$ real) and $\checco_0 $ small enough,
we can ``continue'' also the functions  $X^0_{2j}, X^0_{2j-1},$
obtaining
\begin{eqnarray}
X_{2j}(\cdot,\hat J) 
&:& 
\big[E_{2j-1}(\hat J ), E_{2j}(\hat J )\big]
\to \big[x_{2j-1}(\hat J ),x_{2j}(\hat J )\big]\,,
\nonumber
\\
X_{2j-1}(\cdot,\hat J ) 
&:& 
\big[E_{2j-1}(\hat J ), E_{2j-2}(\hat J )\big]
\to \big[x_{2j-2}(\hat J ),x_{2j-1}(\hat J )\big]\,,
\label{foggia}
\end{eqnarray}
solving the implicit function equations
\begin{equation}\label{premiata}
F \big(\hat J, X_i(E,\hat J ) \big)=E\,,
\qquad
X_i\big(F(\hat J,\psi_n ),\hat J \big)=\psi_n\,,
\qquad
\forall\, 1\leq i\leq 2N
\,.
\end{equation}
Note that $X_i$ is increasing, resp. decreasing (as a function of $E$), if 
$i$ is even, resp. odd.
Note also that
\begin{equation}\label{9till5}
\partial_E X_i(E,\hat J )=
1/\partial_{\psi_n} F \big(\hat J, X_i(E,\hat J ) \big)
\end{equation}
and
\begin{eqnarray}\nonumber
&&X_{2j-1}(E_{2j-2}(\hat J ),\hat J )
=x_{2j-2}(\hat J )\,,
\qquad
X_{2j-1}(E_{2j-1}(\hat J ),\hat J )
=x_{2j-1}(\hat J )\,,
\\
&&X_{2j}(E_{2j-1}(\hat J ),\hat J )
=x_{2j-1}(\hat J )\,,
\qquad
X_{2j}(E_{2j}(\hat J ),\hat J )
=x_{2j}(\hat J )\,.
\label{cook}
\end{eqnarray}

%%%%%%%%%%%%%%%%%
%%%%%%%%%%%

%%%%%%%%%%%%%%%
%%%%%%%%%%%%%

\subsection{The integrating transformation}

\begin{proposition}\label{barbabarba}
Let $\Hpend $ be as in \equ{guglielmo}, \equ{mazinga}, \equ{tolomeo}.
There exist   a suitably large constant $C>1,$ which, when \textbf{(A3)} holds, 
depends only on $n,s_0,r_0$, otherwise depend also on $F^0$
(introduced in \eqref{sgravone}) 
\marginpar{ mi sembra che $C$
non dipende solo da $\pgot$ ma pure da $F^0$!!!!} such that, if
 \begin{equation}\label{caviale2}
\checco\leq 1/C 
\,,
\qquad \qquad
0\le \theta\leq 1/C \,,
\end{equation} 
  the following holds.
 There exist:
 \\
i)  disjoint open connected sets 
\begin{equation}\label{sax}
\mathcal C^i (\theta)=\check{\mathcal C}^i (\theta)\times \T^{n-1}\,,
\quad 
0\leq i\leq 2N \,,
\end{equation} 
decreasing w.r.t. $\theta$ and 
satisfying (recall \equ{sontuosa})
\begin{equation}\label{gonatak}
\hat D \times (-R_{0}/2,R_{0}/2)\times\T^n
\subset
\bigcup_{0\leq i\leq 2N}\overline{\mathcal C^i (0)}
\subset
\hat D \times (-R_{0},R_{0})\times\T^n\,,
\end{equation}
 \begin{equation}\label{inthecourtk}
\meas\left(
\Big(\hat D \times (-R_{0}/2,R_{0}/2)\times\T^n\Big)\ \setminus\ 
\bigcup_{0\leq i\leq 2N}\mathcal C^i (\theta)\right)\leq C \theta |\log \theta|\,;
\end{equation}
\\
ii) open connected sets  $\Pgot^i (\theta)$ decreasing w.r.t. $\theta$ 
with
\begin{equation}\label{pratola}
{\rm diam} \,\big(\Pgot^i(0)\big)
\leq 2\big(R_0 +{\rm diam}(\hat D)\big)\,,
\qquad \forall\, 0\leq i\leq 2N\,;
\end{equation}
iii) symplectomorphisms 
 \begin{equation}\label{red}
\Fgot^i :\Pgot^i (0)\times\T^n\ni(P,Q)\to (J,\psi)\in \mathcal C^i (0)
\end{equation}
in $\mathcal G$
such that\footnote{Recall the notation introduced in \eqref{islands}.} 
\begin{equation}\label{ofdelirium}
\Fgot^i (\Pgot^i (\theta)\times\T^n)=
\mathcal C^i (\theta)\ ,\qquad
\check\Fgot^i (\Pgot^i (\theta)\times\T^1)=
\check{\mathcal C}^i (\theta)
\end{equation}
with $\check\Fgot^i $ injective.
\\
Furthermore, $\Fgot^i $ have the following form
\begin{eqnarray}\label{bruford2}
\hat J=\hat P\,,\quad
J_n=\mathtt v^i (P,Q_n )\,,\quad
\hat \psi=\hat Q +\mathtt z^i (P,Q_n ) \,,\quad
\psi_n=\mathtt u^i (P,Q_n )\,,
\qquad \text{for}\ 1\leq i\leq 2N-1\,,
\\
\hat J=\hat P\,,\quad
J_n=\mathtt v^i (P,Q_n )\,,\quad
\hat \psi=\hat Q +\mathtt z^i (P,Q_n ) \,,\quad
\psi_n=Q_n+\mathtt u^i (P,Q_n )\,,
\qquad\,\, \text{for}\ i=0,2N\,,
\label{brufordbis2}
\end{eqnarray}
with $\mathtt v^i ,\mathtt z^i ,\mathtt u^i $,
$2\pi$-periodic in $Q_n,$
$|\mathtt z^i |\leq C \checco$
and, for $i=0,2N,$
$\sup|\partial_{Q_n}\mathtt u^i |<1$;
for $\theta>0,$ $\Fgot^i $ have holomorphic 
extension
\begin{equation}\label{pediatra2bis}
\Fgot^i :\big(\Pgot^i (\theta)\big)_\rho\times\T^n_\s\ \to\ 
D_{r_0}\times\T^n_{s_0}\,,
\end{equation}
with 
\begin{equation}\label{blueeyes}
\rho	= \theta/C\,,\quad \s= 1/C|\log \theta|\,.
\end{equation}
Finally, $\Fgot^i $ ``integrates'' $\Hpend$,
namely\footnote{The function $\mathtt E^{(i)}(P)$ is actually the inverse
of the action function $E\to 
P_n^{(i)}(E,\hat P)$; see \eqref{enrico}
below and Appendix~\ref{gerusalemme}.}:
\begin{equation}\label{red2}
\Hpend\circ\Fgot^i (P,Q)=\Hpend\circ\check{\Fgot}^i (P,Q_n)
=:\mathtt E^{(i)} (P)\,.
\end{equation}
\end{proposition}
Proposition~\ref{barbabarba}
is proved  in \cite{BCaa};
some more details are given in 
 Appendix \ref{gerusalemme}.

\begin{remark}\label{cleopatra} (i) Actually, as standard in the theory of integrable systems, one first introduces the action 
$P_n^{(i)}$ through  line integrals $\oint pdq$ as 
function of energy $E$ and then defines the integrated Hamiltonian  $\mathtt E^{(i)} $ inverting such function; in particular, one has 
\begin{equation}
\label{enrico}
P_n^{(i)}\Big(\mathtt E^{(i)} (P_n, \hat P), \hat P \Big)=P_n\,;
\end{equation}
for more details, see Appendix~\ref{gerusalemme}.

\nl
(ii) Recalling the definition of $E^{(i)}_\pm$  in  
\eqref{timur2} and setting 
\begin{eqnarray}
&&a^{(2j-1)}_-:=0\,,\qquad
a^{(2j-1)}_+:=
P_n^{(2j-1)}\big(E^{(2j-1)}_+(\hat P)-2\theta,\hat P\big)
\,,
\qquad 1\leq j\leq N \,,
\nonumber
\\
&&
a^{(2j)}_-:=
P_n^{(2j)}\big(E^{(2j)}_-(\hat P)+2\theta,\hat P\big)
\,,
\qquad
a^{(2j)}_+:=
P_n^{(2j)}\big(E^{(2j)}_+(\hat P)-2\theta,\hat P\big)
\,,\qquad 1\leq j< N \,,
\nonumber
\\
&&
a^{(0)}_-:=
P_n^{(0)}\big(R_{0}^2-M-2\theta,\hat P\big)
\,,\qquad
a^{(0)}_+:=P_n^{(0)}\big(E^{(2N)}_-(\hat P)+2\theta,\hat P\big)
\,,
\nonumber
\\
&&
a^{(2N)}_-:=P_n^{(2N)}\big(E^{(2N)}_-(\hat P)+2\theta,\hat P\big)
\,,\qquad
a^{(2N)}_+:=
P_n^{(2N)}\big(R_{0}^2-M-2\theta,\hat P\big)\,,
\label{jolebisk}
\end{eqnarray}
one easily recognizes  that 
\begin{equation}\label{playmobilk}
\Pgot^i (\theta):=\Big\{ P=(\hat P,P_n)\ |\ \hat P\in\hat D,\ \  a^{(i)}_-(\hat P,\theta)<P_n
<a^{(i)}_+(\hat P,\theta)
\Big\}\subseteq  \hat D\times \R\subseteq\R^n\,.
\end{equation}

\end{remark}

\subsection{Properties of the actions as functions of the energy}

Next proposition, which is proved in \cite{BCaa}),
contains the fundamental properties 
of $P_n^{(i)}$ (defined in \eqref{enrico}), which will be heavily   exploited in the following. 

Let $P_n^{(i),0}$ be the function in \eqref{enrico}
when $\checco=0,$ namely when
$b,J_n^*,G$ in \eqref{guglielmo} vanish
(see \eqref{sunday} below).

\begin{proposition}\label{glicemiak}
There exist    suitably small, resp. large, constant $\mathtt r>0$,  resp.
$C>1,$ which, when \textbf{(A3)} holds, 
depend only on $n,s_0,r_0$,
otherwise depend also on $F^0$
(introduced in \eqref{sgravone}) 
\marginpar{ mi sembra che $C$
non dipende solo da $\pgot$ ma pure da $F^0$!!!!}
 such that, if
 $\checco\leq 1/C$ then the following holds.
There exist
  real-analytic functions $\phi^{(i)}_\pm(\z,\hat J),$ $\chi^{(i)}_\pm(\z,\hat J)$, 
defined for $|\z|< \mathtt r$, $\hat J\in \hat D_{r_0}$, 
with
\begin{equation}\label{pappagallok}
\sup_{|\z|<\mathtt r,\,\hat J\in \hat D_{r_0}}
\Big(
|\phi^{(i)}_\pm|
+|\chi^{(i)}_\pm|
\Big) < C\,,
\qquad
\sup_{|\z|<\mathtt r,\,\hat J\in \hat D_{r_0/2}}
\Big(
|\partial_{\hat J}\phi^{(i)}_\pm|
+|\partial_{\hat J}\chi^{(i)}_\pm|
\Big) < C \checco\,,
\end{equation}
such that\footnote{For $0\leq i\leq 2N $ except
$i=0,2N $ and the $+$ sign, since 
$E_+^{(0)}(\hat J)=E_+^{(2N )}(\hat J)=+\infty$
 (recall \eqref{timur2}).}
\begin{equation}\label{LEGOk}
P_n^{(i)}\Big(E_\pm^{(i)}(\hat J)\mp\z, \,\hat P\Big)=
\phi^{(i)}_\pm(\z,\hat J) +\z \log \z\, \chi^{(i)}_\pm(\z,\hat J)
 \,,\qquad \text{for}\ \ \ 0<\z<\mathtt r\,,\ \ 
 \hat J\in \hat D\,.
\end{equation}
Moreover
\begin{equation}\label{lamponek}
\chi^{(2j-1)}_-=0\,,\qquad {}
\end{equation}
and\footnote{Recall the definition of 
$F$ in \eqref{guglielmo}.}
\beq{ciofecak}
|\chi^{(i)}_\pm(0,\hat J)|
\geq\frac{1}{4\pi
\sqrt{\|\partial_{\psi_n \psi_n} F \|_{\hat D,r_0,s_0}}}  
\geq 1/C >0
\,.
 \eeq
  Notice that \equ{lamponek} implies that
 $P_n^{(2j-1)}(E,\hat P)$
 has  holomorphic extension 
 on  $\{ |E-E_-^{(2j-1)}(\hat J)|<\mathtt r\}\times \hat D_{r_0}$, as well as
 $P_n^{(2j-1),0}(E)$
 has  holomorphic extension 
 on\footnote{Recall \eqref{timur}.}  $\{ |E-E_-^{(2j-1),0}|<\mathtt r\}.$
\nl 
 Furthermore
\begin{equation}\label{viglianok}
\hat D\times (E^{(i),0}_- +\mathfrak r/4,E^{(i),0}_+-\mathfrak r/4)
\subset \{(\hat P, E)\ {\rm s.t.}\  \hat P\in\hat D, \ E^{(i)}_- (\hat P)<E<E^{(i)}_+(\hat P)\}
\end{equation}
and
\begin{equation}\label{lecavek}
\sup_{\hat D_{r_0}\times \big(E_-^{(i),0}+(-1)^i\mathtt r/2,
E_+^{(i),0}-\mathtt r/2\big)_{\mathtt r/4}}
|P_n^{(i)}(\hat P, E)-P_n^{(i),0}(E)|
< C \checco\,.
\end{equation}
\end{proposition}

\nl
Finally by Lemma \ref{sorellenurzia}
we get
\begin{equation}\label{vallinsu}
\|\partial_{PP} \mathtt E^{(i)} \|_{ \Pgot^i (\theta), \rho}\leq C/\rho\,,\qquad
{\rm for }\ \ \ 0\leq i\leq 2N\,.
\end{equation}

\subsection{The final canonical transformation}
Let us define
\begin{equation}\label{kyoto}
\mathcal C^i_{*} (\theta):= \Phi^{(2{\rm bis})}\big(\mathcal C^i (\theta)\big)\,,\qquad
\check{\mathcal C}^i_{*} (\theta):= \check{\Phi}^{(2{\rm bis})}\big(
\check{\mathcal C}^i (\theta)\big)\,.
\end{equation}
We have that
\begin{equation}\label{kyoto2}
\mathcal C^i_{*} (\theta)=\check{\mathcal C}^i_{*} (\theta)\times\T^{n-1}\,,
\end{equation}
since
$$
\mathcal C^i_{*} (\theta)= \Phi^{(2{\rm bis})}\big(\mathcal C^i (\theta)\big)
\eqby{sax}
\Phi^{(2{\rm bis})}\big(\check{\mathcal C}^i(\theta)\times \T^{n-1}\big)
\eqby{Ventura}
\check{\Phi}^{(2{\rm bis})}\big(\check{\mathcal C}^i(\theta)\big)
\times \T^{n-1}
\eqby{kyoto}
\check{\mathcal C}^i_{*} (\theta)\times \T^{n-1}\,.
$$
Let us also define
\begin{equation}\label{osaka}
\Fgot^i_{*}:=\Phi^{(2{\rm bis})}\circ \Fgot^i\,. 
\end{equation}
By \eqref{bruford2},\eqref{brufordbis2}
we have that
 $\Fgot^i_{*}$ has the form
\begin{eqnarray}\label{bruford3}
\hat J=\hat P\,,\quad
J_n=\mathtt v^i_{*}(P,Q_n )\,,\quad
\hat \psi=\hat Q +\mathtt z^i_{*}(P,Q_n ) \,,\quad
\psi_n=\mathtt u^i(P,Q_n )\,,
\qquad \text{for}\ 1\leq i\leq 2N-1\,,
\\
\hat J=\hat P\,,\quad
J_n=\mathtt v^i_{*}(P,Q_n )\,,\quad
\hat \psi=\hat Q +\mathtt z^i_{*}(P,Q_n ) \,,\quad
\psi_n=Q_n+\mathtt u^i(P,Q_n )\,,
\qquad \text{for}\ i=0,2N\,,
\label{brufordbis3}
\end{eqnarray}
with $\mathtt v^i_{*},\mathtt z^i_{*},\mathtt u^i$,
$2\pi$-periodic in $Q_n,$
$|\mathtt z^i_{*}|\leq C \checco$.

\Giu
Let us  define the  linear symplectic transformation  of the form in \eqref{grillo}
$\Phi_{\rm lin}:\R^{2n}\ni(J'',\psi'')
\mapsto (J',\psi')\in\R^{2n}$ generated by the generating function
$\hat J'\cdot\hat \psi''+
\big( J_n'+\frac{1}{\kappa}(\hat Ak)\cdot \hat J'\big)\psi_n''$ namely (recalling \eqref{centocelle})
\begin{equation}\label{elfi}
J'=U_k J''\,,\quad {\rm with}\quad
\hat J'=\hat J''\,,\quad
J'_n=J_n'' - \frac{1}{\kappa}(\hat A k)\cdot \hat J''
%\ \ {\rm (namely} \ J'=\tilde L J {\rm )}
\,,\quad
\hat\psi'=\hat\psi'' + \frac{\psi_n''}{\kappa} \hat Ak\,,\quad
\psi'_n= \psi_n''\,.
\end{equation}
\begin{remark}\label{incampana}
 Note that such map is only $2\pi\kappa$-periodic in $\psi_n''.$
\end{remark}
Note also that its inverse is
\begin{equation}\label{elfi2}
\Phi_{\rm lin}^{-1}\ \ \ \text{with}\ \ \ 
\hat J''=\hat J'\,,\quad
J''_n=J_n' + \frac{1}{\kappa}(\hat A k)\cdot \hat J'
%\ \ {\rm (namely} \ J'=\tilde L J {\rm )}
\,,\quad
\hat\psi''=\hat\psi' - \frac{\psi_n'}{\kappa} \hat Ak\,,\quad
\psi''_n= \psi_n'
\end{equation}
and that
the operatorial norms of $\Phi_{\rm lin},\Phi_{\rm lin}^{-1}$ are bounded 
by some constant $c(n)>0$
\begin{equation}\label{sfacteria}
\|\Phi_{\rm lin}\|\,, \ \|\Phi_{\rm lin}^{-1}\|
\ \leq \ c(n)\,.
\end{equation}
Recalling the notation in \eqref{islands}
we introduce the volume preserving map
$\check\Phi_{\rm lin}:\R^n\times\T^1.$
Recalling \eqref{spigolo}
we have that $H^*$ defined in 
\eqref{spigola}
satisfies
\begin{equation}\label{spigolo2}
H^*(J'',\psi_n'')=(\tilde H\circ 
\check\Phi_{\rm lin})(J'',\psi_n'')
- \frac{1}{\kappa} \|\proiezione^\perp \hat A^T \hat J''\|^2\,,
\qquad \text{(namely \ \eqref{spigolo})}
\end{equation}
where $\tilde H$ was defined in \eqref{colosseum3}.
Moreover, recalling \eqref{guglielmo},
we have
\begin{equation}\label{spigolo4}
(\tilde H\circ 
\check\Phi_{\rm lin}\circ
\check\Phi^{(2{\rm bis})}
)(J,\psi_n)
=\Hpend(J,\psi_n)
+ \frac{1}{\kappa} \|\proiezione^\perp \hat A^T \hat J\|^2
\end{equation}
and, recalling \eqref{red2}
\begin{equation}\label{spigolo5}
(\tilde H\circ 
\check\Phi_{\rm lin}\circ
\check\Phi^{(2{\rm bis})}
\circ\check{\Fgot}^i)(P,Q_n)
=\mathtt E^{(i)} (P)
+ \frac{1}{\kappa} \|\proiezione^\perp \hat A^T \hat P\|^2
\end{equation}

Recall also the definition of
$
{\mathcal C}^i_{*}(\theta)=\check{\mathcal C}^i_{*}(\theta)\times\T^{n-1}
$
given in  \equ{kyoto}, \equ{kyoto2} 
and of $\Pgot^i(\theta)$ given in 
Proposition~\ref{barbabarba}.

%%%%%%%%%%%%%%%%%%

\nl
Set
\begin{equation}\label{cruel3}
{\mathtt C}^i(\theta):=\check{\mathtt C}^i(\theta)
\times\T^{n-1}\,,\qquad\qquad
\check{\mathtt C}^i(\theta)	:=\check\Phi_{\rm lin}
\big(\check{\mathcal C}^i_{*}(\theta)\big)\times\T^{n-1}
\end{equation}
and
\begin{equation}\label{gimme}
\Cgot^i(\theta):=\Phi^{(1)}
\big({\mathtt C}^i(\theta)\big)\,.
\end{equation}
By \eqref{inthecourtk} we get\footnote{
Note that $U_k\big(\hat D\times (-R_0,R_0)\big)
=U_k \check D=\tilde D^k$ defined in 
 \eqref{formentera}.
}
\begin{equation}\label{inthecourt2}
\meas\Big(
\Big(U_k\big(\hat D\times (-R_0/2,R_0/2)\big)\times\T^n\Big)\ \setminus\ 
\bigcup_{0\leq i\leq 2N}{\mathtt C}^i(\theta)\Big)\lessdot \theta |\log \theta|.
\end{equation}
Recalling \equ{playmobilk}, one sees that one can define the sets $\mathtt B^i(\theta)$ appearing in 2) of Theorem~\ref{porretta} as 
\begin{equation}\label{closetotheedge}
\mathtt B^i(\theta):=
\left\{ \begin{array}{ll} \Pgot^i(\theta) &
\ \ {\rm if} \ \ 1\leq i\leq 2N-1 \\ 
U_k \Pgot^i(\theta) & \ \ {\rm if} \ \ i=0,2N
 \end{array}\right.\,.
\end{equation}
Note that \eqref{pratola2} follows  by \eqref{pratola}
and \eqref{verruca}.

\subsubsection*{1) The oscillatory case: $1\leq i\leq 2N-1$}
\nl
Fix  
 $1\leq i\leq 2N-1$.
Let us define the  symplectomorphism of the form \eqref{grillo}
\begin{eqnarray}
&&\Phi_i\ :\ \mathtt B^i(\theta)\times \T^n\ \to\ 
\mathtt C^i(\theta)\,,
\qquad \text{defined \ as}
\nonumber
\\
&&\hat J'= \hat p
\,,\qquad
J_n'=
\mathtt v^i_{*}(p,q_n )
 - \frac{1}{\kappa}(\hat A k)\cdot \hat p
 \,,
\nonumber
\\
&&
\hat \psi'=
\hat q + \mathtt z^i_{*} (p,q_n ) 
+ \frac{\mathtt u^i(p,q_n )}{\kappa} \hat Ak 
\,,\qquad
\psi_n'=\mathtt u^i(p,q_n )\,,
\label{PFM}
\end{eqnarray}
with 
$\mathtt v^i_{*},\mathtt z^i_{*},\mathtt u^i$
defined in \eqref{bruford3},\eqref{bruford2}.
The fact that it is symplectic can be seen directly by \eqref{PFM}  but also noting that {\sl locally}
\begin{equation}\label{thunder}
\Phi_i =\Phi_{\rm lin}\circ 
\Phi^{{2\rm bis}}\circ
\Fgot^i
=\Phi_{\rm lin}\circ 
\Fgot^i_{*}
\end{equation}
with 
$\Phi_{\rm lin}$
defined in \eqref{elfi},
$\Phi^{{2\rm bis}}$ defined in 
\eqref{valentinesuite2}
 and
$\Fgot^i$ defined in \eqref{red}
(recall also \eqref{bruford2}).
By \eqref{cruel3},\eqref{ofdelirium} and \eqref{closetotheedge},
 $\Phi_i$ is surjective on
$\mathtt C^i(\theta),$
namely
\begin{equation}\label{babele}
\Phi_i\big( \mathtt B^i(\theta)\times \T^n
\big)=
\mathtt C^i(\theta)\,.
\end{equation}
 The injectivity is obvious.  
Note that, by Lemma \ref{schiena},
\eqref{pediatra2bis}, 
\eqref{blueeyes}, \eqref{sfacteria},
 $\Phi_i$
has a holomorphic extension on
\begin{equation}\label{tales}
\Phi_i\ :\ \big(\mathtt B^i(\theta)\big)_{\rho_0}\times \T^n_{\s_0}\ \to \ D_{r_0}\times \T^n_{s_0}\,,
\qquad
\rho_0:=c\frac{\theta}{C} \,, \quad  \ \s_0:=c\frac{1}{C\log \theta}\,,
\end{equation}
where $c$ is a (small) constant 
depending only on $n.$
Applying again Lemma \ref{schiena},
we also prove that
\begin{equation}\label{tales2}
\Phi_i\ :\ \big(\mathtt B^i(\theta)\big)_{\rho_*}\times \T^n_{\s_*}\ \to \ D_{\tilde r}\times \T^n_{\tilde s}\,,
\end{equation}
with $\rho_*,\s_*$, resp. $\tilde r,\tilde s$, defined in 
\eqref{blueeyes3}, resp. \eqref{pelotas},
 taking $\mathtt c$ large enough.
%%%%%%%%%%

\subsubsection*{2) The libration case: $i=0,2N$}
\nl
Fix $i=0,2N$.
Let us define the  symplectomorphism in $\mathcal G$ (recall \eqref{cruel3} and \eqref{closetotheedge})
\begin{eqnarray}
&&\Phi^i\ :\ \mathtt B^i(\theta)\times \T^n\ \to\ 
\mathtt C^i(\theta)\,,
\qquad \text{defined \ as}
\nonumber
\\
&&\hat J'= \hat p
\,,\qquad
J_n'=
\mathtt v^i_{*,k}(U_k^{-1}p,q_n )
 - \frac{1}{\kappa}(\hat A k)\cdot \hat p
 \,,
\nonumber
\\
&&
\hat \psi'=
\hat q +\mathtt  z^i_{*,k} (U_k^{-1}p,q_n ) 
+ \frac{\mathtt u^i(U_k^{-1}p,q_n )}{\kappa} \hat Ak 
\,,\qquad
\psi_n'=q_n+\mathtt u^i(U_k^{-1}p,q_n )\,,
\label{PFM2}
\\
&&
{\rm where,\ recall\ \eqref{centocelle},}\qquad 
U_k^{-1}p=\Big(\hat p, p_n+\frac{1}{\kappa}(\hat A k)\cdot \hat p\Big)\,,
\nonumber
\end{eqnarray}
with
$\mathtt v^i_{*},\mathtt z^i_{*},\mathtt u^i$
defined in \eqref{bruford3},\eqref{bruford2}.
The fact that it is symplectic can be seen directly by \eqref{PFM2} but also noting that {\sl locally}
\begin{equation}\label{thunder2}
\Phi^i=\Phi_{\rm lin}\circ 
\Phi^{{2\rm bis}}\circ
\Fgot^i\circ\Phi_{\rm lin}^{-1}
=
\Phi_{\rm lin}\circ 
\Fgot^i_{*}\circ\Phi_{\rm lin}^{-1}
\,,
\end{equation}
with 
$\Phi_{\rm lin}$
defined in \eqref{elfi},
$\Phi^{{2\rm bis}}$ defined in 
\eqref{valentinesuite2}
 and
$\Fgot^i$ defined in \eqref{red}
(recall also \eqref{bruford2}).
\\
Note that $\Phi^i$ is injective, as it directly follows 
by the fact that so are 
$\Phi^{{2\rm bis}}$ and
$\Fgot^i$,
and also surjective, namely
\eqref{babele} holds;
indeed
\begin{eqnarray*}
&&\Phi^i\big( \mathtt B^i(\theta)\times \T^n
\big)\eqby{Ventura}
\check{\Phi}^i\big( \mathtt B^i(\theta)\times \T^1
\big)\times\T^{n-1}
\eqby{peggylee}
\check{\Phi}_{\rm lin}\big( 
\check{\Phi}^{{2\rm bis}}\big(
\check{\Fgot}^i\big(
\check\Phi_{\rm lin}^{-1}\big(
\mathtt B^i(\theta)\times \T^1
\big)\big)\big)\big)
\times\T^{n-1}
\\
&&
\eqby{closetotheedge}
\check{\Phi}_{\rm lin}\big( 
\check{\Phi}^{{2\rm bis}}\big(
\check{\Fgot}^i\big(
\Pgot^i(\theta)\times \T^1
\big)\big)\big)
\times\T^{n-1}
\eqby{ofdelirium}
\check{\Phi}_{\rm lin}\big( 
\check{\Phi}^{{2\rm bis}}\big(
\mathcal C^i(\theta)
\big)\big)
\times\T^{n-1}
\\
&&
\eqby{kyoto}
\check{\Phi}_{\rm lin}\big( 
\check{\mathcal C}^i(\theta)
\big)
\times\T^{n-1}
\eqby{cruel3}
\check{\mathtt C}^i(\theta)
\times\T^{n-1}
=\mathtt C^i(\theta)
\end{eqnarray*} 
\\
Reasoning as in the case
$1\leq i<2N,$ we get that
 $\Phi^i$
has a holomorphic extension as in  \equ{tales}.
\eqref{tales2} holds as well.

\subsection{Conclusion of the
proof of part one of the Structure Theorem}

\nl
We set
\begin{equation}\label{mammamia}
\Psi^i:=\Phi^{(1)}\circ\Phi_i\,,
\end{equation}
where $\Phi^{(1)}$
was defined in \eqref{talktothewind}; therefore
\eqref{blueeyes3} holds by \eqref{tales2}
and \eqref{fare2}.

\nl
Then \eqref{fujiyama}
follows by \eqref{inthecourt2}
recalling \eqref{LPDbis}, \eqref{calimero2}
and since $\Phi^{(1)}$
preserve volume being symplectic.

\nl
\eqref{tokyo} follows 
by \eqref{babele}.

\nl
Recalling \eqref{colosseum3}
and \eqref{spigolo5} we have that
$
h^{(i)}:=\cH\circ\Phi^{(1)}\circ \Phi_i
$
can be written as
\begin{equation}\label{matteo2}
h^{(i)}(p)=
\left\{ \begin{array}{ll} \mathtt E^{(i)}(p)+
\hat h_k(\hat p) &
\ \ {\rm if} \ \ 1\leq i< 2N\,, \\ 
 (\mathtt E^{(i)}+\hat h_k)(U_k^{-1}p) 
 & \ \ {\rm if} \ \ i=0,2N\,.
 \end{array}\right.
\end{equation}
where
\begin{equation}\label{tyr} 
\hat h_k(\hat p):=
\kappa^{-1}\|\proiezione^\perp \hat A_k^T p\|^2\,.
\end{equation}
 \eqref{matteo}
follows.
Finally \eqref{vallinsu2} follows by 
\eqref{vallinsu},
\eqref{matteo2}, \eqref{tyr},
\eqref{scimmia}, \eqref{verruca}.

\nl
 This concludes the proof of part 
one of the Structure Theorem.

\section{Proof of Part II of the Structure Theorem}\label{lupin2}

\nl
A crucial fact is that the integrable Hamiltonian
$h^{(i)}$ defined in \eqref{matteo2}
twists, namely \eqref{gadamer} (toghether with the measure estimate \eqref{gibboso2}) holds.
This will be a direct consequence of the following

\begin{proposition}\label{frey}
Under the assumptions of Theorem \ref{porretta}
(in particular \eqref{gaetanga}),
for any $\mu$ satisfying \eqref{gaetanga2},  there exists a subset\footnote{Decreasing w.r.t. $\mu.$}  $\widetilde \Pgot^i (\mu) \subseteq\Pgot^i (0),$
satisfying\footnote{$\mathtt c$ defined in
Theorem \ref{porretta}.} 
\begin{equation}\label{gibboso}
\meas(\Pgot^i (0)\bks \widetilde \Pgot^i (\mu))\le 
\mathtt c \|k\|^{4n} \mu^{1/\mathtt c}\,,
\end{equation}
such that
\begin{equation}\label{gea}
\left| {\rm det}\left[\partial_{PP}
\left(
\mathtt E ^{(i)}(P)+\hat h_k(\hat P)
\right)
    \right]  \right|
> \mu\,,\qquad
\forall\, \ 0\leq i\leq 2N \,,\qquad \forall \ P\in  \widetilde \Pgot^i (\mu)\ .
\end{equation}
\end{proposition}

\noindent
\textbf{Proof of \eqref{gadamer}}. 
One sees that, in analogy to 
\eqref{closetotheedge},
 one can define the sets 
$\tilde{\mathtt B}^i (\mu)$ appearing in 3) of Theorem~\ref{porretta} as 
\begin{equation}\label{closetotheedge3}
\tilde{\mathtt B}^i (\mu):=
\left\{ \begin{array}{ll} \tilde{\Pgot}^i (\mu) &
\ \ {\rm if} \ \ 1\leq i\leq 2N-1\,, \\ 
U_k \tilde{\Pgot}^i (\mu) & \ \ {\rm if} \ \ i=0,2N\,.
 \end{array}\right.
\end{equation}
\nl
Recalling  \equ{closetotheedge}and \eqref{verruca} , by \equ{gibboso}, we get \eqref{gibboso2}.
Finally \eqref{gadamer} follows by \eqref{gea}, Lemma 
\ref{cippa} and noting that det$\, U_k=1$ by \eqref{centocelle}.

\noindent
\textbf{Proof of Proposition \ref{frey}}.
\\
First we note that
 \begin{equation}\label{lovebuzz}
\det\big( \partial_{\hat P \hat P}\hat h_k \big)=
2^{n-1} \kappa^{-n}\,.
\end{equation}
Indeed by \eqref{farfalla}
$$
\hat h_k(\hat P)
+ (P_n)^2
=\kappa^{-1} \|L_k P\|^2\,,
$$
(with $L_k$ defined in \eqref{LPD}) and,
by Lemma \ref{cippa}  and \eqref{diego}, we get
\begin{equation}\label{strunz}
\mathtt d_k:=
\det\Big( \partial_{PP}\big(\hat h_k(\hat P)
+ (P_n)^2 \big) \Big)=2^n \kappa^{-n}
\end{equation}
and, therefore, \eqref{lovebuzz} follows.

\subsubsection*{The ``twist'' determinant as a function of the energy}

\nl
Let us fix $0\leq i\leq 2N$.
Consider the analytic 
function
\begin{equation}\label{nghe}
d^i_k(E,\hat P):=
{\rm det}\left[\partial_{PP} \hat h_k (\hat P)+ 
\partial_{PP}\mathtt E^{(i)}
   \big(\hat P, P_n^{(i)}(\hat P, E)\big) \right] \,,
\end{equation}
where,  $P_n^{(i)}(E,\hat P)$ is, by 
definition\footnote{
This map can be obviously explicitely constructed, see 
subsection \ref{hocuspocus} below.},
the inverse map of $\mathtt E^{(i)}(\hat P, P_n)$
(recall \eqref{enrico}).

\nl
For brevity we will often omit to write the indexes $_k$ and/or $^{(i)}.$

\subsubsection*{The case $i$ odd; close to a maximum
of the potential}

Instead of use the variable $E$
we use, in the case of odd $i=2j-1$,
\begin{equation}\label{chiappa}
\z:=E_+^{(2j-1)}(\hat P)-E\,,
\end{equation}
where $E_+^{(2j-1)}(\hat P)$
was defined in \eqref{timur2}.
In the variable $\z$, (note that $\partial_E$ is equal to
$\partial_\z$, up to sign), recalling \eqref{LEGOk},
\begin{equation}\label{zabaione}
P_n(E,\hat P)=P_n\big(E_+^{(2j-1)}(\hat P)-\z,\hat P\big)=1\oplus \z\log\z \,,
\end{equation}
where by $g=g_1\oplus g_2$
we mean that there exist two  functions
$\f_1(\z,\hat P),\f_2(\z,\hat P)$, analytic in a complex neighborhood
of zero times $\hat D$ (depending only on $\pgot$ defined in \eqref{pgot}), such that
$g(\z,\hat P)=g_1(\z) \f_1(\z,\hat P)+g_2(\z) \f_2(\z,\hat P).$
%ATTENZIONE ALLA DIFFERENZA TRA
%$$\partial_{\hat P}
%P_n(E,\hat P)\qquad
%\text{e}\quad
%\partial_{\hat P}[P_n\big(E_+^{(2j-1)}(\hat P)-\z,\hat P\big)]
%=\partial_{E}P_n\big(E_+^{(2j-1)}(\hat P)-\z,\hat P\big)\partial_{\hat P}E_+^{(2j-1)}(\hat P)
%+\partial_{\hat P}P_n\big(E_+^{(2j-1)}(\hat P)-\z,\hat P\big)
%$$  
%MI SA CHE IO QUI SOTTO HO MESSO LA SECONDA INVECE DELLA PRIMA!!!!
We get (for $i,j=1,\ldots,n-1$)
\begin{eqnarray}\label{gorgo}
&&\partial_E P_n=1\oplus \log\z \,,\quad
\partial_{\hat P_i}P_n=\checco(1\oplus \z\log\z )\,,\quad
\partial_{EE} P_n= \log\z \oplus  \z^{-1}\,,
\nonumber
\\
&&\partial_{E\hat P_i} P_n=\checco(1\oplus \log\z )\,,\quad
\partial_{\hat P_i\hat P_j}P_n=\checco(1\oplus \z\log\z )\,.
\end{eqnarray}
\\
Recalling \eqref{enrico}, by the chain rule we get\footnote{By $\partial_{\hat P}$ we mean the row vector
$(\partial_{P_1},\ldots,\partial_{P_{n-1}})$, then $\partial_{\hat P}^T$ is a column vector.} 
\begin{eqnarray}
&&\partial_{P_n} \mathtt E  
=
\frac{1}{\partial_E P_n }\,,\qquad
\partial_{\hat P} \mathtt E  
=
-\frac{\partial_{\hat P} P_n }{\partial_E P_n }\,,
\qquad
\partial_{P_n P_n} \mathtt E  
=
-\frac{\partial_{EE} P_n }{(\partial_E P_n )^3}\,,
\nonumber
\\
&&\partial_{P_n \hat P} \mathtt E  
=
\frac{\partial_{EE} P_n  \partial_{\hat P}P_n }{(\partial_E P_n )^3}-
\frac{\partial_{E\hat P} P_n  }{(\partial_E P_n )^2}\ \in\ \R^{n-1}
\,,
\nonumber
\\
&&\partial_{\hat P \hat P} \mathtt E  
=
-\frac{\partial_{\hat P\hat P} P_n }{\partial_E P_n }
+2\frac{\partial_{\hat P}^TP_n\ \partial_{\hat P}(\partial_{E} P_n) }{(\partial_E P_n )^2}
-\frac{\partial_{EE} P_n \  \partial_{\hat P}^T P_n\  \partial_{\hat P}P_n}{(\partial_E P_n )^3}
\ \in\ {\rm Mat}_{(n-1)\times(n-1)}\,,
\label{daitarn3}
\end{eqnarray}
where $\mathtt E $ and  $P_n $ are evaluated in  $\big(P_n (E,\hat P),\hat P\big)$
and $(E,\hat P),$ respectively\footnote{
Or, which is equivalent, in $P$ and
$\big(\mathtt E (P),\hat P\big),$
respectively.}. 
\\
By \eqref{gorgo} and \eqref{daitarn3} 
we get  (for $i,j=1,\ldots,n-1$)
\begin{eqnarray}
\z (\partial_E P_n )^3\partial_{P_n P_n} \mathtt E  
&=&
-\z \partial_{EE} P_n 
=
1\oplus \z \log \z 
\,,
\nonumber
\\
\z (\partial_E P_n )^3\partial_{P_n \hat P_i} \mathtt E  
&=&
\z \partial_{EE} P_n  \partial_{\hat P_i}P_n -
\z  \partial_{E\hat P_i} P_n  \partial_E P_n
=
\checco(1\oplus \z \log \z \oplus \z \log^2 \z )
\,,
\nonumber
\\
\z (\partial_E P_n )^3\partial_{\hat P_i \hat P_j} \mathtt E  
&=&
-\z (\partial_E P_n )^2 \partial_{\hat P_i\hat P_j} P_n 
+2\z  \partial_E P_n \partial_{\hat P_i}  P_n\ \partial_{\hat P_j E} P_n 
-\z \partial_{EE} P_n \  \partial_{\hat P_i} P_n\  \partial_{\hat P_j}P_n
\nonumber
\\
&=&
\checco(1\oplus \z \log \z \oplus \z \log^2 \z 
\oplus \z^2\log^3 \z )
\,.
\label{fausti}
\end{eqnarray}

\subsubsection*{The ``rescaled'' determinant $\mathtt f$}

Recalling that, by \eqref{chiappa} we have
$$
E=E_+^{(2j-1)}(\hat P)-\z \,,
$$
we set
\begin{equation}\label{mattoni}
\mathtt f(\z ,\hat P)=\mathtt f^{(i)}_k(\z ,\hat P):=
\z^{n} 
\big(\partial_E P_n(E_+^{(2j-1)}(\hat P)-\z ,\hat P)\big)^{3n}
d^i_k(E_+^{(2j-1)}(\hat P)-\z ,\hat P)
\,.
\end{equation}
We will omit the dependence on $k,i$.
Note that 
\begin{eqnarray*}
d^i_k
&=&
{\rm det}\left(\partial_{PP} \hat h_k + 
\partial_{PP}\mathtt E^{(i)}_k\right)
=
  \text{det}\left(
  \begin{array}{cc}
   \partial_{\hat P\hat P} h +  \partial_{\hat P\hat P}\mathtt E 
  &
    \quad\partial_{\hat P}^T(\partial_{P_n}\mathtt E ) \\ 
 \partial_{\hat P}(\partial_{P_n}\mathtt E )    
    &
  \partial_{P_n P_n}\mathtt E    
  \\ 
  \end{array}
  \right)
  \\
&=&
   \partial_{P_n P_n}\mathtt E
   \ \text{det} (\partial_{\hat P\hat P} h +  \partial_{\hat P\hat P}\mathtt E)
   +
    \text{det}\left(
  \begin{array}{cc}
   \partial_{\hat P\hat P} h +  \partial_{\hat P\hat P}\mathtt E 
  &
    \quad\partial_{\hat P}^T(\partial_{P_n}\mathtt E ) \\ 
 \partial_{\hat P}(\partial_{P_n}\mathtt E )    
    &
 0    
  \\ 
  \end{array}
  \right)
\end{eqnarray*}
developing the determinat w.r.t. the last 
column.
By \eqref{fausti}  we get
\begin{eqnarray*}
&&\z^{n} (\partial_E P_n)^{3n} \partial_{P_n P_n}\mathtt E
   \ \text{det} (\partial_{\hat P\hat P} h +  \partial_{\hat P\hat P}\mathtt E)
\\
&&= 
\z^{n-1} \sum_{\ell =0}^{3n-3}\f_\ell^{(1)} (\z ,\hat P)
\log^\ell  \z \ +\ 
\checco \sum_{\ell=0}^{3n-5}
 \f_{\ell }^{(2)} (\z ,\hat P)  \log^\ell  \z
\ +\ R^{(1)} 
\,,
\end{eqnarray*}
where 
\begin{equation}\label{eagles}
\inf_{\hat P\in\hat D}|\f_{3n-3}^{(1)}(0,\hat P)|\geq c(\pgot) \mathtt d_k>0\,,
\end{equation}
with $\mathtt d=\mathtt d_k$ defined in \eqref{strunz}
and
$$
R^{(1)}=O(\z^n)=\z^n\sum_{\ell =0}^{3n-2}
\f_\ell^{(3)} (\z,\hat P)
\log^\ell  \z\,,
$$
for suitable $\f_{\ell }^{(1)}, \f_{\ell }^{(2)}, 
\f_{\ell }^{(3)},$ analytic and uniformly
bounded\footnote{Here and in the following of this section by ``uniformly bounded''
we mean ``bounded by a constant depending only on  $\pgot$ in }  in a uniform neighborhood of zero (and $\hat P\in \hat D$). 
 
\begin{remark}\label{controllare}
 We stress that, as a consequence of \eqref{pappagallok},
 all the estimates of this section are uniform
 in $\hat P\in \hat D.$
\end{remark}

\nl
Moreover
\begin{eqnarray*}
&&\z^{n} (\partial_E P_n)^{3n}\
  \text{det}\left(
  \begin{array}{cc}
   \partial_{\hat P\hat P} h +  \partial_{\hat P\hat P}\mathtt E 
  &
    \quad\partial_{\hat P}^T(\partial_{P_n}\mathtt E ) \\ 
 \partial_{\hat P}(\partial_{P_n}\mathtt E )    
    &
 0    
  \\ 
  \end{array}
  \right)
\\
  &=&\checco^2 (1\oplus \z  \log \z  \oplus \z  \log^2 \z )^2
(1\oplus \z  \log \z  \oplus \z  \log^2 \z \oplus 
\z\log^3 \z )^{n-2}
\\
&=&  
\checco^2 \sum_{\ell=0}^{3n-4}
 \f_{\ell }^{(4)} (\z ,\hat P)  \log^\ell  \z
\ +\ R^{(2)} 
\,,
\end{eqnarray*}
where 
$$
R^{(2)}=O(\z^n)=\z^n\sum_{\ell =0}^{3n-2}
\f_\ell^{(5)} (\z,\hat P)
\log^\ell  \z\,,
$$
for suitable
$\f_{\ell }^{(4)}, 
\f_{\ell }^{(5)},$
 analytic and uniformly bounded in a uniform neighborhood of zero (and $\hat P\in \hat D$).

\nl
Therefore
\begin{equation}\label{takeiteasy}
\mathtt f=\mathtt f^{(i)}_k=
\z^{n-1} \sum_{\ell =0}^{3n-3}\f_\ell (\z ,\hat P)
\log^\ell  \z \ +\ 
\checco \sum_{\ell=0}^{3n-4}
 \chi_{\ell } (\z ,\hat P)  \log^\ell  \z
\ +\ R
\end{equation}
where, by \eqref{eagles},
\begin{equation}\label{eagles2}
\inf_{\hat P\in\hat D}|\f_{3n-3}(0,\hat P)|\geq c(\pgot) \mathtt d_k>0\,,
\end{equation}
and
\begin{equation}\label{croccante}
R=O(\z^n)=\z^n\sum_{\ell =0}^{3n-2}
\psi_\ell (\z,\hat P)
\log^\ell  \z\,,
\end{equation}
for suitable $\f_{\ell }, \chi_{\ell }, \psi_\ell$ analytic and uniformly bounded in a uniform neighborhood of zero (and $\hat P\in \hat D$). 

\subsubsection*{A class of useful  linear operators}

We now consider the linear operator
$$
L:=\z \partial_\z
$$
and, recursively, $L^{\ell}:=L\circ L^{\ell-1}.$
We have
\begin{eqnarray}
L \left(\z^h
 \sum_{\ell=0}^{\ell_0} g_{\ell }(\z ,\hat P)  \log^\ell  \z
  \right)
  &=&
\z^h \sum_{\ell=0}^{\ell_0} \tilde g_{\ell }(\z,\hat P)  \log^\ell  \z\,,
 \qquad
 \text{with}\ \  \tilde g_{\ell_0 }(0,\hat P)=h g_{\ell_0 }(0,\hat P)\,,
 \label{mostaccioli}
 \\  
L^{\ell_0} \left(
 \sum_{\ell=0}^{\ell_0} g_{\ell }(\z,\hat P)  \log^\ell  \z
  \right)
  &=&
  \ell_0!\, g_{\ell_0}(\z,\hat P)\, +\, 
\z  \sum_{\ell=0}^{\ell_0}  \tilde g_{\ell }(\z,\hat P)  \log^\ell  \z\,,
\label{mostaccioli2}
\\  
L^{\ell_1} \left(
 \sum_{\ell=0}^{\ell_0}  g_{\ell }(\z,\hat P)  \log^\ell  \z
  \right)
  &=&
\z  \sum_{\ell=0}^{\ell_0} \tilde g_{\ell }(\z,\hat P)  \log^\ell  \z\,,\qquad \text{when}\ \ \ell_1>\ell_0\,,
\label{mostaccioli3}
\end{eqnarray}
for (different) suitable $\tilde g_{\ell }.$
Moreover
\begin{eqnarray}
&&\partial_\z 
\left(\z^h
\sum_{\ell=0}^{\ell_0}  g_{\ell }(\z,\hat P)  \log^\ell  \z
  \right)\ =
  \label{mostaccioli4}
 \\ 
  &&
  h   g_{\ell_0 }(\z,\hat P) \z^{h-1} \log^{\ell_0}\ +\
  \z^{h-1} \sum_{\ell=0}^{\ell_0-1}\Big(h  g_{\ell }(\z,\hat P)+(\ell+1)
    g_{\ell+1 }(\z,\hat P) \Big)\log^\ell  \z
+\ \z^h
\sum_{\ell=0}^{\ell_0} \tilde g_{\ell }(\z,\hat P)  \log^\ell  \z\,.
\nonumber
\end{eqnarray}
In particular by \eqref{mostaccioli},
\eqref{mostaccioli3},\eqref{mostaccioli4}
we get
\begin{equation}\label{mostaccioli5}
(\partial_\z \circ L^{\ell_1}) \left(
 \sum_{\ell=0}^{\ell_0} 
   g_{\ell }(\z,\hat P) \z^h  \log^\ell  \z
  \right)
  =
\sum_{\ell=0}^{\ell_0} 
  \tilde g_{\ell }(\z,\hat P) \z^h  \log^\ell  \z\,,\qquad \text{when}\ \ \ell_1>\ell_0\,,
\end{equation}
for suitable $ \tilde g_{\ell }.$

\nl
Introduce the linear differential operator (w.r.t. $\z$) of order $3 n^2 - 3n:$
\begin{equation}\label{plenty}
\mathcal L:=L^{3n-3} (\partial_\z \circ L^{3n-3})^{n-1}\,.
\end{equation}

\subsubsection*{Non-degeneracy of the derivatives of $\mathtt f$}

Let us decompose $\mathtt f$
in \eqref{takeiteasy} as
$$
\mathtt f=
\z^{n-1} \f_{3n-3} (0 ,\hat P)
\log^{3n-3}  \z
+\tilde{\mathtt f} +\tilde R\,,
$$
with
$$
\tilde{\mathtt f}
:=
\z^{n-1} \sum_{\ell =0}^{3n-4}\f_\ell (\z ,\hat P)
\log^\ell  \z \ +\ 
\checco \sum_{\ell=0}^{3n-4}
 \chi_{\ell } (\z ,\hat P)  \log^\ell  \z\,,
 $$
 and\footnote{Note that the function
 $\frac{\f_{3n-3} (\z ,\hat P)
-\f_{3n-3} (0 ,\hat P)}{\z}$ is analytic.}
 $$
\tilde R:=\z^{n} \left(
\frac{\f_{3n-3} (\z ,\hat P)
-\f_{3n-3} (0 ,\hat P)}{\z}
 \right)
\log^{3n-3}  \z +R
$$
We want to evaluate $\mathcal L \mathtt f.$
We have
$$
\mathcal L\, \Big(\z^{n-1} \f_{3n-3} (0 ,\hat P)
\log^{3n-3}  \z\Big)= 
(3n-3)! \big((n-1)!\big)^{3n-2} \f_{3n-3}(0,\hat P)
+
\z \sum_{\ell =0}^{3n-2}
\tilde \f_\ell (\z,\hat P)
\log^\ell  \z\,,
$$
for suitable $\tilde \f_\ell.$
By \eqref{mostaccioli3} and 
 \eqref{mostaccioli5}
we get
$$
\mathcal L \, \tilde{\mathtt f}\, =\,
\z \sum_{\ell=0}^{3n-4}
\tilde \chi_{\ell } (\z ,\hat P)  \log^\ell  \z\,,
$$
for suitable $\tilde \chi_{\ell }.$
Finally, by \eqref{mostaccioli} and 
\eqref{mostaccioli4},
 we have that
$$
\mathcal L \tilde R=
\z\sum_{\ell =0}^{3n-2}
\tilde\psi_\ell (\z,\hat P)
\log^\ell  \z\,,
$$
(with $R$ defined in \eqref{croccante})
for suitable $\tilde\psi_\ell$.
 Recollecting we get
 \begin{equation}\label{caciottone}
(\mathcal L \mathtt f)(\z,\hat P)
=
\mathtt c_k(\hat P)
\ + \ \z \sum_{\ell =0}^{3n-2}
\tilde{\mathtt f}_\ell (\z,\hat P)
\log^\ell  \z\,,
\end{equation}
for suitable 
$\tilde{\mathtt f}_\ell (\z,\hat P)$
analytic and uniformly bounded in a uniform neighborhood of zero (and $\hat P\in \hat D$)
and
where
\begin{equation}\label{sonole3}
\mathtt c_k(\hat P):=
(3n-3)! \big((n-1)!\big)^{3n-2} \f_{3n-3}(0,\hat P)\,,
\end{equation}
satisfies, by \eqref{eagles2},
\begin{equation}\label{eagles3}
\inf_{\hat P\in\hat D}|\mathtt c_k(\hat P)|\geq c(\pgot) \mathtt d_k>0\,.
\end{equation}
By \eqref{caciottone}
\begin{equation}\label{eagles4}
\inf_{0<\z\leq \z_0}
\inf_{\hat P\in\hat D}\big|(\mathcal L \mathtt f)(\z,\hat P)\big|\ \geq \ c(\pgot) \mathtt d_k/2>0\,,
\qquad {\rm for\ a\ suitable}\ \ \ 
\z_0=\z_0(\pgot)>0\,.
\end{equation}
By \eqref{plenty} we have that
\begin{equation}\label{ciaone}
(\mathcal L \mathtt f)(\z,\hat P)
=
\sum_{d=0}^{m_n}
a_d(\z)
\partial_\z^d \mathtt f(\z,\hat P)\,,
\qquad
\text{where}\ \ m_n:=3 n^2 - 3n\,,
\end{equation}
for suitable polynomials $a_d(\z).$
Then by \eqref{eagles4}, taking in case $\z_0$ smaller, 
we get
\begin{equation}\label{loosing}
\inf_{0<\z\leq \z_0}
\inf_{\hat P\in\hat D}
\max_{ 1\leq d\leq m_n}|\partial_\z^d \mathtt f(\z,\hat P)|/d!
\geq
c(\pgot) \mathtt d_k>0
\end{equation}
(with a smaller $c(\pgot)$).

\subsubsection*{The measure of sublevels of $\mathtt f$}

We need the following result, which is proved in the appendix

\begin{lemma}\label{girotondi}
Let $f\in C^{m+1}([a,b])$ 
and  assume that for some $m\geq 1$
\begin{equation}\label{macine}
\min_{x\in[a,b]}\max_{ 1\leq d\leq m}|\partial_x^{d}f(x)|/d!=:\xi_m>0\,.
\end{equation}
Then for $0<\mu <1$
$$
{\rm meas} \big( \{ x\in[a,b]\ \ :\ \ |f(x)|\leq \mu  \} \big)
\leq \frac{m(M+1)(b-a+2\mu  ^{1/m+1})}{\xi_m} \mu  ^{\frac{1}{m(m+1)}}\,,
$$
with
$M:=\max_{x\in[a,b],\ 2\leq d\leq m+1}|\partial_x^{d}f(x)|/d!$
\end{lemma}
We apply Lemma \ref{girotondi} with
$$
f\rightsquigarrow \mathtt f\,,\quad
 m\rightsquigarrow m_n\,,\quad
  a\rightsquigarrow \z_1\,,\quad
   b\rightsquigarrow \z_0\,,\quad
    \xi_m\rightsquigarrow c(\pgot)\mathtt d_k  \ \ \text{(in\ \eqref{loosing})}\,,\quad
     x\rightsquigarrow \z\,,\quad
$$
with $\z_1$ to be chosen later.
By \eqref{takeiteasy} we have that\footnote{Denoting,
as usual, the $\z_1/2$-complex-neighborhood of the real interval
$[\z_1,\z_0]$ by $[\z_1,\z_0]_{\z_1/2}$.}
$$
\sup_{\hat P\in \hat D}\sup_{[\z_1,\z_0]_{\z_1/2}} |\mathtt f(\z,\hat P)|\lessdot
|\log^{2n-2} \z_1|\leq 1/\z_1
$$
for $\z_1$ small enough.
Then, by Cauchy estimates, we get that $M$ in Lemma \ref{girotondi} 
satisfies
$$
M\lessdot 1/\z_1^{m_n+2}\,.
$$
By Lemma \ref{girotondi} we get,
for $0<\mu<1,$
$$
{\rm meas} \big( \{ \z\in[\z_1,\z_0]\ \ :\ \ 
|\mathtt f(\z,\hat P)|\leq \mu  \} \big)
\lessdot \mu  ^{\frac{1}{m_n(m_n+1)}}/\mathtt d_k 
\z_1^{m_n+2}\,,
$$
for every\footnote{Note that the hidden constant 
in $\lessdot$ is independent of $\hat P\in \hat D.$
} $\hat P\in \hat D.$
We can optimize the choice of $\z_1$ taking$$
\z_1=
\mu  ^{\frac{1}{m_n(m_n+1)}}/ 
\z_1^{m_n+2}\,,\qquad \text{namely}
\ \ \  \z_1:= 
\mu  ^{\frac{1}{m_n(m_n+1)(m_n+3)}}\,,
$$
so that
\begin{equation*}
{\rm meas} \big( \{ \z\in(0,\z_0]\ \ :\ \ 
|\mathtt f_k(\z,\hat P)|\leq \mu  \} \big)
\lessdot \mu  ^{\frac{1}{m_n(m_n+1)(m_n+3)}}/\mathtt d_k\,,
\end{equation*}
for every $\hat P\in \hat D.$
Since by the first equality in \eqref{gorgo}
$$
\z^{n} \big(\partial_E P_n
(E_+^{(2j-1)}(\hat P)-\z ,\hat P)\big)^{3n}\leq 1
$$
for $\z_0$ small,
we get, recalling \eqref{mattoni},
\begin{equation}\label{missionbell0}
{\rm meas} \big( \{ \z\in(0,\z_0]\ \ :\ \ 
|d^i_k(E_+^{(2j-1)}(\hat P)-\z ,\hat P)|\leq \mu  \} \big)
\lessdot \mu  ^{\frac{1}{m_n(m_n+1)(m_n+3)}}/\mathtt d_k\,,
\end{equation}
for every $\hat P\in \hat D,$
where $d_k^i$ was defined in \eqref{nghe}.
Recalling \eqref{chiappa}, we get \begin{equation}\label{missionbell}
{\rm meas} \big( \{ E\in\big[E_+^{(2j-1)}(\hat P)-\z_0,E_+^{(2j-1)}(\hat P)\big)\ \ :\ \ 
|d^i_k(E,\hat P)|\leq \mu  \} \big)
\lessdot \mu  ^{\frac{1}{m_n(m_n+1)(m_n+3)}}/\mathtt d_k\,,
\end{equation}
for every 
$\hat P\in \hat D$.

\subsubsection*{Conslusion of the proof in the case
$i$ odd, close to maxima of the potential}

Recalling \eqref{nghe} and \eqref{red2}, we have
\begin{equation}\label{nocciola}
{\rm det}\left[\partial_{PP} \big(\hat h_k (\hat P)+ 
\mathtt E^{(i)}_k(P)\big)
  \right] 
   =
   d^i_k(\mathtt E^{(i)}_k(P),\hat P)\,.
\end{equation}
Recalling that, by \eqref{jolebisk}, we have
$
E_+^{(2j-1)}(\hat P)=\mathtt E^{(2j-1)}_k(\hat P, a_+^{(2j-1)}(\hat P))\,,
$ we set
\begin{equation}\label{pistacchio}
\mathtt r_0(\hat P):=
a_+^{(2j-1)}(\hat P) - 
P_n\Big( E_+^{(2j-1)}(\hat P)-\z_0 ,\hat P \Big)
=
P_n\Big( E_+^{(2j-1)}(\hat P),\hat P \Big)
 - 
P_n\Big( E_+^{(2j-1)}(\hat P)-\z_0 ,\hat P \Big)
\,.
\end{equation}
Note that by \eqref{pistacchio}
and \eqref{vana}
(which relies on Proposition
\ref{glicemiak}) there exists $\mathtt r_0>0$,
which, when \textbf{(A3)} holds, 
depends only on $n,s_0,r_0$, otherwise depend also on $F^0$,
such that 
\begin{equation}\label{peppapig}
\mathtt r_0(\hat P)\geq \mathtt r_0\,,\qquad \hat P\in\hat D\,.
\end{equation}
For every $\hat P\in\hat D,$
we now consider the change of variable 
$P_n=P_n (E,\hat P).$
By \eqref{missionbell} and \eqref{nocciola},
and noting that (recall \eqref{zabaione})
\begin{equation}\label{stracciatella}
0<\partial_E P_n(E,\hat P)
\lessdot
\big |\log (E_+^{(2j-1)}(\hat P)-E)|
\,,
\end{equation}
 we 
get\footnote{  Let 
$\mathcal E(\hat P)\subseteq \big[E_+^{(2j-1)}(\hat P)-\z_0,E_+^{(2j-1)}(\hat P)\big)$ 
and 
$\mathcal P(\hat P)\subseteq 
\big[a_+^{(2j-1)}(\hat P)-\mathtt r_0,
a_+^{(2j-1)}(\hat P)
\big)$ be
the sets whose measures are estimated in \eqref{missionbell}
and \eqref{missionbell2}; so that 
$\mathcal P(\hat P)=P_n(\mathcal E(\hat P),\hat P).$
Then, denoting by $\mathcal Z(\hat P)\subseteq(0,\z_0]$ the 
set whose measure 
$\mu_0:={\rm meas}(\mathcal Z(\hat P))$
 is estimated in \eqref{missionbell0}, we have
\begin{eqnarray*}
&&
{\rm meas}(\mathcal P(\hat P))=\int_{\mathcal P(\hat P)}dP_n
=\int_{\mathcal E(\hat P)} \partial_E P_n(E,\hat P)dE
\stackrel{\eqref{stracciatella}}\lessdot
\int_{\mathcal E(\hat P)} 
 |\log (E_+^{(2j-1)}(\hat P)-E)|dE
\\
&&
\lessdot
\int_{\mathcal Z(\hat P)}  |\log\z|d\z
=\int_{\mathcal Z(\hat P)\cap(0,\mu_0]} |\log\z|d\z
+\int_{\mathcal Z(\hat P)\cap(\mu_0,\z_0]} |\log\z|d\z
\leq \mu_0|\log\mu_0|+|\log\z_0| \mu_0
\lessdot \mu_0|\log\mu_0|
\,.
\end{eqnarray*}
Therefore by \eqref{missionbell0} 
$$
{\rm meas}(\mathcal P(\hat P))\lessdot
\left( \mu  ^{\frac{1}{m_n(m_n+1)(m_n+3)}}/\mathtt d_k\right)
\left|\log
\left( \mu  ^{\frac{1}{m_n(m_n+1)(m_n+3)}}/\mathtt d_k\right)
\right|\,.
$$
} 
\begin{eqnarray}
{\rm meas} \big( \{ P_n\in\big[a_+^{(2j-1)}(\hat P)-\mathtt r_0,
a_+^{(2j-1)}(\hat P)
\big)\ \ :\ \ 
|{\rm det}\left[\partial_{PP} \big(\hat h_k (\hat P)+ 
\mathtt E^{(2j-1)} (P)\big)
  \right] |\leq \mu  \} \big)
\lessdot \mu^{\mathtt a_n}\mathtt d_k^{-2}\,,
\nonumber
\\
\mathtt a_n:=\frac{1}{27 n^6}<
{\frac{1}{m_n(m_n+1)(m_n+3)}}\qquad
\label{missionbell2}
\end{eqnarray}
(recall \eqref{ciaone}) uniformly in $\hat P\in\hat D.$
Then by Fubini theorem we get
\begin{equation}\label{paradise}
{\rm meas} \big(\{ \ P\  | \  P_n\in\big[a_+^{(2j-1)}(\hat P)-\mathtt r_0,
a_+^{(2j-1)}(\hat P)\big),\ \hat P\in\hat D\,,\ \ 
|{\rm det}[\partial_{PP} (\hat h_k + 
\mathtt E^{(2j-1)} )
  ] |\leq \mu  \} \big)
\lessdot \mu^{\mathtt a_n}\mathtt d_k^{-2}\,.
\end{equation}

This conclude the proof close to the maxima of the odd case
$i=2j-1$.

\subsubsection*{The case $i$ even close to maxima}

The proof is analogous to the odd case,
leading to
\begin{equation}\label{paradise1}
{\rm meas} \big(\{ \ P\  | \  P_n\in\big[a_+^{(2j)}(\hat P)-\mathtt r_0,
a_+^{(2j)}(\hat P)
\big),\ \hat P\in\hat D\,,\ \ 
|{\rm det}[\partial_{PP} (\hat h_k + 
\mathtt E^{(2j)} )
  ] |\leq \mu  \} \big)
\lessdot \mu^{\mathtt a_n}\mathtt d_k^{-2}\,,
\end{equation}
\begin{equation}\label{paradise2}
{\rm meas} \big(\{ \ P\  | \  
P_n\in\big(a_-^{(2j)}(\hat P),
a_-^{(2j)}(\hat P)+\mathtt r_0\big]
,\ \hat P\in\hat D\,,\ \ 
|{\rm det}[\partial_{PP} (\hat h_k + 
\mathtt E^{(2j)} )
  ] |\leq \mu  \} \big)
\lessdot \mu^{\mathtt a_n}\mathtt d_k^{-2}\,.
\end{equation}

\bigskip

\subsubsection*{Far away from maxima}

We now study the point far away from maxima, where we note that the second derivatives of
the functions $\mathtt E^{(i)} $   are uniformly bounded
(see Lemma \ref{sorellenurzia2}).
We have
\begin{equation}\label{sugo}
{\rm det}\left[\partial_{PP} \big(\hat h_k (\hat P)+ 
\mathtt E^{(i)} (P)\big)
  \right]
  =
  {\rm det}\left[\partial_{\hat P \hat P} \hat h_k (\hat P)\right]\cdot 
\partial_{P_n P_n}\mathtt E^{(i)} (P)
\ +\ O(\tetta)
\stackrel{\eqref{lovebuzz}}=2^{n-1}\kappa^{-n}
\partial_{P_n P_n}\mathtt E^{(i)} (P)
\ +\ O(\checco)
\,,
\end{equation}
valid in the sets\footnote{Recall that $a_-^{(2j-1)}(\hat P)=0$ as defined in\eqref{jolebisk}.} 
\begin{eqnarray}
\Big\{\ P=(\hat P,P_n)\ | \ 
P_n\in\big(
0,
a_+^{(2j-1)}(\hat P)-\mathtt r_0/2
\big]\,, \ \hat P\in \hat D\ \Big\}\,,
\label{pioggia}
\\
\Big\{\ P=(\hat P,P_n)\ | \ 
P_n\in\big[
a_-^{(2j)}(\hat P)+\mathtt r_0/2,
a_+^{(2j)}(\hat P)-\mathtt r_0/2
\big]\,, \ \hat P\in \hat D\ \Big\}\,,
\label{pioggia2}
\end{eqnarray}
when $i$ is odd, respectively even.
  
\nl  
We have to distinguish the cases
when {\bf (A2)} or {\bf (A3)} hold.

\subsubsection*{The case when {\bf (A2)} holds}

We note that, by Proposition \ref{glicemiak}, we can extend the Hamiltonian 
$\mathtt E^{(2j-1)} (P)$
in a complex
$4\mathtt r_0$-neighborhood\footnote{Reducing, in case, $\mathtt r_0=\mathtt r_0(\pgot)>0$.}   of zero so that $\mathtt E^{(2j-1)} (P)$  has holomorphic extension on the complex
 domain 
\begin{equation}
\hat D_{r_0}\, \times \, 
\big[
-2\mathtt r_0,\
a_+^{(2j-1),0}-\mathtt r_0/4
\big]_{\mathtt r_0/8}\,,
\label{pioggia3}
\end{equation}
Then 
\eqref{sugo}
holds in the domain \eqref{pioggia3}.

\nl
Let us define the intervals
\begin{equation}\label{ciclette}
\mathtt I^{(2j-1)}:= \big[
-\mathtt r_0/4,
a_+^{(2j-1),0}-\mathtt r_0/4
\big]\,,\qquad
\mathtt I^{(2j)}:=
\big[
a_-^{(2j),0}+3\mathtt r_0/4,
a_+^{(2j)}(\hat P)-3\mathtt r_0/4
\big]
\end{equation}
Since we are far away from hyperbolic equilibria, by  Lemma \ref{sorellenurzia2}
 we get that 
\begin{equation}\label{soclose}
\sup_{\hat D_{r_0} \times 
\mathtt I^{(i)}_{\mathtt r_0/8}}
\|\partial_{PP}\mathtt E^{(i)} (P)\|\lessdot 1
\,,\qquad \forall\,
0\leq i\leq 2 N\,.
\end{equation}
Note that \eqref{sugo}
holds in 
$
\hat D \times\mathtt I^{(i)}
$
for every $i$
and also that, for $\checco$ small enough,
\begin{eqnarray}\label{ciclette2}
&&\Big(\hat D\times\mathtt I^{(2j-1)}\Big)\ \cup\
\Big\{ \ P\  | \  P_n\in\big[a_+^{(2j-1)}(\hat P)-\mathtt r_0,
a_+^{(2j-1)}(\hat P)\big),\ \hat P\in\hat D\Big\}
\ \supseteq \ 
\Pgot^{(2j-1)}(0)
\\
&&
\Big(\hat D\times\mathtt I^{(2j)}\Big)\ \cup\
\Big\{ \ P\  | \  P_n\in
\big(a_-^{(2j)}(\hat P),
a_-^{(2j)}(\hat P)+\mathtt r_0\big]\cup 
\big[a_+^{(2j)}(\hat P)-\mathtt r_0,
a_+^{(2j)}(\hat P)
\big)
,
\ \hat P\in\hat D\Big\}
\ = \ 
\Pgot^{(2j)}(0)
\nonumber
\end{eqnarray}
(defined in \eqref{playmobilk}).

\nl
Let us consider now the finitely many non constant analytic functions\footnote{$\mathtt E^{(i),0} $
being the inverse of $P_n^{(i),0}.$}
$$
P_n \to \partial_{P_n P_n}\mathtt E^{(i),0} (P_n)\,,\qquad
0\leq i\leq 2 N\,.
$$
By analyticity we have that there exist
$m=m(F^0)\geq 1$ and $\xi=\xi(F^0)>0,$
depending on the function 
$F^0$ introduced in \eqref{sgravone},
 such that
\begin{equation*}
\min_{P_n\in \mathtt I^{(i)}}\max_{ 1\leq d\leq m}|\partial_{P_n}^{d}\partial_{P_n P_n}
\mathtt E^{(i),0} (P_n)|/d!\geq 2\xi>0\,, \qquad \forall
0\leq i\leq 2 N \,.
\end{equation*}
By Lemma \ref{sorellenurzia2}
\begin{equation}\label{pettino}
\sup_{\hat D_{r_0} \times 
\mathtt I^{(i)}_{\mathtt r_0/8}}
|\partial_{P_n}\mathtt E^{(i)} (P)-
\partial_{P_n}\mathtt E^{(i),0} (P_n)|
\leq C\checco\,,
\qquad \forall
0\leq i\leq 2 N 
\,,
\end{equation}
where $C>1$ was defined in Proposition 
\ref{glicemiak}.
By Cauchy
estimates
we get, for $\checco$ small enough depending on $F^0$,
\begin{equation}\label{ventaglietti}
\inf_{P\in \hat D\times \mathtt I^{(i)}}\max_{ 1\leq d\leq m}|\partial_{P_n}^{d}\partial_{P_n P_n}
\mathtt E^{(i)} (P)|/d!\geq \xi>0\,, \qquad \forall
0\leq i\leq 2 N \,.
\end{equation}
Recalling \eqref{sugo}
 we have
$$
{\rm det}\left[\partial_{PP} \big(\hat h_k (\hat P)+ 
\mathtt E^{(i)} (P)\big)
  \right]
  =
2^{n-1}\kappa^{-n}
\partial_{P_n P_n}\mathtt E^{(i)} (P)
\ +\ O(\tetta)
=
2^{n-1}\kappa^{-n}\Big(
\partial_{P_n P_n}\mathtt E^{(i)} (P)
\ +\ O(\kappa^n\checco)
\Big)\,.
$$
By \eqref{gaetanga} we get that the term\footnote{More precisely
$|O(\kappa^n\checco)|\leq C
\kappa^n\checco.$} $O(\kappa^n\checco)$ above is negligible, together  with its derivatives of any order; then by 
\eqref{ventaglietti} we obtain
\begin{equation}\label{ventaglietti2}
\inf_{P\in \hat D\times \mathtt I^{(i)}}\max_{ 1\leq d\leq m}
\left|\partial_{P_n}^{d}
{\rm det}\left[\partial_{PP} \big(\hat h_k (\hat P)+ 
\mathtt E^{(i)} (P)\big)
  \right]
\right|/d!\geq 2^{n-2}\kappa^{-n} \xi>0\,, \qquad \forall
0\leq i\leq 2 N \,.
\end{equation}
We can apply Lemma \ref{girotondi}
uniformly in $\hat P\in\hat D$
with
$$
f(\cdot)\rightsquigarrow 
{\rm det}\left[\partial_{PP} \big(\hat h_k (\hat P)+ 
\mathtt E^{(i)} (\hat P, \cdot)\big)
  \right]\,,\ \ x\rightsquigarrow P_n\,,\ \ 
  [a,b]\rightsquigarrow \mathtt I^{(i)}\,,\ \ 
  \xi_m
  \rightsquigarrow 
2^{n-2}\kappa^{-n} \xi 
$$
and $M\lessdot 1$ (by \eqref{soclose}
and Cauchy estimates); then we get
$$
{\rm meas} \left( \left\{
P_n\in \mathtt I^{(i)}\ | \ 
\left|
{\rm det}\left[\partial_{PP} \big(\hat h_k (\hat P)+ 
\mathtt E^{(i)} (P)\big)
  \right]
\right|\leq \mu
\right\}\right)
\leq C
\kappa^n \mu  ^{\frac{1}{m(m+1)}}\,,
$$
uniformly in $\hat P\in\hat D.$
Then by Fubini theorem we have
\begin{equation}\label{city}
{\rm meas} \left( \left\{
P\in \hat D\times\mathtt I^{(i)}\ | \ 
\left|
{\rm det}\left[\partial_{PP} \big(\hat h_k (\hat P)+ 
\mathtt E^{(i)} (P)\big)
  \right]
\right|\leq \mu
\right\}\right)
\leq C
\kappa^n \mu  ^{\frac{1}{m(m+1)}}\,.
\end{equation}
Recalling the definition of $\mathtt a_n$ in \eqref{missionbell2}
we choose
$$
\mathtt c\geq \max\{ 6^3 n^6, m(m+1)\}
$$
Recalling 
\eqref{ciclette2}
by \eqref{paradise},\eqref{paradise1},
\eqref{paradise2} and \eqref{city}
we get \eqref{gibboso} when 
{\bf (A2)} holds.

\subsubsection*{The case when {\bf (A3)} holds}

We are in the cosine-like case (recall \eqref{drago}); by Lemma \ref{lapajata}
below
(see \eqref{thor}), \eqref{sugo} and
\eqref{gaetanga}
we have that
\begin{equation*}
\kappa^{-n}\lessdot 
\left|{\rm det}\left[\partial_{PP} \big(\hat h_k (\hat P)+ 
\mathtt E^{(i)} (P)\big)
  \right]\right|
\end{equation*}
for all the values of $P$ in \eqref{pioggia},  \eqref{pioggia2} respectively.
Then, by \eqref{gaetanga2}
\begin{equation}\label{sugo2}
\left|{\rm det}\left[\partial_{PP} \big(\hat h_k (\hat P)+ 
\mathtt E^{(i)} (P)\big)
  \right]\right|
  > \mu
\end{equation}
again for all the values of $P$ in \eqref{pioggia},  \eqref{pioggia2} respectively.
Then
\eqref{missionbell2}, \eqref{paradise}, 
\eqref{paradise2}, \eqref{sugo2}
prove Proposition~\ref{frey} in the case when 
{\bf (A3)} holds.
This concludes the proof of 
Proposition~\ref{frey}.
%%%%%%%%%%%%%

\section{Proof of the Main Theorem}
\label{lupin3}

\noi
In this final section we will show the existence of a high density (almost exponential) of Kolmogorov's tori of different topologies in 
all neighbourhoods of simple resonances. 

\subsection{Application of the Structure Theorem}

Fix $k\in\Z^n_{*,K}.$
In this subsection we will apply 
the Structure Theorem \ref{porretta}
to the ``effective part'' of the Hamiltonian $H_k$ in \eqref{disfare}.
Namely we apply Theorem 
\ref{porretta} with
\begin{eqnarray}
&&\cH(I',\f') \rightsquigarrow
 \frac{1}{\kappa} \|I'\|^2+\frac{1}{\d_k}
 G^k(\riscala I',k\cdot\f')\,,\quad
\Ggot (I',x)  \rightsquigarrow
 \frac{1}{\d_k}
 G^k(\riscala I',x)\,,\quad
F^0 \rightsquigarrow
  \frac{1}{\d_k} F^k\,,
 \nonumber
\\
&& 
 \cD \rightsquigarrow
D^{'k}\,,\quad
\hat D \rightsquigarrow
\frac{1}{\riscala}\hat Z_k\,,
\quad
R_0  \rightsquigarrow
R_{0,k}:=
\frac{K^{\nu+1}}{\sqrt{2\d_k}\|k\|}\,,
\quad
r_0 
 \rightsquigarrow
 4
 \,,
 \quad
 s_0 
 \rightsquigarrow
 s/4
 \,.
 \label{rana}
\end{eqnarray}
where $\d_k,$ $G^k$ and $F^k,$
$\riscala,$
$D^{'k},$
$\hat Z_k,$
were defined in 
\eqref{cenerentola},
\eqref{sugna},
\eqref{capocotta},
\eqref{fare},
\eqref{caspio}.
We  note that with such positions
we have that 
$r'_k$ defined in \eqref{fare}
satisfies
\begin{equation}\label{trooper}
r'_k=K^{\nu+1}/\sqrt {8\d_k}
\geq r'=c n |k|_\infty r_0\,,
\end{equation}
where $r'$ was defined in 
\eqref{santostefano}.
Indeed \eqref{trooper} follows 
by \eqref{sassa} taking $\e$
small enough (recall \eqref{lapparenza}).
We also note that
$s_*$ defined in \eqref{canarino}
satisfies
\begin{equation}\label{trooper2}
s_*\geq 3 s_0\,,
\end{equation}
taking $\e$ small enough
(recall again \eqref{lapparenza}).
Note that, by \eqref{battiato2}
\begin{equation}\label{battiato4}
\|g\|_{s_0}\leq \coth^n\left(\frac{s'-s_0}{2}\right)|g|_{s'}\,,\qquad
\forall\, g\in\cA^n_{s'}\,.
\end{equation}

\begin{remark}
 Here we take $\e\leq \e_0$
 small enough. Obviously
 $\e_0$ is small 
 {\sl uniformly}
 on $k.$
\end{remark}

\nl
We now verify that the hypotheses of Theorem \ref{porretta} hold.

\nl
Let us start with {\bf (A1)}.
We see that \eqref{Ccristina5} follow
by\footnote{Recall also  \eqref{fare} and \eqref{battiato4} 
with $n=1$.}
 \eqref{cristina}, resp.
\eqref{neftali}, if $|k|\leq K_s(\d),$
resp. $K_s(\d)<|k|\leq K,$
taking
\begin{equation}\label{pizzoli}
\checco_*:=\coth(s/24)
\frac{2^{10} n e^s}{\d s K^{\frac{4\nu-n-11}{2}}}\,.
\end{equation}
Moreover the inequality in \eqref{calimero} holds 
(recall the definition of $R_0$ in \eqref{rana})
taking $\e$
small enough (recall \eqref{lapparenza}).

\nl
We now check that {\bf (A3)}
holds
when $ K_s(\d)<|k|\leq K,$
taking the constant $\const$ in \eqref{miguel}
large enough.
Indeed we have that 
by \eqref{allblacks} and \eqref{battiato4}
(with $n=1$, also recall \eqref{cenerentola})
$\d_k^{-1}F^k$ is $\g$-cosine-like
taking 
\begin{equation}\label{larocca}
\g:=\coth(s/24)\left(
\frac{2|k|^{\frac{n+3}{2}} e^{-|k|s/4}}{\d}+
\frac{2^{10} n e^s}{\d s K^{\frac{4\nu-n-11}{2}}}\right)\,.
\end{equation}
Then \eqref{drago} holds taking
$c$ in \eqref{miguel}
large enough 
depending on $n$ in order to adjust the first addendum and, then, $\e$ small enough
in order to adjust the second addendum.

\nl
We finally check that {\bf (A2)}
holds
when $|k|\leq K_s(\d).$
We first note that in this case
$\d_k^{-1}F^k=F^k$ since 
$\d_k=1$ (recall \eqref{cenerentola}).
By \eqref{trilli} and \eqref{battiato2} we get
$$
\|F^k\|_{s_0} \leq 
\coth(3s/8)=:M\,.
$$
Moreover   every
$F^k$ is $(\b_k,M)$-Morse-non-degenerate for some $\b_k>0$
by (P2) and (P3) of Definition
\ref{albertone}.
Then we can take
$$
\b:=\min_{|k|\leq K_s(\d)}
 \b_k>0\,,
$$
since we are taking the minimum over a {\sl finite} set.
We conclude that
every
$F^k$ is $(\b,M)$-Morse-non-degenerate.

\begin{remark}\label{chestress}
A crucial fact is that we can choose the constant 
$\mathtt c$ appearing in Theorem \ref{porretta}, {\sl uniformly}
in $k\in Z^n_{*,K}$.
Indeed for all the cases
with $K_s(\d)<|k|\leq K,$
$\mathtt c$ only depends
on $n,s_0,r_0,$ namely, in view of \eqref{rana}, on 
$n$ and $s.$
In the {\sl finite} number of cases
$|k|\leq K_s(\d)$ the constant
$\mathtt c$ actually depends
on $k$ (since it depends on the particular for of $F^k$);
however we can the minimum
over all this finite number of constants obtaining a
constant $\mathtt c=\mathtt c(n,s)>0$
for which Theorem \ref{porretta} holds 
{\sl uniformly}
in $k\in\Z^n_{*,K}.$
\end{remark}

\nl
Finally we have that 
\eqref{gaetanga} and, a fortiori
\eqref{tolomeo}, are 
simultaneously, namely for every $k\in\Z^n_{*,K},$ satisfied  taking $\e$
small enough
(recall \eqref{lapparenza}).
Then, as a corollary of Theorem 
\ref{porretta}, we get the following
result (where we denote by $N_k$
the number of maxima/minima of
$F^k$).

\begin{proposition}\label{cermone}
Let $\e$ be small enough.
Let $\theta>0$ and\footnote{The constant $\mathtt c$ is defined in Remark \ref{chestress}.}
\begin{equation}\label{gaetanga2!}
0<\mu\leq 1/\mathtt c K^{2n}\,,
\end{equation}
 where $\mathtt c=\mathtt c(n,s)>0$
 was defined in Remark \ref{chestress}.
 For every $k\in\Z^n_{*,K}$ and
$0\leq i\leq 2N_k$ there exist
\\
i) disjoint subsets
$\Cgot^i_k(\theta,\mu)\subseteq D^{'k}\times\T^n$ 
decreasing w.r.t. $\theta$ and $\mu$, with\footnote{$D^{'k}_\sharp$ was defined in \eqref{foggydew}.}
\begin{equation}\label{fujiyama!}
\meas\Big(
 \big(D^{'k}_\sharp\times \T^n\big) \setminus 
 \bigcup_{0\leq i\leq 2N_k}
\Cgot^i_k(\theta,\mu) \Big) \leq \mathtt c \big(\theta |\ln\theta|+K^{4\nu}
\mu^{1/\mathtt c}\big)\,;
\end{equation}
ii)  $\mathtt B^i_k(\theta,\mu)\subseteq\R^n$, decreasing w.r.t. $\theta$
and $\mu$, with\footnote{With
$R_{0,k}$ defined in \eqref{rana} and $c$
in \eqref{verruca} respectively.}
\begin{equation}\label{pratola!}
{\rm diam} \big( \mathtt B^i_k(\theta,\mu)\big)\leq 2c \big(R_{0,k} +\riscala^{-1}
{\rm diam}(\hat Z_k)\big)\,;
\end{equation}
\\
iii)  holomorphic 
symplectomorphisms\footnote{
$r'_k$ defined in \eqref{fare}.}  
\begin{equation}\label{blueeyes3!}
\Psi^i_k \ :\ \big(\mathtt B^i_k (\theta,\mu)\big)_{\rho'}\times\T^n_{\s'}
\ \to \ D^{'k}_{r'_k}\times\T^n_{s/4}
\,,\qquad
\mbox{with}\quad
\rho':=\frac\theta{\mathtt c K^{n-1}} \,,
\quad  \ \s':=\frac1{{\mathtt c}K^{n-1}
|\log \theta|}\,,
\end{equation}
with
\begin{equation}\label{tokyo!}
\Psi^i_k  
\Big(\mathtt B^i_k (\theta,\mu)\times\T^n \Big)=
\Cgot^i_k (\theta,\mu)\,,
\end{equation}
such that\footnote{$H_k$
defined in \eqref{disfare}.} 
 \begin{equation}\label{matteo!}
H_k
\circ \Psi^i(p,q)=:
h^{(i)}_k(p)+f^{(i)}_k(p,q)\,.
\end{equation}
with\footnote{By \eqref{avezzano2}
and \eqref{battiato4}, taking $\e$ small enough,
depending on $s$ and $\d$.}
\begin{equation}\label{avezzano}
\|f^{(i)}_k\|_{\mathtt B^i_k (\theta,\mu),\rho',\s'}
\leq
 \e^{\frac{s}5 |\log \e|^3}
\end{equation}

\noindent
Furthermore 
\begin{equation}\label{gadamer!}
\left| {\rm det}\left(\partial_{pp} h ^{(i)}_k (p)
    \right)  \right|
\geq \mu\,,\qquad
\forall\, \ 0\leq i\leq 2N \,,\ \ |k|\leq K\, , \qquad \forall \ p\in  {\mathtt B}^i_k (\theta,\mu)\ .
\end{equation}
Finally\footnote{Recall \eqref{vallinsu2}. 
}
\begin{equation}\label{vallinsu!}
\|\partial_{pp} h^{(i)} \|_{ \mathtt B^i_k (\theta,\mu), \rho'}\leq \mathtt c/\theta\,,\qquad
\text{for }\ \ \ 0\leq i\leq 2N\,.
\end{equation}

\end{proposition}

\subsection{Application of the KAM theorem}

Let us start by stating a quantitative version, suitable for our purposes, of the classical KAM Theorem; for references, discussions and extensions we refer to \cite{AKN} and references therein; see also \cite{KAMstory} for a nice divulgative account of KAM theory. \\

 \begin{theorem}\label{KAM}
Fix $n\ge 2$ and $\t>n-1$. Let 
$\ttD$ be any non--empty, bounded subset of $\real^n$. 
Let
\beqno%{insonnia}
\ttH(p,q):=h(p)+f(p,q)
\eeqno
  real--analytic on $\ttD_{r_0}\times\T^n_s,$
 for some $r_0>0$ and $0<s\le 1$, and having finite norms:
\beq{bellini}
\ttM:= \|\partial_{pp} h\|_{\ttD, r_0} \,,\qquad\qquad
\|f\|_{\ttD,r_0,s} \,.
\eeq
Assume that the frequency map $p\in \ttD \to\o=\partial_{p} h$ is a local diffeomorphism, namely, assume:
\begin{equation}\label{dupa}
d:=  \inf_\ttD|\det \partial_{pp} h| >0\,.
\end{equation}
 Define
\beq{bove}
\ttm:=\frac{d}{\ttM^n}\le 1\ .
\eeq
Then there exists a positive constant $c<1$, depending only on $n$ and $\t$, such that, if   
\beq{enza} 
\epsilon:= \frac{\|f\|_{\ttD,r_0,s} }{\ttM r_0^2}\le c \,\ttm^8\ s^{4\tau+4}\ , 
\eeq
then the following holds. 
Define
\begin{equation}\label{nicaragua}
\a:=\frac{c}{  \ttm\, s^{3\tau+3}} \, (\ttM r_0)\,\sqrt\epsilon\ ,\qquad
\hat r:= \, \ttm^2 r_0\ ,\qquad r_\epsilon:=\frac{1 }{c \,\ttm}\, \sqrt\epsilon\, r_0\ .
\end{equation}
Then, there exists a positive measure 
%``non--torus set'' 
set ${\cal T}_\a\subseteq (\ttD_{\hat r}\cap\R^n)\times \torus^n$ formed by ``primary''  Kolmogorov's tori; more precisely, 
for any point $(p,q)\in{\cal T}_\a$, $\phi^t_\ttH(p,q)$ covers densely an $\ttH$--invariant, analytic, Lagrangian torus,  with $\ttH$--flow analytically conjugated to a linear flow with $(\a,\t)$--Diophantine frequencies 
$\o=h_p(p_0)$, for a suitable $p_0\in\ttD$; each of such tori is a graph over $\torus^n$ $r_\epsilon$--close
to the unperturbed  trivial graph  $\{(p,\theta)=(p_0,\theta)|\ \theta\in \torus^n\}$.\\
Finally, 
the Lebesgue outer measure of $(\ttD\times\torus^n)\bks {\cal T}_\a$ is bounded by:
\beq{gusuppo}
\meas \big((\ttD\times\torus^n)\bks {\cal T}_\a\big) \le  C\, \sqrt\epsilon
\eeq
with
\beq{ciccione}
C:=
\big(\max\big\{ \ttm^2 r_0\, ,\, \diam\, \ttD\big\}\big)^n \cdot \frac{1}{c\, \ttm^{n+5}\ s^{3\tau+3}}\,.
\eeq

\end{theorem}

\begin{remark}
 The statement of the above quantitative KAM theorem
 is as Theorem 1 in \cite{BCKAM} with the following minor simplification. In Theorem 1 of \cite{BCKAM}
 appear the quantity 
 $\l:=\ttL \ttM$, where $\ttM$ is defined in 
 \eqref{bellini} and $\ttL$ denotes a suitable uniform Lipschitz constant of the local complex inverse of the ``frequency map'' $p\mapsto \o=\partial_p h(p)$ (compare formula (9) of \cite{BCKAM}); since one can show that $1\leq \l\le 2\cdot n! \, \ttm^{-1}$ (see formula 
 (14) of \cite{BCKAM}), we  substitute  everywhere
 $\l$ with $1$ or $2\cdot n! \, \ttm^{-1}$ in Theorem \ref{KAM},
obtaining a slightly weaker formulation of  Theorem 1 of \cite{BCKAM}.
  
\end{remark}

%%%%%%%%%%%
\Giu

\subsubsection*{KAM tori in $\O^0$}

We now apply Theorem~\ref{KAM} to the Hamiltonian
$H_{\{0\}}$
in \eqref{prurito}.
It is immediato to see, thanks to \eqref{552},
that KAM tori cover all $\O^0$ (defined in \eqref{neva})
up to a set of measure
$\e^{\frac{s}{5}|\log \e|}.$
\Giu

\subsubsection*{KAM tori in $\O^1$}
We, now, want to apply Theorem~\ref{KAM} to the Hamiltonians $h^{(i)}_k(p)+f^{(i)}_k(p,q)$ defined in \eqref{matteo!}, for all $k\in \integer^n_{*,K},$
$0\leq i\leq 2N_k.$

\nl
The objects appearing in Theorem~\ref{KAM} have to be replaced by the following:
\beq{cipolla}
 h \vain h^{(i)}_k\,,\qquad f \vain f^{(i)}_k\,,
 \qquad   
 \mathtt D\vain 
 \mathtt B^i_k(\theta,\mu)\,,
 \qquad
 r_0 \vain \rho'=\frac\theta{\mathtt c K^{n-1}} \,,
\qquad   
 s \vain \s':=\frac1{{\mathtt c}K^{n-1}}\,,  
\eeq
(recall  \equ{blueeyes3!}). 

\nl
By \equ{avezzano}, it follows immediately that
\beq{bureca1}
\|f\|_{\ttD,r_0,s}\, \le  \e^{\frac{s}5 |\log \e|^3}\ .
\eeq
By \eqref{vallinsu!}
 we get
\beq{bureca2}
\mathtt M\, \le \const /\theta\,,
\eeq 
where, here and in the following,
$$
c=c(n,s)\geq 1\,,
$$ 
are suitably large (different) constants depending only on $n$ and $s.$
By \eqref{gadamer!} we have that
$d$ defined in \eqref{dupa}
satisfies 
$$
d\geq \mu\,.
$$
So we get that $\ttm$ in \eqref{bove}
satisfies
\begin{equation}\label{bovissimo}
\ttm=\frac{d}{\ttM^n}\geq \mu\theta^n/\const\,.
\end{equation}
Now, we choose the parameters 
 $\m$ and $\theta$ as follows\footnote{Where
 $\mathtt c$ is the constant defined in Proposition
 \ref{cermone}.}
\beq{sesemorotti}
\m=\theta:= \e^{|\log \e|^2}
\,.
\eeq
With such choices the condition of the KAM Theorem \ref {KAM} are met: 
in particular \eqref{bove} follows by \eqref{bovissimo}
and also \eqref{enza},
which is implied by the stronger condition
$$
\|f\|_{\ttD,r_0,s} \le c \,
\mu^8\theta^{8n+1}K^{2-2n} s^{4\tau+4}\ ,
$$
which holds by \eqref{bureca1} (recall also
\eqref{lapparenza}), taking $\e$ small enough.

\nl
Noting that, by \eqref{pratola!}, \eqref{rana}
and 	\eqref{capocotta}, 
 \begin{equation}\label{pratola5}
{\rm diam} (\mathtt D)\leq \frac{cK^{\nu+1}}{\sqrt{\d_k}\|k\|}
+\riscala^{-1}
{\rm diam}(\hat Z_k)
\leq
c\frac{K^{\nu+1}\sqrt \e}{\riscala}
+\riscala^{-1}
{\rm diam}(\hat D)
\stackrel{\eqref{arrosticini}}
\leq
c\frac{K^{\nu+1}\sqrt \e}{\riscala}
+\riscala^{-1} K^n
\,,
\end{equation}
the maximum in \eqref{ciccione} 
can be estimated by
$c\frac{K^{\nu+1}}{\riscala}.$ 
By \eqref{gusuppo} and \eqref{ciccione}
we get that the measure of the non torus set 
in every $ \mathtt B^i_k(\theta,\mu)$
is bounded by
$$
c \frac{K^c}{\riscala^n} \sqrt{\|f\|_{\ttD,r_0,s}}
\frac{\ttM^{n(n+5)-1/2}}{d^{n+5}s^{3\tau+3}}
\leq  \frac{c}{\riscala^n \e^c}\e^{\frac{s}{10} |\log \e|^3}
\leq  \frac{1}{\riscala^n }\e^{\frac{s}{11} |\log \e|^3}\,,
$$
for $\e$ small enough.
Then by \eqref{fujiyama!} (and \eqref{tokyo!})
we get that the measure of the non-torus 
set in 
$\big(D^{'k}_\sharp\times \T^n\big) $
is estimated by\footnote{Note that
$\riscala\leq 1$, see \eqref{capocotta}.}
$$
\frac{1}{\riscala^n }\e^{ 2|\log \e|}\,,
$$
for $\e$ small enough.
Therefore
the measure of the non-torus 
set in $\bigcup_{k\in\Z^n_{*,K}}D^{'k}_\sharp\times\T^n$  
is estimated by
$$
\frac{1}{\riscala^n }\e^{ |\log \e|}\,,
$$
for $\e$ small enough (recall \eqref{lapparenza}).
By \eqref{foggydew3}, \eqref{prospettiva} and
\eqref{capocotta}
we get that {\sl the measure of the non-torus 
set in $\O^1\times\T^n$ is bounded by
$\e^{ |\log \e|}.$}

%%%%%%%%%%%%%%%

\giu
{\bf Acknowledgment.}
We are indebted to V. Kaloshin, G. Loddi, A. Neishtadt and A. Sorrentino.

\appendix

\section{Properties of the class of non--degenerate potentials}\label{pogba}
{\bf Proof} of Proposition~\ref{grizman}.

\giu
$\bullet$ {\sl $\cP_s\, \cap \, \bB_s^n\in \cB$ and $\m_s(\cP_s\cap \bB_s^n)=1$}

\giu
We shall prove that, for every $\d>0,$ 
the measure of the sets of
potentials $f$ that do not satisfy, respectively, (P1), (P2), (P3), (P4) is,
respectively,  $O(\d^2),0,0,0$, the result will follow letting $\d\to 0$.

\nl
First, by the identification \eqref{miserere}, the
measure of the set of
potentials $f$ that do not satisfy (P1)
with a given $\delta$
 is  bounded by $\d^2\, \sum_{k\in\integer^n} |k|^{-n-3}$.

\nl
Next, recall that properties (P2), (P3) and (P4) concern only a {\sl finite} number of $k$,
i.e.,  $k\in\integer^n_*,\  |k|\leq  K_s(\d)$.

\nl 
To show that the set of potentials that do not satisfy (P2) has $\m_s$-measure zero it is enough to check that, for every $k\in\integer^n_*,\  |k|\leq  K_s(\d)$, the set 
${\cal E}^{(k)}$ 
of $f$'s
for which\footnote{Recall the definition of $F^k$ in  \equ{dec}.} $F^k$ has a degenerate critical point has zero $\m_s$-measure. \\
Fix $k\in\integer^n_*,\  |k|\leq  K_s(\d)$ and 
denote points in  $\cE^{(k)}$ by $(\zeta,\f)$, where $\zeta=f_k$ and $\f=\{f_h\}_{h\neq k}$.
Write
\begin{equation}\label{efesta}
F^k(\x)=
\zeta e^{\ii \x} + \bar\zeta e^{-\ii \x} + G(\x)\,,\quad
{\rm where}\ \ 
\zeta:= f_k\ \ {\rm and}\ \ 
 G(\x):=\sum_{|j|\geq 2} f_{jk}e^{\ii j\x}\,.
\end{equation}
Now, one checks immediately that $\partial_\xi F^k(\x_0)=0=\partial^2_\xi F^k(\x_0)$ is equivalent to 
$\zeta=\zeta(\x_0,\f)=\frac12 e^{-\ii \x_0} \big( \ii G'(\x_0)+G''(\x_0) \big)$, which, as $\x_0$ varies in $\torus$, describes a smooth closed ``critical'' curve in $\complex$,
as a side remark, notice that  $\zeta$ depends on $\f$ only through the Fourier coefficients $f_{jk}$ with $|j|\ge 2$.  Thus the section $\cE^{(k)}_\f=\{\zeta\in D: (\zeta,\f)\in\cE^{(k)}\} $ is (a piece of) a smooth curve in $D=\{z\in\complex: |z|\le 1\}$, hence meas$(\cE^{(k)}_\f)=0$ for every $\f$ and by Fubini's theorem $\m_s (\cE^{(k)})=0$, as claimed.

\nl
An analogous  result\footnote{In this case  the critical curve is given by
$\{
\zeta=(-b(\x)\pm \sqrt{b^2(\x)-c(\x)}+\ii G'(\x)) e^{-\ii \x}/2, \ \x\in\real\,, b^2(\x)\geq c(\x)\},$
where $b(\x):=(G''''(\x)-G''(\x))/2$ and $c(\x):=-G''(\x) G''''(\x) +5(G'(\x)+G'''(\x))^2/3$.} 
holds true  for 
(P3). 

\nl
Regarding (P4) we have that 
the three {\sl real} equations
$$
\partial_\xi F^k(\xi_1)=\partial_\xi F^k(\xi_2)=0\,,\qquad
F^k(\xi_1)-F^k(\xi_2)=0\,, \qquad {\rm for}\ \ \ \xi_1,\xi_2\in\T\,,
$$
can be rewritten as (recall \eqref{efesta}) the {\sl complex}
equation
\begin{equation}\label{celebration}
\zeta=\zeta(\x_1,\xi_2,\f)= \frac{\ii}{2(e^{\ii \xi_1}-e^{\ii \xi_2})}
\Big(
G'(\xi_1)-G'(\xi_2) +\ii G(\xi_1) -\ii G(\xi_2)
\Big)
\end{equation}
and the real one
\begin{eqnarray}
\frac12(e^{\ii (\xi_1-\xi_2)}- e^{-\ii (\xi_1-\xi_2)})
\Big(
G'(\xi_1)-G'(\xi_2) +\ii G(\xi_1) -\ii G(\xi_2)
\Big)
-
(e^{\ii \xi_1}-e^{\ii \xi_2})
\Big(
e^{-\ii \xi_2} G'(\xi_1) - e^{-\ii \xi_1} G'(\xi_2)
\Big)
\nonumber
\\
=:g(\xi_1,\xi_2,\f)
=
\big(1-\cos(\xi_1-\xi_2)\big) \Big( G'(\xi_1)+G'(\xi_2) \Big)
- \sin (\xi_1-\xi_2)\Big( G(\xi_1)-G(\xi_2) \Big) 
=0\,.
\label{celebration2}
\end{eqnarray}
We claim that, for every fixed $\f$, the analytic function 
$(\xi_1,\xi_2)\mapsto g(\xi_1,\xi_2,\f)$ 
is not identically zero and, therefore, 
the set
$R_\f$ of its zeros has zero measure. 
Assume by contradiction that $g$ is identically zero.
Then $g(\xi_2+\e,\xi_2,\f)\equiv 0$ for every $\xi_2$ and 
$\e,$ in particular, evaluating the order fourth term
of the Taylor expansion in $\e$ around 
$\e=0,$ 
we get
$
\frac{1}{12} \Big(
G'''(\xi_2)+G'(\xi_2)
\Big)=0\,,\  \forall\, \xi_2\,.
$
The general (real) solution of the above equation 
is  $G(\xi_2)= c e^{\ii \xi_2} + \bar c e^{-\ii \xi_2}+c_0,$
with $c\in\C,$ $c_0\in\R, $
which contradicts  the expression of $G$
in \eqref{efesta}.
Therefore, for every fixed $\f,$
the image of the zero measure set $R_\f$
through the Lipschitz function  $(\xi_1,\xi_2)\mapsto
\zeta(\x_1,\xi_2,\f)$ (defined in \eqref{celebration}) 
has zero measure in $D.$  
Then we conclude as in the case (P2) above.

\Giu
$\bullet$  {\sl  $\cP_s$ contains an open subset $\cP_s'$ which is dense in the unit ball of  $ \cA_s^n$.}

\nl
Let us define $\cP_s'$  as $\cP_s$ but with the difference that  (P1) is replaced by the
{\sl stronger condition}\footnote{Note that $\m_s(\cP'_s)=0$.}

\nl
(P$1'$)
$\exists\, \d>0$ s.t. 
$\displaystyle{  |f_k|\geq\d\ e^{-|k|s} \,,\ \ 
\forall\,k\in\integer^n_*,\  |k|> K_s(\d)}$.

\nl
Let us first prove that $\cP_s'$ is open. Let $f\in\cP_s'$.
We have to show that there exists $\rho>0$ such that if $|g|_s< \rho,$ then $f+g\in\cP_s'$.
Fix $\d>0$ such that (P$1'$) holds and choose $\rho<\d$ small enough such that
$
[K_s(\d)] > K_s(\d')-1\,,$ where $ \d':=\d-\rho
$
and $[\cdot]$ denotes integer part.
Then, it is immediate to verify that $|k|> K_s(\d) \iff |k|> K_s(\d')$.
Moreover
$$
|f_k+g_k|e^{|k|s}\geq |f_k| e^{|k|s} -|g|_s \geq \d-\rho=\d'\,, \qquad
\forall\,k\in\integer^n_*,\  |k|> K_s(\d')\,,
$$
namely $f+g$ satisfies  (P$1'$) (with $\d'$ instead of $\d$).
Since (P2), (P3) and (P4) are ``open'' conditions and regard only a finite  number of $k$
it is simple to see that they are satisfied also by $f+g$ for $\rho$ small enough.
Then $f+g\in\cP_s'$  for $\rho$ small enough.

\nl
Let us now show that $\cP_s'$  is dense in the unit ball of  $ \cA_s^n$.
Take $f$ in the unit ball of  $ \cA_s^n$ and $0<\theta <1$. We have to find $\tilde f\in\cP_s'$
with $|\tilde f-f|_s\leq \theta $.
Let $\d:=\theta /4$ and denote by  $f_k$ and $\tilde f_k$ (to be defined) be the Fourier coefficients of, respectively,  
$f$ and $\tilde f$. We, then, let $\tilde f_k=f_k$ unless one of the following two cases occurs: 

\begin{itemize}
\item $k\in\integer^n_*$, $|k|> K_s(\d)$ and $|f_k|e^{|k|s}< \d$, in which case, $\tilde f_k= \delta e^{-|k|s}$,

\item $k\in\integer^n_*$, $|k|\leq K_s(\d)$ and  $F^k$ (defined as in \eqref{dec}) does not satisfy either (P2), (P3) or (P4), in which
case, $\tilde f_k$ is chosen at a distance less than $\theta e^{-|k|s}$ from $f_k$  but outside the critical curves defined above.

\end{itemize} 
At this point, it is easy to check that  $\tilde f\in \cP'_s$ and is $\theta $--close to $f$. 

\Giu
$\bullet$  {\sl  $\cP_s$ is prevalent}.

\giu
Consider the following compact subset of
$\ell_\io^n$: 
let $\mathcal K:=\{ z=\{z_k\}_{k\in \integer^n_\sharp  } : z_k\in D_{1/|k|} \},$
where $D_{1/|k|}:=\{w\in\complex:\ |w|\le 1/|k|\},$ 
and let  $\nu$ be 
the unique probability measure supported on $\mathcal K$ such that,
given  Lebesgue measurable sets 
$A_k\subseteq D_{1/|k|}$,
with $A_k\neq D_{1/|k|}$
only for finitely many $k$, one has
$$
\nu \Big(\prod_{k\in\Z^n_\sharp} A_k\Big):=
\prod_{\{k\in \integer^n_\sharp:\, A_k\neq D_{1/|k|}\}} 
\frac{|k|^2}{\pi}{\rm meas}(A_k)\,.
$$ 
The isometry $j_s$ in \eqref{miserere}
naturally induces a  probability measure
$\nu_s$ on  $\cA^n_s$ with support in the compact set $\mathcal K_s
:=j_s^{-1}\mathcal K$.
Now, for $\d>0,$ let $\cP_{s,\d}$ be the set of $f$'s
 in the unit ball of $ \cA_s^n$
satisfying (P1)--(P4), so that $\cP_s=\cup_{\d>0}
\cP_{s,\d}$.
Reasoning as in the proof of 
$\mu_s(\cP_s)=1,$ 
one can show that
$\nu_s(\cP_{s,\d})
\geq 1-{\rm const}\, \d^2$.
It is also easy to check that, for every $g\in  \cA_s^n$, the translated 
set $\cP_{s,\d}+g$ satisfies
$\nu_s(\cP_{s,\d}+g)\geq \nu_s(\cP_{s,\d})$.
Thus,  one gets 
$\nu_s(\cP_s+g)= \nu_s(\cP_s)
= 1$, $\forall\, g\in  \cA_s^n$, which means  that
$\cP_s$ is prevalent.
(recall footnote~\ref{nurzia}) \qed

%\appendix

\section{Proof of the Normal Form Lemma~\ref{pesce}}\label{provapesce}
Given a function $\phi$ we denote by $X_\phi^t$
the hamiltonian flow at time $t$ generated by $\phi$
and by ``ad''  the linear operator  $u\mapsto {\rm ad}_\phi u:=\{u,\phi\}$ and ${\rm ad}^\ell$ its iterates:  
$$
{\rm ad}^0_\phi u:=u\,,  \qquad
{\rm ad}^\ell _\phi u:=\{ {\rm ad}^{\ell -1}_\phi u,\phi\}\,, \qquad \ell \geq 1\,,
$$
as standard, $\{\cdot, \cdot\}$ denotes Poisson bracket\footnote{Explicitly, 
$\dst \{u,v\}= \sum_{i=1}^n (u_{x_i}v_{y_i}- u_{y_i}v_{x_i})$.}.

\nl
Recall the identity (``Lie series expansion'') 
\begin{equation}\label{ellade}
u\circ X_\phi^1 =\sum_{\ell \geq 0}
 \frac{1}{\ell !} {\rm ad}^\ell_\phi u=
\sum_{\ell=0}^\io \frac{\partial_t^\ell   (u\circ X_\phi^t)}{\ell!} \Big|_{t=0} 
 \, ,
\end{equation}
valid for analytic functions  and small $\phi$.

\nl
By standard Cauchy estimates,
we get (compare, e.g., Lemma~B4 of \cite{poschel})
\begin{lemma}\label{messene}
For $0<r-\rho<r_0,$ $0<s-\s<s_0,$ $\rho,\s>0$
\begin{equation}\label{tirinto}
|\{f,g\}|_{r-\rho,s-\s}
\leq
\frac ne
\left(
\frac{1}{(r_0-r+\rho)\s} +
\frac{1}{(s_0-s+\s)\rho }
\right)
|f|_{r_0,s_0} |g|_{r,s}\,.
\end{equation}
\end{lemma}
Summing the Lie series in \eqref{ellade} and using  Lemma B5 of \cite{poschel}, we get, also, 
\begin{lemma}\label{olimpiabis}
Let $0<\rho<r\leq r_0-\rho$ and 
$0<\s<s\leq s_0-\s.$
Assume that
\begin{equation}\label{tebe}
\hat\tetta:=\frac{4n|\phi|_{r_0,s_0}}{\rho\s}\leq 1\,.
\end{equation}
Then
\begin{equation}\label{delfi}
\big| u\circ X_\phi^1 - u  \big|_{r-\rho,s-\s}
\leq 
\sum_{\ell \geq 1} \frac{1}{\ell !}
 \big| {\rm ad}^\ell_\phi u\big|_{r-\rho,s-\s}
\leq
\hat\tetta |u|_{r,s}\,.
\end{equation}
\end{lemma}
Given $K\ge 2$ and a lattice $\L$, recall the definition of $f^\flat$ in \equ{senzanome} and define 
\beq{fK}
f^K:=f-f^\flat=T_{K}\proiezione_\L^\perp  f\ ,
\eeq
so that we have the decomposition (valid for any $f$):
\beq{decomposizione}
f=f^\flat+f^K\,,\qquad
f^\flat:= P_\L f+\Tp\proiezione_\L^\perp  f \,,\qquad
f^K:= T_{K}\proiezione_\L^\perp  f\,.
\eeq

\begin{lemma}\label{megara}
Consider a real--analytic Hamiltonian 
\begin{equation}\label{olinto}
H=H(y,x)=h(y)+f(y,x)\qquad\mbox{analytic \ on \ } D_r\times \T^n_s\,.
\end{equation}
Suppose that $D_r$ is ($\a$,$K$) non--resonant modulo $\L$
for $h$ (with $K\ge 2$).
Assume that
\begin{equation}\label{gricia*}
\check\tetta := \frac{2^5 n K^3}{\a r s}\,  |f^K|_{r,s}\leq 1\,.
\end{equation}
Then there exists a real--analytic 
 symplectic change of coordinates
$$
\Psi: D_{r_+}\times \T^n_{s_+} 
\to D_r \times \T^n_s \,,\qquad
r_+:=r(1-1/2K)\,,\ \ \ s_+:=s(1-1/K^2)\,,
$$
such that 
\begin{equation}\label{olintobis}
H\circ\Psi=h(y)+f_+(y,x)\,, \qquad f_+:=f^\flat+f_*
\end{equation}
with
\begin{equation}\label{salamina}
|f_*|_{r_+,s_+} 
\leq\ 
2\check\tetta |f|_{r,s}
\,.
\end{equation} 
\end{lemma}
Notice that, by \equ{decomposizione} and \equ{salamina} (and the fact that $|f-f^K|_{r,s}\le |f|_{r,s}$), one has
\beq{nuvole}
f_+^K=f_*^K\ , \quad |f_+|_{r_+,s_+}=|f_*+f-f^K|_{r_+,s_+}\le 
|f_*|_{r_+,s_+}+|f|_{r,s}\le
(1+2\check\tetta)|f|_{r,s}\ .
\eeq
Notice also that 
\beq{bemolle}
f_+^\flat-f^\flat\eqby{olintobis} f_*^\flat\quad \Longrightarrow \quad |f_+^\flat-f^\flat|_{r_+,s_+}\le
|f_*|_{r_+.s_+}\stackrel{\equ{salamina}}\le 
 2 \check\tetta |f|_{r,s}\ .
\eeq
\proof (of Lemma~\ref{megara})
Let us define
$$
\phi=\phi(y,x):=\sum_{|m|\leq K, m\notin\L}
\frac{f_m(y)}{\ii h'(y)\cdot m} e^{\ii m\cdot x}
\,,\qquad
\Psi:=X^1_\phi
\, ,
$$
and note that $\phi$ solves the homological equation
\begin{equation}\label{tessalonica}
\{ h,\phi\}+f^K=0\,.
\end{equation}
Since $D_r$ is ($\a$,$K$) non--resonant modulo $\L$
\begin{equation}\label{maratona}
|\phi|_{r,s}\leq |f^K|_{r,s}/\a\,.
\end{equation}
Then, one has 
$$
H\circ\Psi=h+f^\flat+f_*
$$
with
$$
f_* 
=
 (h\circ\Psi-h-\{h,\phi\})+
 (f\circ\Psi-f) \,.
 $$
In order to estimate $f_*$
we now use Lemma \ref{olimpiabis} with parameters
$$r_0\rightsquigarrow r\,,\ \  
s_0 \rightsquigarrow s\,,\ \ 
r \rightsquigarrow r(1-1/4K)\,,\ \ 
s \rightsquigarrow s(1-1/2K^2)\,,\ \ 
\rho \rightsquigarrow r/4K\,,\ \  
\s \rightsquigarrow s/2K^2\,.
$$
With these choices  it is  
$\hat\tetta= \check\tetta$, and, by \eqref{gricia*}  $\check\tetta\leq 1$.
Thus,  \eqref{tebe} holds and Lemma \ref{olimpiabis} applies.
By \eqref{delfi} we get
\eqref{salamina} noting that
$$
h\circ\psi-h-\{h,\phi\}
=\sum_{\ell \geq 2}
 \frac{1}{\ell !} {\rm ad}^\ell_\phi h
 =
 \sum_{\ell \geq 1}
 \frac{1}{(\ell +1)!} {\rm ad}^\ell_\phi \{h,\phi\}
 \stackrel{\eqref{tessalonica}}=
- \sum_{\ell \geq 1}
 \frac{1}{(\ell +1)!} {\rm ad}^\ell_\phi 
f^K\,,
$$
which implies  (again by \eqref{delfi})
that 
$$|h\circ\Psi-h-\{h,\phi\}|_{r_+,s_+}
\leq \check\tetta |f^K|_{r,s}\le \check\tetta |f|_{r,s}\ .$$
Finally, applying again Lemma~\ref{olimpiabis} with $u=f$, by \equ{delfi}, we get $|f\circ\Psi-f|_{r_+,s_+}\le \check\tetta |f|_{r,s}$, concluding the proof of Lemma~\ref{megara}. 
%\qed
\eproof

\Giu
\textbf{Proof of the Normal Form Lemma \ref{pesce}}
Denote by
\beq{tetto}
\bar K:= \lceil K\rceil:=\min \{n\in\integer:\ n\ge K\}\ ,
\eeq
the ceiling function of $K$.
The idea is to construct $\Psi$ by applying $\bar K$ times Lemma~\ref{megara}. 

\nl
To do this, fix $1\le j< K$ and
make the following {\sl inductive assumption}: \\
{\sl Let 
\beqa{cornetto}
&& f_0:=f \ ,  \quad H_0:=h+f_0=H\,, \quad
\rho:=\frac{r}{4\bar K}\,, \qquad
\s:=\frac{s}{2K\bar K}
\,,
\nonumber\\
&& 
\nonumber\\
&&r_i:=r- 2i \rho\,,\qquad
s_i:=s- 2i \s\,,\qquad
|\cdot|_i:=|\cdot|_{r_i,s_i}\,,
\eeqa
and  assume that there exist, for $1\le i\le j$, real--analytic symplectic transformations
$$
\Psi_{i-1} \ :\ D_{r_i}\times \T^n_{s_i}\to 
D_{r_{i-1}}\times \T^n_{s_{i-1}}\,,
$$
such that
\begin{equation}\label{olintoj}
H_i:=H_{i-1}\circ \Psi_{i-1}
=:h+f_i
\end{equation}
satisfies, for $1\leq i\leq j$,  the estimates
\begin{equation}\label{pontina}
\tetta_i\leq (4\d |f|_{r,s})^{i+1}\,,
\qquad
|f_i^\flat-f^\flat_{i-1}|_i\le 2 \tetta_{i-1}\, |f_{i-1}|_{i-1}\ ,
\end{equation}
where
\begin{equation}\label{corcira}
\tetta_i:=\d |f_i^K|_i\qquad {\rm with}\qquad 
\d:=
\frac{2^5\, n \, K^3}{\a r s}
\ .
\end{equation}
}
Notice that, recalling \equ{gricia}, it is
\beq{cappuccino}
\tetta_* =2^4 \d |f|_{r,s} \quad \Longrightarrow\quad
4\d |f|_{r,s} = \frac{\tetta_*}{4} < \frac1{4}<1\ .
\eeq
%%%%% passo 1
Let us first show that  the inductive hypothesis is true for $j=1$. Indeed, by \equ{cappuccino}, 
$\d|f^K|_0\le \d |f|_0<1/16<1$, therefore, by the definition of $\d$ and $\check\tetta$ in, respectively,  \equ{corcira} and \equ{gricia*},   
we see that  we
can apply Lemma~\ref{megara} with $f=f_0$, being $\check\tetta=\tetta_0=\d|f^K|_0$. Thus, we obtain the existence of $\Psi_0$ so that 
$H_1:=H_0\circ \Psi_0=h+f_1$ and, by \equ{nuvole} and \equ{salamina},
$$
\tetta_1=\d |f_1^K|_1 \le \d (2  \tetta_0 |f|_0)=2\, \d^2 |f^K|_0 \, |f|_0\le 2 (\d |f|_0|)^2\le (4\d |f|_0)^2\ ,
$$
showing that the first inequality in \equ{pontina} holds for $i=1$, the second inequality follows from \equ{bemolle}.

\nl
Now, let us assume that the inductive hypothesis holds true for $1\le i\le j<K$ and let us prove that it holds also for $i=j+1$.
First, let us check that
\beq{ausoni}
|f_i|_i\le 2|f|_{r,s}\ ,\qquad \forall\ 1\le i\le j\ .
\eeq
Indeed, by the estimate in \equ{nuvole}, one has that $|f_i|_i\le (1+2\tetta_i)|f_{i-1}|_{i-1}$, for all $1\le 1\le j$, which, iterated, yields
\begin{eqnarray}
|f_i|_i
&\leq&\dst
 |f_0|_0 \prod_{\ell=1}^{i}(1+2 \tetta_\ell )
=|f|_{r,s}\exp\big(\sum_{\ell=1}^i\log(1+2 \tetta_\ell) \big)
\leq 
|f|_{r,s}\exp\big(2\sum_{\ell=1}^i \tetta_\ell \big)
\nonumber
\\
&\stackrel{\eqref{pontina}}\leq&\dst
|f|_{r,s}\exp\big({2\sum _{\ell=1}^i  (4\d|f|_{r,s})^\ell }\big)
\stackrel{\equ{cappuccino}}{\le}|f|_{r,s}\exp\big({2\sum _{\ell=1}^\infty  2^{-2\ell}\big)}
\le 2 |f|_{r,s}\,.
\nonumber
\end{eqnarray}
Now, by \equ{corcira}, \equ{pontina} with $i=j$ (inductive assumption) and \equ{cappuccino}, we have that $\tetta_j<1$. Thus, we can apply Lemma~\ref{megara} to $f_j$ (with $\check\tetta =\tetta_j$) and get a symplectic transformation $\Psi_j$ such that $H_{j+1}:=H_{j}\circ \Psi_{j}
=h+f_{j+1}$ satisfies 
\beqno
\tetta_{j+1}:= \delta |f^K_{j+1}|_{j+1}\stackrel{\equ{salamina}}\le  \d\, (2\tetta_j |f_j|_j)
\stackrel{\equ{ausoni}}\le (4\d |f|_{r,s}) \, \tetta_j
\stackrel{\equ{pontina}_j}\le 
(4\d |f|_{r,s}|)^{j+2}\ ,
\eeqno
which is the first inequality in \equ{pontina} with $i=j+1$, the second inequality comes from \equ{bemolle}. This completes the proof of the induction.

\nl
Now,  we can conclude the proof of Lemma \ref{pesce}: recall \equ{tetto} and define 
$$\Psi:=\Psi_0\circ\cdots\circ\Psi_{\bar K-1}\ .
$$
Notice that, by \equ{cornetto},  $r_{\bar K}=r/2=r_*$ and $s_{\bar K}=s(1-1/K)=s_*$ and notice that, by the induction, it is
\beq{bic}
H\circ \Psi=H_{{\bar K}-1}\circ\Psi_{{\bar K}-1}\stackrel{\equ{olintoj}_{\bar K}}=h+f_{\bar K}=:h+f^\flat+f_*\ .
\eeq
But, then, since $T_K P_\L^\perp f^\flat=(f^\flat)^K=0$ (for any $f$), using  that ${\bar K}\ge 2$, we have
\beqa{biro}
|T_K P_\L^\perp f_*|_{r_*,s_*} &=& |f_{\bar K}^K|_{\bar K}\stackrel{\equ{corcira}}= \d^{-1} \tetta_{\bar K}\stackrel{\equ{pontina}}\le \d^{-1} (4\d |f|_0)^{{\bar K}+1}
=4 (4\d |f|_0)^{\bar K} \, |f|_0\nonumber \\
&\le& (2^3\d |f|_0)^{\bar K} |f|_0 \stackrel{\equ{cappuccino}}= \big(2^{-1} \tetta_* \big)^{\bar K} |f|_{r,s}
\le
\big(2^{-1} \tetta_* \big)^K |f|_{r,s}
\ ,
\eeqa
proving the second estimates in \equ{pirati}. 

\nl
Finally, (using again that ${\bar K}\ge 2$ and that $\tetta_* <1$)
\beqano
|f_*|_{r_*,s_*}&\stackrel{\equ{bic}}=&|f_{\bar K}-f^\flat|_{\bar K}\stackrel{\equ{decomposizione}}= |f_{\bar K}^K+f_{\bar K}^\flat-f^\flat|_{\bar K}\le
|f_{\bar K}^K|_{\bar K}+|f_{\bar K}^\flat-f^\flat|_{\bar K}\\
&\stackrel{\equ{biro}}\le& \frac{\tetta_*}4\,  |f|_0  + \sum_{i=1}^{\bar K} |f_i^\flat-f_{i-1}^\flat|_i\stackrel{\equ{pontina}}\le
\frac{\tetta_*} 4\, |f|_0 + \sum_{i=1}^{\bar K} 2 \tetta_{i-1} |f_{i-1}|_{i-1}
\\
&\stackrel{\equ{ausoni},\equ{pontina}}\le& \frac{\tetta_*}4\,  |f|_0  +  4 |f|_0 \sum_{i=1}^{\bar K} (4\d |f|_0)^i
\stackrel{\equ{cappuccino}}=
\frac{\tetta_*}4 \, |f|_0  +   4 |f|_0 \sum_{i=1}^{\bar K} (\tetta_* /4)^i\le 2 \tetta_* \, |f|_0\ ,
\eeqano
which proves also the first estimate in \equ{pirati}.
\eproof

\section{On action--angle variables for 1D mechanical systems with parameters}\label{gerusalemme}

We will use the notations of sections
\ref{erezione} and  \ref{riso},  in particular
subsections \ref{farro} and  \ref{orzo}.

\subsection{The ``unperturbed case''}

Consider the ``unperturbed case'' when
$\checco=\checco_*=0$ (recall \eqref{tolomeo2}).
 Namely consider the one dimensional Hamiltonian
\begin{equation}\label{goffredo}
\Hpend^0(J_n,\psi_n)=J_n^2+F^0(\psi_n)\,,
\quad \text{with}  \ F^0 \ \text{satisfying }\ \eqref{legna}	\,.
\end{equation}
In the particular important case in which $F^0$
is minus cosine  we can explicitly evaluate
\begin{equation}\label{lontra}
F^0(x)=-\cos x\quad\Longrightarrow\quad
M=\cosh s_0\,,\ \
N=1\,,\ \ 
x^0_1=0\,,\ x^0_2=\pi\,,\ \
E^0_1=-1\,,\ E^0_2=1\,,\ \
\b=1\,.
\end{equation}
For $E\in (E^{(i),0}_-,E^{(i),0}_+)$, let us define the functions $P_n^{(i),0}(E)$ as
\begin{eqnarray}
P_n^{(2j-1),0}(E)&:=&
\frac{1}{\pi}
\int_{X_{2j-1}^0(E )}^{X_{2j}^0(E )}
\sqrt{E-F^0( x )}\, dx\,,
\nonumber
\\
P_n^{(2j),0}(E)&:=& 
\frac{1}{\pi}
\int_{X_{2j_-+1}^0(E )}^{X_{2j_+}^0(E )}
\sqrt{E-F^0( x )}\, dx\,,
\nonumber
\\
P_n^{(2N),0}(E)&:=& 
\frac{1}{2\pi}
\int_{-\pi}^{\pi}
\sqrt{E-F^0( x )} \, dx\,,
\nonumber
\\
P_n^{(0),0}(E)&:=&
-\frac{1}{2\pi}
\int_{-\pi}^{\pi}
\sqrt{E-F^0( x )} \, dx\,.
\label{sunday}
\end{eqnarray}
In the following we will use the notations
$\pgot$ and $\lessdot$
introduced in \eqref{pgot} and \eqref{pgot2}.

\begin{lemma}\label{maso}
 For real $E$, we have that
  \begin{equation}\label{vana0}
\min_{1\leq i\leq 2N-1 } \inf_{E\in (E^{(i),0}_-,E^{(i),0}_+)}
\partial_E P_n^{(i),0}(E)
=: C_{F^0}>0
\,.
\end{equation}
In particular\footnote{
 In the special case in which $F^0(\psi_n)=-\cos \psi_n$
(note that $N=1$), the minimum is $1/\sqrt 2.$ }
\begin{equation}\label{malandrino}
\text{if}\ \ F^0 \ \ \text{satisfies}\ \ 
\text{({\bf A3})}\ \ \text{then}\ \ 
C_{F^0}\geq 1/2\,.
\end{equation}
\end{lemma}
\proof See \cite{BCaa}.
\eproof

\

\Giu

%%%%%%%%%%%%%%%%%
%%%%%%%%%%%%%%%%%%

\subsection{The action as a function of the angle at constant energy}

Let us consider now the  Hamiltonian $\Hpend$
defined in
\eqref{guglielmo}.

For $\checco $ small enough we can solve, w.r.t. $J_n$,  the implicit function equation
\begin{equation}\label{meaux}
 J_n^*(\hat J )+\frac{z}{\sqrt{1+ b(J,\psi_n )}}-J_n=0\,,
\end{equation}
finding
$$
J_n=\mathcal J_n(z,\psi_n,\hat J )\,,
$$
where $\mathcal J_n$ is the analytic function $\mathcal J_n:(-R_0,R_0)_{r_0/4}\times\mathbb T_{s_0}\times \hat D_{r_0}$
\begin{equation}\label{spiderman}
\mathcal J_n(z,\psi_n,\hat J )=
J_n^*(\hat J )+
z+\tilde{\mathcal J_n}(z,\psi_n,\hat J )\,
\end{equation}
whit $\tilde{\mathcal J_n}$ solving 
the fixed point equation
\begin{equation}\label{meauxbis}
\tilde{\mathcal J_n}=
\Phi(\tilde{\mathcal J_n};z,\psi_n,\hat J):=
\frac{1}{\sqrt{1+ b\big(\hat J,
  J_n^*(\hat J )+z+\tilde{\mathcal J_n},\psi_n \big)}}-1\,.
\end{equation}
We are going to solve 
 \eqref{meauxbis} 
in the closed ball 
\begin{equation}\label{3holes}
\|\tilde{\mathcal J_n}\|_{\mathcal B}
\leq  \checco\,, 
\end{equation}
of the Banach space
$\mathcal B$
of analytic functions $\phi:(-R_0,R_0)_{r_0/4}\times\mathbb T_{s_0}\times \hat D_{r_0}\to \C$
endowed with the sup-norm
$$
\|\phi\|_{\mathcal B}:=
\sup_{z\in (-R_0,R_0)_{r_0/4}} 
\|\phi(z,	\cdot,\cdot)\|_{\hat D,r_0,s_0}\,.
$$
 We first note that
by \eqref{spiderman}, \eqref{3holes}, \eqref{mazinga}
we get
\begin{equation}\label{antigono}
\|\mathcal J_n(z,	\cdot,\cdot)\|_{\mathcal B}
\leq r_0 \checco+ R_0+ r_0/4+ \checco
\leq R_0+3r_0/8\,,
\end{equation}
assuming
\begin{equation}\label{suppersready}
\checco\leq  \min\{1,r_0\}/32\,.
\end{equation}
 For $|t|\leq 1/4$,
we have that $\left|\frac{d}{dt}\frac{1}{\sqrt{1+t}}\right|\leq 1,$ then
we get by \eqref{mazinga}, \eqref{suppersready},
\eqref{antigono} and Cauchy estimates 
\begin{equation}\label{itsasin}
\|\Phi(\tilde{\mathcal J_n})\|_{\mathcal B}\leq \checco
\,,\qquad
\|D_{\tilde{\mathcal J_n}}\Phi(\tilde{\mathcal J_n})
\|_{\mathcal L(\mathcal B,\mathcal B)}
\leq
\left\|
\frac{\partial_z b(\hat J,
  J_n^*+z+\tilde{\mathcal J_n},\psi_n )}{2\big(
  1+b(\hat J,
  J_n^*+z+\tilde{\mathcal J_n},\psi_n )\big)^{3/2}}
\right\|_{\mathcal B}
\leq 
 \frac{8 \checco}{r_0}
 \leq
 \frac{1}{4}
\end{equation}
in the closed ball of  $\tilde{\mathcal J_n}$ satisfying \eqref{3holes}.
\\
Obviously\footnote{For real values of $\hat J,$ $\psi_n,$ $E.$}
\begin{equation}\label{ummagamma}
J_n=\mathcal J_n\Big(\pm\sqrt{E-F(\hat J,\psi_n )},\psi_n,\hat J \Big)\quad
{\rm solves \ (w.r.t.} \ J_n) \quad
\Hpend(\hat J,J_n,\psi_n )=E\,,
\end{equation}
according to $\pm \big( J_n- J_n^*(\hat J ) \big)\geq 0,$
for every (real) $E$ such 
that\footnote{So that $\mathcal J_n\Big(\pm\sqrt{E-F(\hat J,\psi_n )},\psi_n,\hat J \Big)$
is well defined. Recall \eqref{legna}.}, 
\begin{equation}\label{adcazzum}
E+M< R_0^2\,.
\end{equation} 
By \eqref{meauxbis} we get
$$
\partial_z \tilde{\mathcal J}_n=-\left(
1+\frac{\partial_z b}{2(1+ b)^{3/2}}\right)^{-1}
\frac{\partial_z b}{2(1+ b)^{3/2}}
$$
so that, recalling \eqref{itsasin},
\begin{equation*}
\| \partial_z \tilde{\mathcal J}_n \|_{\mathcal B}\leq
\frac43 \frac{8\checco}{r_0}\leq \frac13\,.
\end{equation*}
Then $\mathcal J_n$ is an increasing function of (real) $z$,
indeed by \eqref{spiderman} we obtain
$\partial_z \mathcal J_n=1+\partial_z \tilde{\mathcal J}_n$.

\Giu

\begin{remark}\label{badedas}
In the following we will often omit the explicit dependence on $\hat J$, for brevity. 
\end{remark}

\subsection{The domains of definition of action angle variables}

Outside the zero measure set formed by the connected components in the set of  critical energies $\{ \Hpend=E_i \},$
$1\leq i\leq 2N,$ containing the critical points  $x_i$,
the phase space $\R^n\times\T^n$
is composed by $2N+1$ open connected components
$
\mathcal C^i,$ $0\leq i\leq 2N,$
defined as 
$$
\mathcal C^{i}:=\check{\mathcal C}^{i}\times\T^{n-1}\,,
$$
where
$$
\check{\mathcal C}^{i} 
\subseteq\hat D\times\R\times\T^1
\subseteq\R^n\times\T^1
$$
are defined as follows\footnote{Omitting to write, for brevity, 
the explicit dependence of $\mathcal J_n, F, X_i$ on
$\hat J$}.
\\
For $i=2j-1$ odd, $1\leq j\leq N,$
$\check{\mathcal C}^{2j-1} $ is a normal set with respect to the variable
$J_n,$
\begin{eqnarray}\label{skys}
&&\check{\mathcal C}_{2j-1}   :=
\\
&&\Big\{
\mathcal J_n
\Big(-\sqrt{E^{(2j-1)}_+(\hat J) -F( \psi_n   )},\psi_n \Big)
< J_n<
\mathcal J_n
\Big(\sqrt{E^{(2j-1)}_+(\hat J) -F( \psi_n   )},\psi_n \Big)\,,
\nonumber
\\
&&\qquad X_{2j-1}\big(E^{(2j-1)}_+(\hat J)   \big)
< \psi_n <
X_{2j}\big(E^{(2j-1)}_+(\hat J)   \big)\,,\quad 
\hat J\in\hat D\ 
\Big\}
\nonumber
\\ &&\setminus \ 
\Big\{  
J_n= J_n^* \,, \psi_n=
x_{2j-1} 
\Big\}\,.
\nonumber
\end{eqnarray}

For $i=2j$ even, $1\leq j\leq N-1,$
$\check{\mathcal C}^{2j}  $ is still a normal set with respect to the variable $J_n$:
 
 \begin{eqnarray}\label{skys2}
&&\check{\mathcal C}_{2j}   :=
\\
&&\Big\{
\mathcal J_n
\Big(-\sqrt{E^{(2j)}_+(\hat J) -F( \psi_n   )},\psi_n \Big)
< J_n<
\mathcal J_n
\Big(\sqrt{E^{(2j)}_+(\hat J) -F( \psi_n   )},\psi_n \Big)\,,
\nonumber
\\
&&\quad
X_{2j_-+1}\big(E^{(2j)}_+(\hat J)   \big)
< \psi_n <
X_{2j_+}\big(E^{(2j)}_+(\hat J)   \big)\,,
\quad \hat J\in\hat D\ 
\Big\}\  
\nonumber
\\
&&\setminus \ \Big\{
\mathcal J_n
\Big(-\sqrt{E^{(2j)}_-(\hat J) -F( \psi_n   )},\psi_n \Big)
\leq J_n\leq
\mathcal J_n
\Big(\sqrt{E^{(2j)}_-(\hat J) -F( \psi_n   )},\psi_n \Big)\,,
\nonumber
\\
&&\qquad X_{2j_-+1}\big(E^{(2j)}_-(\hat J)   \big)
\leq \psi_n \leq
X_{2j_+}\big(E^{(2j)}_-(\hat J)   \big)\,,
\quad \hat J\in\hat D\ 
\Big\}\,,
\nonumber
\end{eqnarray}
where   $j_-,j_+$ were defined in \eqref{21stcentury}.
\\
Finally
\begin{eqnarray}\label{skys3}
\check{\mathcal C}_{2N}   
&:=&
\Big\{
J_n>
\mathcal J_n
\Big(\sqrt{E^{(2N)}_-(\hat J) -F( \psi_n   )},\psi_n \Big)\,,\ \ \ \psi_n\in\T \,,
\quad \hat J\in\hat D\ 
\Big\}
\\
\label{skys4}
\check{\mathcal C}_{0}   
&:=&
\Big\{
J_n<
\mathcal J_n
\Big(-\sqrt{E^{(2N)}_-(\hat J) -F( \psi_n   )},\psi_n \Big)\,,\ \ \ \psi_n\in\T\,,
\quad \hat J\in\hat D\  
\Big\}
\end{eqnarray}
Note that actually in $\check{\mathcal C}_i  $ with $1\leq i<2N,$
$\psi_n$ is not an angle!

Let us introduce the (small) parameter 
\begin{equation}\label{moda}
\theta\geq 0\,.
\end{equation}
Recalling \eqref{adcazzum}, we define the following subsets of $\check{\mathcal C}^i  $ (defined in \eqref{skys},\eqref{skys2},\eqref{skys3}, \eqref{skys4})
\begin{eqnarray}
 \check{\mathcal C}^{2j-1}(\theta ) &:=&
  \check{\mathcal C}^{2j-1}  \cap \{ 
 E^{(2j-1)}_-(\hat J) <\Hpend<E^{(2j-1)}_+(\hat J)-2\theta   \}\,,\qquad
  \text{for}\ \ \  1\leq j\leq N\,,
  \nonumber
\\
 \check{\mathcal C}^{2j}(\theta ) &:=&
  \check{\mathcal C}^{2j}  \cap \{ E^{(2j)}_-(\hat J)+2\theta <\Hpend<E^{(2j)}_+(\hat J) -2\theta  \}\,,\qquad
  \text{for}\ \ \  1\leq j<N\,,
  \nonumber
\\
 \check{\mathcal C}^i(\theta ) &:=&
  \check{\mathcal C}^i  \cap \{ E^{(2N)}_-(\hat J)+2\theta <\Hpend 
  <R_0^2-M-2\theta  \}\,,
  \qquad
  \text{for}\ \ \  i=0,\, 2N\,,
  \label{frigia}
\end{eqnarray}
where
$\Hpend$ was defined in
\eqref{guglielmo}. Note that $\check{\mathcal C}^i(0 )=\check{\mathcal C}^i$ for $1\leq i<2N.$

Finally we define\footnote{With a little abuse of notation we invert the order of
$\psi_n$ and $\hat \psi.$} 
\begin{equation}\label{siberia}
{\mathcal C}^i(\theta)=\check{\mathcal C}^i(\theta)\times\T^{n-1}\,,\qquad 
\check{\mathcal C}^i(\theta)\ni(J,\psi_n)\,,
\quad
\T^{n-1}\ni\hat \psi\,.
\end{equation}
The sets $\check{\mathcal C}^i(\theta)$ and, therefore, ${\mathcal C}^i(\theta)$ have different homotopy: for every fixed $\hat J$
the set
$$
\check{\mathcal C}^i_{\hat J}(\theta):=\{
(J_n,\psi_n) \ |\ (J,\psi_n)\in \check{\mathcal C}^i(\theta)
\}\subseteq \R^1\times\T^1
$$ 
is contractible for $1\leq i\leq 2N-1$
and is not contractible for $i=0,2N.$
Note that, recalling \eqref{legna},
\begin{equation}\label{gonata}
\hat D\times (-R_0/2,R_0/2)\times\T^n
\subset
\bigcup_{0\leq i\leq 2N}\overline{\mathcal C^i(0)}
\subset
\hat D\times (-R_0,R_0)\times\T^n\,.
\end{equation}

\subsection{Definition of action variables}\label{hocuspocus}

On the above connected components $\mathcal C_i,$
$0\leq i\leq 2N,$
we want to define action angle variables integrating
$\Hpend.$
We first define the action variables 
as a function of the energy $E$ and of the dummy
variable $\hat J$.
More precisely, for $0\leq i\leq 2N,$ we are going to define the functions
$$
P_n^{(i)}\ :\ (E, \hat J )\in \mathcal E^{i}\ \to \R\,,
\qquad {\rm where}\ \ \ \mathcal E^{i}:=
\mathcal E^{i}(0)
$$
and
\begin{eqnarray}
\mathcal E^{2j-1}(\theta) &:=&
\{ (E,\hat J) \ \ \text{s.t.}\ \   E^{(2j-1)}_-(\hat J)
<E<E^{(2j-1)}_+(\hat J)-2\theta\,,\
\hat J\in\hat D\}\,,
\qquad 1\leq j\leq N\,,
\nonumber
\\
\mathcal E^{2j}(\theta) &:=&
\{ (E,\hat J) \ \ \text{s.t.}\ \   E^{(2j)}_-(\hat J)
+2\theta<E<E^{(2j)}_+(\hat J)-2\theta
\,,\
\hat J\in\hat D\}\,,
\qquad 1\leq j< N\,,
\nonumber
\\
\mathcal E^{2N}(\theta)=\mathcal E^0(\theta) &:=&
\{ (E,\hat J) \ \ \text{s.t.}\ \   E^{(2N)}_-(\hat J)+2\theta<E<R_0^2-M-2\theta\,,\
\hat J\in\hat D\}\,,
\label{vispo}
\end{eqnarray}
where  the positive parameter $\theta$ was introduced in \eqref{moda}
(and recall \eqref{adcazzum}).
We also introduce the complex $\theta$-neighborhoods\footnote{Recall the notation on page \pageref{giallo}.} 
\begin{equation}\label{vispa}
\mathcal E^i_\theta(\theta):=
\big(\mathcal E^i(\theta)\big)_\theta
\subseteq \C^n\,.
\end{equation}
\\
The functions $P_n^{(i)}$ are defined as 
follows\footnote{Sometimes omitting, for brevity,
to write the explicit dependence on $\hat J$.}.
\\
For $i=2j-1$ odd, $1\leq j\leq N,$ 
and $E^{(2j-1)}_-(\hat J) <E<E^{(2j-1)}_+(\hat J) $, we set
\begin{eqnarray}\label{musicalbox}
&&
P_n^{(2j-1)}(E)=P_n^{(2j-1)}(E, \hat J )
\\
\nonumber
&&
:=
\frac{1}{2\pi}
\int_{X_{2j-1}(E )}^{X_{2j}(E )}
\Big[
\mathcal J_n
\Big(\sqrt{E-F( x )},x \Big) 
-
\mathcal J_n
\Big(-\sqrt{E-F( x )},x \Big)
\Big]
\, dx
\\
\nonumber
&&
=\frac{1}{\pi}
\int_{X_{2j-1}(E )}^{X_{2j}(E )}
\sqrt{E-F( x )}\Big(
1
+ b_\sharp\big(\sqrt{E-F( x )},x \big)
\Big)\, dx\,,
\end{eqnarray}
where the last equality holds recalling \eqref{meaux}
and defining $b_\sharp=
b_\sharp(\hat J,z,x )$ as follows:
\begin{equation}\label{4holes}
2+2 b_\sharp(\hat J,z,x ):=
\frac{1}{\sqrt{1+ b\big(\hat J,\mathcal J_n(z,x,\hat J ), x \big)}}
+
\frac{1}{\sqrt{1+ b\big(\hat J,\mathcal J_n(-z,x,\hat J ), x \big)}}
\end{equation}
($b$ defined in \eqref{guglielmo}).
Note that 
\begin{equation}\label{3holesbis}
b_\sharp \ \ \text{is even w.r.t.} \ z \ \ \text{and}\ \ 
\sup_{z\in (-R_0,R_0)_{r_0/2}} \|b_\sharp(\hat J,z,x)\|_{\hat D,r_0,s_0}\lessdot \checco\,.
\end{equation}
For $i=2j$ even, $1\leq j\leq N-1$
and $E^{(2j)}_-(\hat J) <E<E^{(2j)}_+(\hat J) $, 
we set (recall \eqref{21stcentury})
\begin{eqnarray}\label{musicalbox2}
&&
P_n^{(2j)}(E)=
P_n^{(2j)}(E, \hat J )
\\
\nonumber
&&
:=
\frac{1}{2\pi}
\int_{X_{2j_-+1}(E )}^{X_{2j_+}(E )}
\Big[
\mathcal J_n
\Big(\sqrt{E-F( x )},x \Big) 
-
\mathcal J_n
\Big(-\sqrt{E-F( x )},x \Big)
\Big]
\, dx
\\
\nonumber
&&
=\frac{1}{\pi}
\int_{X_{2j_-+1}(E )}^{X_{2j_+}(E )}
\sqrt{E-F( x )}\Big(
1
+ b_\sharp\big(\sqrt{E-F( x )},x \big)
\Big)\, dx\,,
\end{eqnarray}
where   $j_-,j_+$ were defined in \eqref{21stcentury}.
\\
Finally for $E>E^{(2N)}_-(\hat J) $ we set
\begin{eqnarray}\label{musicalbox3}
P_n^{(2N)}(E)=
P_n^{(2N)}(E, \hat J )
&:=&
\frac{1}{2\pi}
\int_{-\pi}^{\pi}
\mathcal J_n
\Big(\sqrt{E-F( x )},x \Big) \, dx\,,
\\
\label{musicalbox4}
P_n^{(0)}(E)=
P_n^{(0)}(E, \hat J )
&:=&
\frac{1}{2\pi}
\int_{-\pi}^{\pi}
\mathcal J_n
\Big(-\sqrt{E-F( x )},x \Big) \, dx\,.
\end{eqnarray}
%%%%%%%%%%%%%%%%%%%%

\subsection{Properties of the actions as functions of the energy and viceversa}

\begin{lemma}\label{raggio}
 Assume that
\begin{equation}\label{parigibis}
\frac{\checco}{c_{F^0}}\leq \checco_0(\pgot)\,,
\end{equation}
with $c_{F^0}$ defined in \eqref{vana0} and $\checco_0=\checco_0(\pgot)$ small enough.
Then  
for every $ 1\leq i\leq 2N-1$ 
\begin{equation}\label{vana}
\inf \partial_E P_n^{(i)}(E,\hat J)  \geq\  c_{F^0}/2>0\,,
\end{equation}
while\footnote{Recall \eqref{adcazzum}} 
\begin{eqnarray}\label{moldavater}
  \frac{1}{8\sqrt E}\leq \partial_E P_n^{(2N)}(E,\hat J)\,,\ -\partial_E P_n^{(0)}(E,\hat J)
  \leq  \frac{2}{\sqrt E}\,, \qquad \forall\,  2M\leq E \leq R_0^2-M\,,\ \ \hat J\in\hat D\,,
\\
\label{moldava2ter}
\partial_E P_n^{(2N)}(E,\hat J)\,,\ -\partial_E P_n^{(0)}(E,\hat J)
\geq  \frac{1}{8\sqrt{2M}}\,,\qquad
\forall\, E_{2N}<E\leq 2M\ \ \hat J\in\hat D
\end{eqnarray}
($M$ defined in \eqref{goffredo}).
\end{lemma}
\proof It essentially follows from
Proposition
\ref{glicemiak};  see \cite{BCaa}
for details. 
\eproof

By \eqref{vana} we have that, for every fixed $\hat J\in\hat D,$
the function $E\mapsto P_n^{(i)}(E,\hat J)$ is strictly monotone
and, therefore, invertible with inverse
$\mathtt E^{(i)}(\hat J, P_n)$ such that
\begin{equation}\label{13ottavi}
\mathtt E^{(i)}(\hat J, P_n^{(i)}(E,\hat J))=E\qquad
\text{and}\qquad
P_n^{(i)}(\mathtt E^{(i)}(\hat J, P_n),\hat J))=P_n\,.
\end{equation}
As a corollary of Proposition \ref{glicemiak}
and of the cain rule applied to \eqref{13ottavi}, giving
$$
\partial_{P_n} \mathtt E^{(i)}
(P)  
=
\frac{1}{\partial_E P_n }\,,\qquad
\partial_{\hat P} \mathtt E  
=
-\frac{\partial_{\hat P} P_n }{\partial_E P_n }\,,
$$
by Lemma \ref{raggio}
we get the following.
For every $ 1\leq i\leq 2N-1$
 \begin{equation}\label{vana00}
0<\partial_{P_n} \mathtt E^{(i)} 
  \leq\  2/c_{F^0}\,,
\end{equation}
($c_{F^0}$ defined
in \eqref{vana0}). 
Moreover
\begin{eqnarray}\label{moldavabis}
  \frac{\sqrt E}{2}\leq 
  \partial_{P_n} \mathtt E^{(2N)} \,,\ -\partial_{P_n} \mathtt E^{(0)}
  \leq 8\sqrt E\,, \qquad \forall\, 2M\leq E \leq R_0^2-M\,,\ \ \hat P\in\hat D\,,
\\
  \partial_{P_n} \mathtt E^{(2N)} \,,\ -\partial_{P_n} \mathtt E^{(0)}
\leq 8\sqrt{2M}\,,\qquad
\forall\, E^{(0)}_-\,, \ E^{(2N)}_-<E\leq 2M\ \ \hat P\in\hat D
\label{moldava3}
\end{eqnarray}
($M$ defined in \eqref{goffredo}).

The proofs of the following lemmata
essentially follows from
Proposition
\ref{glicemiak};  see \cite{BCaa}
for details.

\begin{lemma}\label{sorellenurzia2}
Let $C>1$ as in Proposition \ref{glicemiak}.
\begin{equation}\label{ofena}
\sup_{\hat D_{r_0} \times 
\mathtt I^{(i)}_{\mathtt r_0/8}}
\|\partial_{PP}\mathtt E^{(i)} (P)\|\leq C\,,\qquad
\sup_{\hat D_{r_0} \times 
\mathtt I^{(i)}_{\mathtt r_0/8}}
|\partial_{P_n}\mathtt E^{(i)} (P)-
\partial_{P_n}\mathtt E^{(i),0} (P_n)|
\leq C\checco\,,
\qquad \forall
0\leq i\leq 2 N 
\,,
\end{equation}
where the intervals   $\mathtt I^{(i)}$ where defined in
 \eqref{ciclette} and $\mathtt r_0>0$
in \eqref{peppapig}.
\end{lemma}

\begin{lemma}\label{sorellenurzia}
 \eqref{vallinsu} holds.
 \end{lemma}

\begin{lemma}\label{lapajata}
 Assume that $F$ 
is cosine-like according to Definition
\ref{sorellastre}, with $\cgot$ (namely $\cgot_*$ defined in \eqref{drago}) small enough.
Then 
\begin{equation}\label{thor}
\inf_{E_1<E<E_2,\,\hat P\in \hat D}
\left|\partial_{P_n P_n} \mathtt E^{(1)}\right|
\, \geq \,
c_\sharp\,,
\qquad
\inf_{E_2<E<R_0^2-2,\,\hat P\in \hat D} 
\left|\partial_{P_n P_n} \mathtt E^{(2)}\right|
\, \geq \,
c_\sharp\,,
\end{equation}
for a suitable (absolute constant)
$c_\sharp>0$.
\end{lemma}

%%%%%%%%%%%%%%

\section{Miscellanea}\label{noccetti}
%\section{Reminders and Technicalities}\label{noccetti}

\subsubsection*{Composition of  maps} 
\begin{lemma}\label{cippa}
Let $y=L y',$ where $L:\R^n\to\R^n$ is  linear.
Then for every function $h:\R^n\to\R$ we have
$$
\det\Big(\partial_{y'y'}\big( h \circ L(y')\big)  \Big)
=(\det L)^2 \det(\partial_{yy}h)_{\vert{y=L y'}}\ .
$$
\end{lemma}
Lemma~\ref{cippa} follows immediately observing that  $\partial_{y'y'}\big( h(Ly')\big)=
L^T (\partial_{yy}h)_{\vert{y=L y'}}L$.

\nl
Given a Hamiltonian $H(J,\psi)$ 
we denote by $\Phi^t_H$ its flow at time $t.$
\begin{lemma}\label{hans} {\rm (i) [Time rescaling]} Let 
 $c>0$. Then
$\Phi^t_H=\Phi^{ct}_{c^{-1}H}$.\\
{\rm (ii) [Action rescaling]} Consider the conformally symplectic change of 
variables 
\beq{macchiaioli}
(J,\psi)=\Phi(\tilde J,\tilde \psi):=(c\tilde J,\tilde\psi)
\eeq
and set 
$\tilde H:=H\circ\Phi$ and  $\hat H:=c^{-1} H\circ\Phi$.
Then  $\Phi\circ \Phi^{t/c}_{\tilde H}=\Phi^t_H\circ\Phi$ and 
$\Phi\circ \Phi^t_{\hat H}=\Phi^t_H\circ\Phi$.
\end{lemma} 
The proof is a straightforward check.
%\subsubsection*{On the norm of the composition of maps}

\begin{lemma}\label{fragile}
 Consider $s_1>0, s>s_2>0$ and an holomorphic map
 $\Phi:\T^n_{s_1}\to\T^n_{s_2}$ and an holomorphic function
 $f$ with $|f|_s<\infty.$ Then
 $$
 |f\circ\Phi|_{s_1}\leq \coth^n\left(\frac{s-s_2}{2}\right)|f|_s
 \leq 
 \left(1+\frac{2}{s-s_2}\right)^n|f|_s
 \,.
 $$
\end{lemma}
\proof
By \eqref{battiato2} we get
$\dst
 |f\circ\Phi|_{s_1}\leq 
 \|f\circ\Phi\|_{s_1}
 \leq \|f\|_{s_2}
 \leq\coth^n\left(\frac{s-s_2}{2}\right)|f|_s\,.$ \qed

\begin{lemma}\label{pergamena}
Given a matrix $M\in {\rm Mat}_{n\times n}(\Z)$ with
$\det M=\pm 1,$ consider the symplectic linear map
$\Phi:\R^n\times \T^n\to\R^n\times \T^n$ defined as
$$
(J',\psi')=\Phi(J,\psi):=(M^T J, M^{-1} \psi)\,.
$$
Let $D\subseteq\R^n$ and $r,s>0.$ Set
$D':=(M^T)^{-1} D.$ Then 
$$
\Phi (D'_{r'}\times \T^n_{s'})\subseteq
D_r\times \T^n_s\,,\qquad
\text{where}\ \ \ 
r':=r/\|M^T\|\,,\ \ 
s':=s/\sup_{1\leq i\leq n}\sum_{1\leq j\leq n}|(M^{-1})_{ij}|\,.
$$ 
Moreover, given a function $f:D_{r_0}\times\T^n_{s_0}$
with $r_0\geq r,$ $s_0>s,$ we have
$$
|f\circ\Phi|_{D',r',s'}\leq 
\coth^n\left(\frac{s_0-s}{2}\right)|f|_{s_0}
 \leq 
 \left(1+\frac{2}{s_0-s}\right)^n|f|_{s_0}
 \,.
$$
\end{lemma}
The first part is obvious; the second part follows from Lemma~\ref{pergamena}.

\subsubsection*{Restrictions of maps}

\begin{lemma}\label{schiena}
 Let $\Phi:D'\times\T^n
 \to D\times\T^n$ be a real analytic map with holomorphic
 extension
 $$
 \Phi:D'_{r'}\times\T^n_{s'}
 \to D_r\times\T^n_s
 $$
 for some $r,r's,s'>0.$
 There exists
  a suitably small
constant  $c$ depending only on $n$ such that
\begin{equation}\label{ABBA}
\Phi \left(D'_{c a r'}\times
\T^n_{c a s'}\right)\ \ \subseteq \ D_{a r}\times\T^n_{a s}\,,
\qquad \forall\, 0<a\leq 1\,.
\end{equation}
\end{lemma}
\proof
  By 
Cauchy estimates,
applied to the various components of $\Phi$. 
\qed

\subsubsection*{\bf A group of parameter--dependent symplectic transformations} 
Let us consider the group $\cal G$
introduced in 
\eqref{grillo}.
 
\begin{lemma}\label{chepalle}
 Given a symplectic transformation of the form
 \eqref{grillo}, we have that, for every fixed $\hat J,$ the restriction
 $$
 (J_n,\psi_n)\mapsto
 \big(I_n(J,\psi_n), \f_n(J,\psi_n)\big) $$
 is also symplectic.
\end{lemma}
\proof
Note that by the conservation of the symplectic form
$d I\wedge d\f=dJ\wedge d\psi$ follows that
$\partial_{J_n} I_n\partial_{\psi_n}\f_n -
\partial_{J_n} \f_n
\partial_{\psi_n}I_n =1.$
\eproof
Recall the definition given in \equ{islands} and note that
\begin{equation}\label{peggylee}
{(\Phi_1\circ\Phi_2)}\check{\phantom A}=\check\Phi_1\circ\check\Phi_2\,.
\end{equation}
Furthermore, obviously, one has
\beq{Ventura}
\check\phi(E)\times \torus^{n-1}= \phi(E\times \torus^{n-1})\ ,\qquad \forall \phi\in{\cal G}\ ,\qquad \forall E\subseteq \real^n\times \torus^1\ .
\eeq
By Lemma \ref{chepalle} we have the following
\begin{lemma}\label{chepalle2}
If  $\Phi\in {\cal G}$, then $\check \Phi$ is volume-preserving.
 \end{lemma}

\subsubsection*{An elementary result in linear algebra}\label{appendicedialfonso}
\begin{lemma}\label{alfonso}
Given $k\in\Z^n,$ $k\neq 0$
there exists a matrix $A=(A_{ij})_{1\leq i,j\leq n}$ with integer entries such that
$A_{nj}=k_j$ $\,\forall\, 1\leq j\leq n$,
$\det A=d:={\rm gcd}(k_1,\ldots,k_n)$, and 
$|A|_\infty= |k|_\infty$.
%\begin{equation}\label{atlantide}
%|A_k^{-1}|_\infty\leq (n-1)^{(n-1)/2}|k|_\infty^{n-1}\,.
%\end{equation}
\end{lemma}

\proof
The argument is by induction over $n$.
For $n=1$ the lemma is obviously true. For $n=2$, it follows at once from\footnote{The first statement in this formulation of Bezout's Lemma is well known 
and it can be found in any textbook on elementary number theory; the estimates on $x$ and $y$ are easily deduced from the well known fact that 
given a solution $x_0$ and $y_0$ of the equation $a x + b y =d$, all other solutions have the form $x= x_0+k (b/d)$ and $y=y_0- k (a/d)$ with $k\in \integer$
and by choosing $k$ so as to minimize $|x|$.
} 

\nl
{\bf Bezout's Lemma} {\sl Given two integers $a$ and $b$ not both zero, there exist two integers $x$ and $y$ such that $a x+ b y = d:={\rm gcd}(a,b)$, and such that 
$\max\{|x|,|y|\}\le \max\{|a|/d, |b|/d\}$.}

\nl
Indeed, if $x$ and $y$ are as in Bezout's Lemma with $a= k_1$ and $b=k_2$ one can take $A= \begin{pmatrix}y & -x\\ k_1 & k_2 \end{pmatrix}$.
Now, assume, by induction for $n\ge 3$ that the claim holds true for $(n-1)$ and let us prove it for $n$.   
Let $\bar k = (k_1,...,k_{n-1})$ and $\bar d={\rm gcd}(k_1,...,k_{n-1})$ and notice that ${\rm gcd}(\bar d, k_n)=d$. By the inductive assumption, there exists a matrix $\bar A=\begin{pmatrix} \tilde A\\ \bar k \end{pmatrix}\in
 {\rm Mat}_{(n-1)\times (n-1)}(\Z)$  with  $\tilde  A\in {\rm Mat}_{(n-2)\times (n-1)}(\Z)$, such that $\det \bar A=\bar d$ and $|\bar A|_\infty=|\bar k|_\infty$.
Now,  let $x$ and $y$ be as in Bezout's Lemma with $a= \bar d$, and $b=k_n$.
We claim that $A$ can be defined as follows:
\beq{giaccherini}
A=\begin{pmatrix} & \tilde k & \ & \tilde x \\ & \bar A &\  & \begin{pmatrix} 0 \\ \vdots \\ 0 \\k_n\end{pmatrix} \end{pmatrix}\ ,\qquad \tilde k= (-1)^n y\,  \frac{\bar k}{\bar d}\ ,\qquad \tilde x:= (-1)^{n+1} x \ .
\eeq
First, observe that since $\bar d$ divides $k_j$ for $j\le (n-1)$, $\tilde k \in \integer^{n-1}$. Then, expanding the determinant of $A$ from last column, we get
\beqano
\det A& =& (-1)^{n+1}\tilde x \det \bar A + k_n \det \begin{pmatrix} \tilde k \\ \tilde A\end{pmatrix}\\ 
& =& 
(-1)^{n+1}\tilde x \, \bar d + k_n (-1)^{n-2} \det  \begin{pmatrix} \tilde A \\ \tilde k\end{pmatrix} \\
&=& (-1)^{n+1}\tilde x\, \bar d + k_n  (-1)^{n-2}  (-1)^n \frac{y}{\bar d} \det \bar A\\
&=& 
x \bar d + k_n y = d\ .
\eeqano
Finally, by Bezout's Lemma, we have that $\max\{|x|,|y|\}\le \max\{\bar d/d , |k_n|/d\}$, so that
$$
|\tilde k|_\infty = |y| \frac{|\bar k|_\infty }{\bar d} \le \frac{|\bar k|_\infty}{d}\le |k|_\infty\ ,\quad
|\tilde x|= |x|\le \frac{|k_n|}{d}\le |k|_\infty\ ,
$$
which, together with $|\bar A|_\infty=|\bar k|_\infty$,  shows that $|A|_\infty=|k|_\infty$. \qed

\subsubsection*{Measure of sub--levels of smooth non--degenerate functions}\label{doveefinitocoglitore}
Here we prove Lemma~\ref{girotondi}. 
We start by recalling an elementary result, whose proof can be found in \cite{BCo}:
\begin{lemma}\label{girotondibis}
Let $g(x)$ a monic polynomial of degree $d$.
Then
$$
{\rm meas} \big( \{ x\in \mathbb R \ \ :\ \ |g(x)|\leq \gamma  \} \big)
\leq 2 d \gamma^{1/d}\,.
$$
\end{lemma}
We now prove Lemma \ref{girotondi}.
Let us divide the interval $[a,b]$ in  disjoint intervals
of length $2r:=2 \mu  ^{1/m+1}$.
Let $I$ one of such intervals and let $x_0$ is middle point.
By \eqref{macine} let $1\leq d\leq m$ such that
\begin{equation}\label{macinebis}
|\partial_x^{d}f(x_0)|/d!\geq \xi_m\,.
\end{equation}
By the Taylor remainder formula
we get, for $x\in I$,
$$
|f(x)-P^d_{x_0}(x)|\leq M r^{d+1}=M \mu  ^{\frac{d+1}{m+1}}\,,
$$
where $P^d_{x_0}(x):=\sum_{0\leq j\leq d}\frac{\partial_x^{j}f(x_0)}{j!}(x-x_0)^j$
is the Taylor polynomial of degree $d$. 
Then we have that
\begin{equation}\label{macineter}
\{x\in I \ \ :\ \ |f(x)|\leq \mu  \}\ 
\subseteq
\ 
\{x\in I \ \ :\ \ |P^d_{x_0}(x)|\leq (M+1)\mu  ^{\frac{d+1}{m+1}} \}
\,.
\end{equation}
We now apply Lemma~\ref{girotondibis} to the monic polynomial 
$g(x):=d! P^d_{x_0}(x)/\partial_x^{d}f(x_0)$
with 
$$
\gamma:=(M+1)\mu  ^{\frac{d+1}{m+1}}/\xi_m
\stackrel{\eqref{macinebis}}\geq
 d!(M+1)\mu  ^{\frac{d+1}{m+1}}/|\partial_x^{d}f(x_0)|\ .
$$
By \eqref{macineter}  and Lemma \ref{girotondibis}  we get
(recall $1\leq d\leq m$)
$$
{\rm meas}\big( \{x\in I \ \ :\ \ |f(x)|\leq \mu  \} \big)
\leq
2 d \gamma^{1/d}
\leq \frac{2m(M+1)}{\xi_m} \mu  ^{1/m}
$$
Since the number of disjoint intervals is smaller that
$\frac{b-a}{2\mu  ^{1/m+1}}+1,$
this concludes the proof of Lemma \ref{girotondi}.
\eproof

\subsubsection*{Canonical form of generalized pendula}
The following lemma describes how to make independent of the action $J_n$ a pendulum depending on parameters.

\begin{lemma}\label{francesco}
Let
$$
H^*(y,x):=y_n^2 + F^0(x_n) +  G^*(y,x_n)\,,
$$
with $\|G^*\|_{D,r_0,s_0}\leq \checco_*$ .
Assume that  
\begin{equation}\label{urea}
\checco_*\leq r_0^2/16\,.
\end{equation}
Then the fixed point equation
\begin{equation}\label{olimpia}
\mathtt y(\hat Y,X_n)
=
-\frac12 \partial_{Y_n} G^*\big( \hat Y,\mathtt y(\hat Y,X_n),X_n   \big)
\end{equation}
has 
a unique solution
 $\mathtt y=\mathtt y(\hat Y,X_n)$ with
\begin{equation}\label{urina}
\|\mathtt y\|_{\hat D,r_0,s_0}\leq 2\checco_*/r_0\leq r_0/8\,.
\end{equation}
 Set
\begin{equation*}
J^*_n(\hat Y)
:= \langle \mathtt y(\hat Y,X_n) \rangle\,,\qquad 
 a_*(\hat Y,X_n ):=\mathtt y(\hat Y,X_n)-\langle \mathtt y(\hat Y,X_n) \rangle\,,
\end{equation*}
where $ \langle \cdot \rangle$ denotes the average w.r.t. $X_n.$ Let $\phi=\phi(\hat Y,X_n)$ the unique function
satisfying $a_*=\partial_{x_n} \phi$
with $\langle \phi\rangle=0.$
Consider the canonical transformation  $\Psi$
\begin{equation}\label{rondo}
y_n=
Y_n+a_*(\hat Y,X_n)
=Y_n-
J^*_n(\hat Y)
+ \mathtt y(\hat Y,X_n)
\,,\qquad x_n=X_n\,,
\qquad \hat y =\hat Y\,,\quad 
\hat x=\hat X +b_*(\hat Y,X_n),
\end{equation}
obtained by the
 generating function
$Y_n x_n+\hat Y \hat x+ \phi(\hat Y,x_n),$
with
$b_*=-\partial_{\hat Y} \phi
$
and
\begin{equation}\label{urina3}
\|J^*_n\|_{\hat D,r_0}\leq 2\checco_*/r_0
\stackrel{\eqref{urea}}\leq r_0/8\,,\quad 
\|a_*\|_{\hat D,r_0,s_0}\leq 4\checco_*/r_0\,,
\quad 
\|b_*\|_{\hat D,r_0/2,s_0}
\leq (16\pi+8)\checco_*/r_0^2\,.
\end{equation}
Note that
\begin{equation}\label{flaviano}
\Psi\ :\ D_{r_0/2}\times
 \T^{n-1}_{\hat s}\times\T_{s_0}
 \ \to\ 
 D_{r_0}\times
 \T^{n-1}_{\hat s+(16\pi+8)\checco_*/r_0^2}\times\T_{s_0}\,.
\end{equation}
Then
\eqref{rondo} casts $H^*$ into
\begin{eqnarray*}
&&\Big(
Y_n-
J^*_n(\hat Y)
+ \mathtt y(\hat Y,X_n)
\Big)^2+
F^0(X_n) +  G^*(\hat Y,Y_n-
J^*_n(\hat Y)
+ \mathtt y(\hat Y,X_n) ,X_n)
\\
&&=
\big(1+ b (Y,X_n)\big) \big(Y_n-J^*_n(\hat Y)\big)^2  
  + F (\hat Y,X_n ) \,,
\end{eqnarray*}
 with
$$
F=F^0+G\,,\qquad 
G:= G^*(\hat Y, \mathtt y(\hat Y,X_n),X_n )
+ (\mathtt y(\hat Y,X_n))^2
$$
and\footnote{Using \eqref{olimpia}.}
omitting, for brevity, the dependence on
$\hat Y, X_n,$
\begin{equation}\label{uforobot}
b=\frac{
G^*(\mathtt y+Y_n- J_n^*)
-
G^*(\mathtt y)
-
\partial_{Y_n}G^*(\mathtt y) ( Y_n- J_n^*)
}
{( Y_n- J_n^*)^2}
=
\int_0^1 (1-t) \partial_{Y_n Y_n}G^*
\big(\mathtt y+t( Y_n- J_n^*) \big)dt
\,.
\end{equation}
Finally
\begin{eqnarray}\label{straussbis}
&&\|G\|_{\hat D,r_0,s_0}\leq \left( 1+4/r_0^2\right)\checco_*\,,
\qquad
\| (1+|Y_n- J_n^*|) b(Y,X_n)\|_{D,r_0/2,s_0}\leq
\left(4+\frac{34}{r_0^2}\right)\checco_*
\,,
\nonumber
\\
&&
\| |Y_n- J_n^*| \partial_{Y_n}b(Y,X_n)\|_{D,r_0/2,s_0}\leq
\frac{48}{r_0^2}\checco_*
\,.
\end{eqnarray}

\end{lemma}
\proof
\eqref{olimpia} is solved by the standard Fixed Point Theorem for $\mathtt y$ in the ball in \eqref{urina}.
\eqref{urina3} and the first estimate in \eqref{straussbis} follow by 
\eqref{urea}, \eqref{urina} and Cauchy estimates.
By \eqref{uforobot}, 
\eqref{urea}, \eqref{urina} and Cauchy estimates
we get
\begin{equation}\label{ilcapo}
\|b\|_{D,r_0/2,s_0}\leq
\frac{16}{r_0^2}\checco_*
\end{equation}
and, therefore,
dividing the cases  
$|Y_n- J_n^*|\leq 1$ and $|Y_n- J_n^*|>1,$
we get the second estimate in 
\eqref{straussbis}.
Finally, by \eqref{uforobot}, we have
\begin{equation*}
\partial_{Y_n} b=\frac{
\partial_{Y_n}G^*(\mathtt y+Y_n- J_n^*)
-
\partial_{Y_n}G^*(\mathtt y)
}
{( Y_n- J_n^*)^2}
- \frac{2b}{ Y_n- J_n^*}
\,. 
\end{equation*}
Then
$$
|\partial_{Y_n} b| |Y_n- J_n^*|
\leq
\frac{|
\partial_{Y_n}G^*(\mathtt y+Y_n- J_n^*)
-
\partial_{Y_n}G^*(\mathtt y)|
}
{|Y_n- J_n^*|} + 2|b|
\leq \frac{48}{r_0^2}\checco_*\,. \qedeq
$$

\end{document}